\documentclass[12pt]{amsart}
\usepackage{amsbsy}
\usepackage{graphicx,epsfig,subfigure}
\begin{document}

%

%
\newtheorem{theorem}{Theorem}
\newtheorem{proposition}[theorem]{Proposition}
\newtheorem{lemma}[theorem]{Lemma}
\newtheorem{corollary}[theorem]{Corollary}
\newtheorem{definition}[theorem]{Definition}
\newtheorem{remark}[theorem]{Remark}
\numberwithin{equation}{section} \numberwithin{theorem}{section}
\newcommand{\be}{\begin{equation}}
\newcommand{\ee}{\end{equation}}
\newcommand{\re}{{\mathbb R}}
\newcommand{\n}{\nabla}
\newcommand{\ren}{{\mathbb R}^N}
\newcommand{\iy}{\infty}
\newcommand{\pa}{\partial}
\newcommand{\ms}{\medskip\vskip-.1cm}
\newcommand{\mpb}{\medskip}
\newcommand{\ssk}{\smallskip}
\newcommand{\BB}{{\bf B}}
\newcommand{\Am}{{\bf A}_{2m}}
\newcommand{\bL}{\BB^*}
\newcommand{\bLs}{\BB}
\renewcommand{\a}{\alpha}
\renewcommand{\b}{\beta}
\newcommand{\g}{\gamma}
\newcommand{\ka}{\kappa}
\newcommand{\G}{\Gamma}
\renewcommand{\d}{\delta}
\newcommand{\D}{\Delta}
\newcommand{\e}{\varepsilon}
\newcommand{\vp}{\varphi}
\renewcommand{\l}{\lambda}
\renewcommand{\o}{\omega}
\renewcommand{\O}{\Omega}
\newcommand{\s}{\sigma}
\renewcommand{\t}{\tau}
\renewcommand{\th}{\theta}
\newcommand{\z}{\zeta}
\newcommand{\wx}{\widetilde x}
\newcommand{\wt}{\widetilde t}
\newcommand{\noi}{\noindent}
\newcommand{\lb}{\left (}
\newcommand{\rb}{\right )}
\newcommand{\lsb}{\left [}
\newcommand{\rsb}{\right ]}
\newcommand{\lab}{\left \langle}
\newcommand{\rab}{\right \rangle }
\newcommand{\gap}{\vskip .5cm}
\newcommand{\bz}{\bar{z}}
\newcommand{\bg}{\bar{g}}
\newcommand{\Ba}{\bar{a}}
\newcommand{\bt}{\bar{\th}}
\def\com#1{\fbox{\parbox{6in}{\texttt{#1}}}}


\title{\bf Regional, single point,\\ and global  blow-up for  the
fourth-order\\
 porous medium type equation with source}

\author
{V.A. Galaktionov} 

\address{ 
Department of Math. Sci., University of Bath,
 Bath, BA2 7AY, UK}
\email{vag@maths.bath.ac.uk}


 \keywords{Higher-order quasilinear
parabolic equation, finite propagation, blow-up, similarity
solutions, variational operators, branching.\\
 {\bf  Submitted
 to:} JPDE, September, 2008}
\date{\today}



\begin{abstract}

Blow-up
 behaviour  for the fourth-order quasilinear
porous medium equation
 with source,
 \be
 \label{00}
 u_t = - (|u|^n u)_{xxxx} + |u|^{p-1}u  \quad \mbox{in} \,\,\, \re \times \re_+,
 \quad n>0, \quad p>1,
  \ee
 is studied.
 Countable and  finite families  of  similarity
blow-up patterns of the form
 $$
  \mbox{$
 u_S(x,t)=(T-t)^{- \frac 1{p-1}}f(y), \quad  \mbox{where} \quad  y=x/(T-t)^{\b},
  \,\,\, \b = \frac {p-(n+1)}{4(p-1)},
   $}
 $$
 which blow-up as $t \to T^- < \infty$, are described.
 These solutions explain key  features of regional (for
$p=n+1$), single point ($p>n+1$), and global ($p \in(1,n+1)$)
blow-up. The concepts and various variational, bifurcation,  and
numerical approaches for revealing the structure and
multiplicities of such blow-up patterns are presented.




\end{abstract}


\maketitle


\section{Introduction: blow-up  reaction-diffusion models}


\subsection{On classic second-order blow-up models and higher-order diffusion}

Blow-up phenomena as intermediate asymptotics and
 approximations
of highly nonstationary processes  are common and well known in
various areas of mechanics and physics.
 The origin of intensive
  systematic
studies of such
nonlinear effects was gas dynamics (since the end of the 1930s and
1940s) supported later in the 1960s by plasma physics (wave
collapse) and nonlinear optics (self-focusing phenomena).
Nevertheless, it was {\em reaction-diffusion theory}
 that exerted
 the strongest influence on mathematical blow-up
research since the 1970s. It is not an exaggeration to say that
reaction-diffusion theory proposed basic and nowadays canonical
models, which eventually led to  qualitative and rigorous
description of principles of formation of blow-up and other
singularities in nonlinear PDEs.

 Finite-time blow-up  singularities
lie in the heart of several principal problems of PDE theory
concerning existence, uniqueness, optimal regularity, and
free-boundary propagation. The role of blow-up analysis in
nonlinear PDE theory will increase when more complicated classes
of higher-order degenerate parabolic, hyperbolic, nonlinear
dispersion, and other equations of interest are involved in the
 framework of massive mathematical research and application.
 For such general classes of equations with typically non-potential,
 non-divergent, and non-monotone operators (see classic monographs by Berger \cite{Berger} and
  Krasnosel'skii--Zabreiko \cite{KrasZ} for fundamentals of nonlinear operator theory),
   applications of many
 known classic techniques associated with some remarkable and  famous specific PDEs become
 non-aplicable, so that a principally new methodology and even
 philosophy of nonlinear PDEs via blow-up analysis should play a part.

\ssk









\noi\underline{\sc Reaction-diffusion models with blow-up.}
 The
second-order quasilinear reaction-diffusion equations from
combustion theory are widely used in the mathematical literature
and many applications. This class of parabolic PDEs of principal
interest in the twentieth century includes  classic models
beginning with the {\em Frank-Kamenetskii equation} with
exponential nonlinearity ({\em solid fuel model}, 1938
\cite{ZBLM}),
 \be
 \label{1.1}
  \begin{matrix}
  u_t= \D u + {\mathrm e}^u, \qquad \qquad\quad \ssk\ssk \\
 u_t= \D u + u^p, \qquad \qquad\quad \ssk\ssk \\
 u_t = \D (u^{n+1}) + u^p, \qquad\,\,\, \ssk\ssk \\
 u_t =  \n \cdot(|\n u|^n \n u) + u^p,
  \end{matrix}
   \ee
   etc.,
   where $n >0$ and $p>1$ are fixed exponents,
    and similar equations with more general nonlinearities. Due to the
   superlinear behaviour of the source terms ${\mathrm e}^u$ or $u^p$ as $u \to
   +\infty$, these PDEs are known to create finite-time blow-up
   in the sense that a bounded solution ceases to exist and
  \be
   \label{bl1}
  \mbox{$
 \sup_{x\in \re} \, u(x,t) \to + \infty \quad \mbox{as} \quad t
 \to T^-<\infty.
  $}
   \ee

 If blow-up happens, the first key question is the behaviour of
 solutions as $t \to T^-$ that reflects both mathematical and
 physical-mechanical essence of these phenomena.
 Such singular limits create a class of one of the most difficult
 asymptotic problems in nonlinear PDE theory.
  The internal
 structure of these blow-up
   singularities can be rather complicated even for simplest
   combustion models in (\ref{1.1}). It is worth mentioning that a full classification of the
 blow-up patterns is still not available for the last two
 {\em quasilinear} equations in this   list.


   In combustion modelling, the blow-up phenomena (\ref{bl1}) are treated as
 {\em adiabatic explosions} (first results were due to Todes, 1933,
 in an ODE model), as an extremal type of instabilities in such nonlinear systems.
 The research in the area of blow-up combustion processes  had been essentially intensified
 since the 1970s, when blow-up ideas first time penetrated
  key principles of the Laser Fusion (E.~Teller's famous Report of 1972
  on possibility of blow-up shockless
   compression of a D-T drop); see surveys in \cite[p.~401]{GVaz}
   and \cite[\S~1.2]{GS1S-V}.
   There exists a huge amount of mathematical papers devoted to
   various aspects of blow-up and other singularity behaviour in
   parabolic equations and systems.
 We refer to monographs \cite{BebEb,  GalGeom, AMGV, Pao, QSupl, SGKM},
     where detailed structures of blow-up are
 mainly studied and described
for second-order reaction-diffusion PDEs. In \cite{MitPoh},
 a systematic nonlinear capacity approach is developed to
 prove finite-time blow-up for wide classes of nonlinear evolution
 PDEs and systems of various  orders and types.
 In
\cite[Ch.~3--5]{GSVR}, blow-up structures for higher-order
parabolic, hyperbolic, and nonlinear dispersion  equations are
treated by means of exact solutions on invariant subspaces for
nonlinear operators.


 First study and classification
 of finite families of
  single point blow-up
patterns ($p>n+1$)
 for the second-order counterpart
 \be
 \label{sec1}
 u_t = (u^{n+1})_{xx} + u^p \quad (u \ge 0)
  \ee
  have been known since 1970s. Basic new ideas of blow-up and
  heat localization for (\ref{sec1})
were developed by Kurdyumov and his Russian School of blow-up
regimes.
   See key references, main results,
 and history in \cite[Ch.~4]{SGKM}.


\ssk

\noi\underline{\sc Higher-order diffusion and structure of
blow-up.} Much less is known still on formation of blow-up
singularities in parabolic PDEs with higher-order diffusion that
is also of essential demand in modern application. Even  simplest
 PDEs such as
 the {\em extended Frank-Kamenetskii equation} in one dimension,
 \be
 \label{1.2}
 u_t = -u_{xxxx} +{\mathrm e}^u \quad \big(\mbox{or}
 \quad  u_t = u_{xxxxxx} +{\mathrm e}^u\, \big),
  \ee
  and their  counterparts with power source
\be
 \label{1.3}
 u_t = -u_{xxxx} +|u|^{p-1}u \quad \big(\mbox{or}
 \quad  u_t = u_{xxxxxx} +|u|^{p-1}u\, \big),
  \ee
reveal several principally new asymptotic blow-up phenomena
demanding novel  mathematical approaches; see details in
  \cite{BGW1, Gal2m}.
  Similar difficulties occur for
   the model  equation from the {\em Semenov--Rayleigh--Benard
  problem} with the leading operator of the form
 \be
 \label{1.4}
 u_t= - u_{xxxx} + \b[(u_x)^3]_x + {\mathrm e}^u \quad (\b \ge 0);
 \ee
see \cite{GW1}. The mathematical difficulties in understanding the
ODE and PDE patterns  increase dramatically with the order of
differential diffusion operators in the equations. Interesting
{\em regional blow-up} (see a  definitions below) and oscillatory
properties \cite{CG2m} are exhibited by a semilinear diffusion
equation with ``almost linear" logarithmic source term
 \be
 \label{1.41}
 u_t= -u_{xxxx} + u \, \ln^4 |u|.
  \ee

There are just a few examples of higher-order parabolic
reaction-diffusion PDEs with reasonably well-understood (but still
not fully proved!) blow-up singularities. On the other hand, as we
have mentioned, the techniques for proving global nonexistence of
solutions are already available and cover a variety of nonlinear
evolution equations of arbitrary orders; see \cite{MitPoh}.
  Note that all
the above models (\ref{1.2})--(\ref{1.41}) are semilinear and do
not admit blow-up patterns with finite propagation (i.e., with
finite interfaces)  and hence with {\em free boundaries} that are
of special interest in PDE and ODE theory.

\subsection{Fourth-order  PME-type equation with source: on the model,
 techniques, and main results}

In the present paper, we propose and study self-similar blow-up
for the following
{\em fourth-order porous medium equation with source} (PME$-4$
with source) :
 \be
 \label{1.5}
  \fbox{$
   \ssk\ssk
u_t = {\bf A}(u) \equiv - \big(|u|^n u \big)_{xxxx} + |u|^{p-1}u
\quad \mbox{in} \,\,\, \re \times \re_+,
   \ssk\ssk
 $}
 \ee
 where $n >0$ and $p>1$. For $n=0$, this gives the semilinear
 equation (\ref{1.3}) describing single point
 blow-up for all $p>1$, \cite{BGW1}.
  For $n>0$, blow-up phenomena and the corresponding mathematics
 for (\ref{1.5}) are unknown and clearly become more involved.

We consider for (\ref{1.5}) the Cauchy problem with  bounded
compactly supported data
 \be
 \label{u01}
 u(x,0)=
 u_0(x) \quad \mbox{in}
 \quad \re.
  \ee
 Unlike the well-known from the 1980s {\em thin film equations}
 (TFEs) (see (\ref{GPP4}) below), the PDE (\ref{1.5}) contains a
 fully divergent, monotone, and potential,  differential operator.
 Moreover,
   (\ref{1.5}) is a gradient
system and admits strong estimates via multiplication by $(|u|^n
u)_t$ in $L^2$.
 The PME operator is potential in the metric of $H^{-2}$ and is also
monotone, so local existence of a unique continuous solution
follows from classic theory of monotone operators; see
\cite[Ch.~2]{LIO}.
 Finite
propagation  in the PDE (\ref{1.5}) 
 is proved by energy estimates via Saint--Venant's principle; see
various techniques in \cite{Bern01} and more recent results in
\cite{Shi2} and also a survey in \cite{GS1S-V} for further
references.

Global nonexistence results for PDEs like (\ref{1.5}) are known
from the beginning of the 1980s; see \cite{Gal82}. Concerning more
recent results, let us mention that any solution of the Cauchy
problem for
 $$
u_t = - \D^2(|u|^n u) + |u|^p \quad \mbox{in} \quad \ren \times
\re_+, \quad u(x,0)=u_0(x) \,\,\, \mbox{in} \,\,\, \ren
 $$
 (recall the difference in the source term $|u|^p$ of the equation,
   which is associated with the nonlinear capacity method) blows up
in the {\em subcritical Fujita range}
 $$
  \mbox{$
 n+1 < p \le p_0= n+1 + \frac 4N, \quad \mbox{provided that}
 \quad \int u_0 >0
  $}
  $$
  (the first condition $n+1<p$ is also purely technical);
 see \cite{Eg40} for a brief simple exposition of the method
 and \cite{MitPoh} for a systematic treatment.

 It turns out that studying  and proper understanding
 of higher-order models such as (\ref{1.5})
  demand  completely different and more difficult
mathematics, and, often, the known mathematical methods of classic
nonlinear operator theory are not sufficient to answer some even
basic and principal questions arisen. In other cases, tremendous
efforts are necessary to justify their applicability for the
nonlinear higher-order non-potential and non-variational operators
involved. As happened before for the classic combustion model
(\ref{sec1}),  complicated blow-up singularities for (\ref{1.5})
hopefully are helping to initiate new mathematical directions that
will take a long time to be proper developed.
 Several open problems appear that are stated throughout the
 text. Some of them are extremely difficult so the author just
 states and discusses these  without specific mathematical efforts
 to solve.

\ssk

\noi\underline{\sc Local problem: oscillatory sign-changing
behaviour at finite interfaces.} It is curious that though the
principal fact on the oscillatory behaviour of solutions close to
interfaces of the pure {\em fourth-order PME} (PME$-4$)
 \be
 \label{p1}
 u_t= - \big(|u|^n u \big)_{xxxx} \quad (n>0)
  \ee
 was rigorously established (together with existence and uniqueness of the
 Zel'dovich--Kompaneetz--Barenblatt source-type solutions) by Bernis \cite{Bern88} (1988) and in
  Bernis--McLeod \cite{BMc91}
 (1991), a detailed generic structure of such oscillations for (\ref{p1}) was
 not still fully
 addressed in mathematical literature. This was one of the principle open problem of local
 theory of higher-order PMEs. Therefore, we begin our analysis
 with this local oscillatory phenomenon that is also key for some
similarity blow-up profiles:

\ssk

{\bf (0)}  oscillatory behaviour of solutions near interfaces
(Section \ref{SectLocR}).

\ssk

\noi\underline{\sc Global blow-up problems: patterns of S-, LS-,
and HS-regimes.} Concerning further evolution properties  of
(\ref{1.5}), one of the main mathematical problems and
difficulties is associated with the description and classification
of blow-up patterns occurring in finite time.
 We reveal three classes of similarity blow-up for
(\ref{1.5}) in the ranges:


 {\bf (i)} $p=n+1$, {\em regional blow-up} or {\em S-regime of blow-up} in
  Kurdyumov's terminology \cite[Ch.~4]{SGKM} (Section
  \ref{SectS}): this means that the solutions blows up as $t \to
  T^-$ in a bounded {\em localization domain};


 {\bf (ii)} $p>n+1$, {\em single point blow-up} or {\em LS-regime} (Section \ref{SectLS}):
 the {\em final-time profile} $|u(x,t)| \to \iy$ as $t \to T^-$ at a single point and
 is uniformly bounded away from it; and


{\bf (iii)} $p \in (1,n+1)$, {\em global blow-up} or {\em
HS-regime} (Section \ref{SectHS}): $|u(x,t)| \to \iy$ as $t \to
T^-$ uniformly on any compact subset in $x$.


In all the three cases, we detect the first $p_0$-branch of
$F_0$-similarity profiles, which are expected to be
 {\em evolutionary} ({\em structurally}) stable, i.e., describe generic
 formation of blow-up singularities for the PDE \ref{1.5}). We
 present no proof (it is assumed to be very difficult) and no numerical evidence (a PDE
 modelling is then necessary that is out of our goals and reach here), and rely
 on our previous experience.

 For $p \in (1,n+1]$, the similarity patterns are shown to be
{\em compactly supported} and oscillatory near interfaces, while
for $p>n+1$ they are not.
It is also crucial that,
 in the case of regional blow-up for  $p=n+1$, the problem for similarity profiles
 becomes {\em variational}. Therefore, in Section \ref{SectS}, we detect
 countable families of various compactly supported similarity solutions by
 Lusternik--Schnirel'man category theory \cite[\S~57]{KrasZ}  and spherical fibering
 approach \cite{Poh0, PohFM}. Our rather natural and  practical
    approach is {\em to use variational families of
 similarity profiles for $p=n+1$ as the origin
 of $p$-branches for $p \not = n+1$ appeared at the branching point
 $p=n+1$.}
In  other words,
  the continuation in the parameter
 $p$ makes it possible to observe and classify several blow-up patterns and $p$-branches in the
 non-variational cases $p>n+1$ and $p<n+1$. In Sections
 \ref{SectLS} and \ref{SectHS}, we present typical
 $p$-bifurcation (branching) diagrams.
For the higher-order ODEs appeared, several results are essentially  based on refined numerical
 experiments (by using standard codes supplied {\tt MatLab}),
  which often become the only tool to check global extensibility properties of
  various branches of similarity profiles. We do not think that, in a certain generality,
  all of these conclusions can be justified completely rigorously.

It is worth mentioning another
 fruitful
  approach
  to use the analogy with the linear bi-harmonic equation
  \be
   \label{bh1}
   u_t= - u_{xxxx} + u \quad \mbox{in} \quad \re\times \re_+,
    \ee
    which is obtained from (\ref{1.5}) by both limits
 $n \to 0$ and $p\to 1$. A countable set of exponential
 patterns for (\ref{bh1}) will be described on the basis of
 spectral theory developed in Section \ref{SectExi}.

 \subsection{On other PDEs with striking similarities of the blow-up results}

The present research  and main approaches are better understood in
conjunction with the results in \cite{GalpLap}, where another
canonical fourth-order models of the $p$-Laplacian operator with
source,
 \be
 \label{pLap1}
 u_t= -(|u_{xx}|^n u_{xx})_{xx} + |u|^{p-1}u, \quad n>0, \quad
 p>1,
  \ee
is studied. It turns out amazingly that several general
conclusions on blow-up behaviour for (\ref{1.5}) and (\ref{pLap1})
coincide, but the necessary mathematics is completely different in
many key places (e.g., the local oscillatory behaviour and the
L--S variational approach), since the diffusion operators works in
distinct functional spaces, $H^{-1}$ and $L^2$ respectively.
However, even existing similarities for these two nonlinear PDEs
are not that easy to detect, since the mathematics involved
essentially differs and the common roots of analogous results are
rather obscure still. We will refer and use some common results in
\cite{GalpLap} when necessary and appropriate. The main principles
of classification of blow-up similarity patterns and branches
 are claimed to remain the same for higher-order PDEs of similar type, e.g., for
the {\em sixth-order PME-type equation with source} (PME$-6$ with
source)
 \be
 \label{661}
 u_t=(|u|^n u)_{xxxxxx} + |u|^{p-1} u.
  \ee
  Moreover,
 further developing this paradoxical phenomenon of analogy,
  the author
claims that a number of local and similarities and  geometrically
equivalent  set of nonlinear eigenfunctions are available for
nonlinear degenerate PDEs such as (these awkward and artificial
examples are chosen for the sake of the argument; note that these
are not that far from real applications currently needed)
 \be
 \label{nnn1}
 u_{ttt}=(|u|^{n-1} u)_{xxxxxxxx} + |u|^{p-1}u
 \quad \mbox{or} \quad
  u_{ttttt} =
(|u_{xxxxx}|^n u_{xxxxx})_{xxxxx} - (|u|^{p-1}u)_{xx}.
   \ee
These PDEs are even hardly classified being a mixture of parabolic
(mainly this), as well as having some features of  wave, and
nonlinear dispersion types. Such equations can be shown to exhibit
blow-up
 patterns with oscillatory finite interfaces, where the
 mathematics and patterns formation/classification remain
 similar but of course more involved.

 \subsection{On  other higher-order PDEs with blow-up or extinction
 and
finite interfaces}

For convenience of the Reader and for completeness of our short
survey on  blow-up in nonlinear higher-order parabolic PDEs,
 the {\em interface} and {\em blow-up phenomena} are natural and most well-known  for the
degenerate  unstable
 {\em thin film equations} (TFEs)
 \be
\label{GPP4}
    u_t = -\big(|u|^n u_{xxx}\big)_x - \big(|u|^{p-1}u \big)_{xx} \quad
\mbox{with} \,\,\, n>0, \,\, p>1.
  \ee
Equations of this form  admit non-negative solutions constructed
by special sufficiently ``singular" parabolic approximations of
nonlinear coefficients that often  lead to free-boundary problems.
This direction was initiated by the pioneering paper \cite{BF1}
and was continued  by many researchers; we refer to \cite{LPugh,
WitBerBer} and the references therein.
 Blow-up similarity solutions of the fourth-order TFE (\ref{GPP4})
 have been also well studied and understood;
 see \cite{BerPugh1, BerPugh2, Bl4, SPugh, WitBerBer}, where
 further references on the mathematical properties of the models
 can be found. Oscillatory solutions and countable sets of blow-up
 patterns for this TFE were described in \cite{Bl4}.

Interface and
 {\em extinction}  behaviour occur for other
reaction-absorption PDEs such as
 \be
 \label{GPP}
\mbox{$
 u_t=  
 -u_{xxxx} - |u|^{p-1}u
 $}
  \ee
in the {\em singular} parameter range  $p \in (- \frac 13,1)$,
where the absorption term $|u|^{p-1}u$ is not Lipschitz continuous
 at $u=0$. This creates interesting evolution extinction phenomena
 and various sets of similarity patterns;
see \cite{Galp1} and references therein.

 \section{Fundamental solution and spectral properties for $n=0$}
 \label{SectExi}

These questions have been discussed in \cite[\S~3.3]{GalpLap}, so
we  briefly report on those auxiliary properties.
 Consider the linear {\em bi-harmonic equation} corresponding to
 $n=0$,
 \be
 \label{bi1}
 u_t=- u_{xxxx} \quad \mbox{in} \quad \re \times \re_+.
  \ee
    Its {\em fundamental
  solution} has the form
   \be
   \label{bi2}
    b(x,t)= t^{-\frac 14} F(y), \quad y = x/t^{\frac 14}, \quad
    \mbox{where}
     \ee
  \be
\label{ODEf}
 {\bf B} F \equiv - F^{(4)} +  \mbox{$\frac 1{4}$}\, y
 F' + \mbox{$\frac 1{4}$}\, F =0 \quad
 \mbox{in} \,\,\, \re,  \quad \mbox{with} \,\,\, \mbox{$\int$} \, F = 1.
 \ee
The rescaled kernel $F=F(|y|)$ is radial, has exponential decay,
oscillates as $|y| \to \infty$,
 and   \cite[p.~46]{EidSys}
\be \label{fbar}
 |F(y)| \le  D {\mathrm e}^{-d|y|^{4/3}}
\quad \mbox{in} \,\,\, \re,
 \ee
  for some
positive constant $D$ and
  $d = 3 \cdot 2^{-11/3}$.
The necessary spectral properties of the linear
operator ${\bf B}$ and the corresponding adjoint operator ${\bf
B}^*$
  can be found in \cite{Eg4} for similar more general $2m$th-order
 operators (see also \cite[\S~4]{Bl4}).
 In particular, ${\bf B}$ is naturally defined
in a weighted space $L^2_\rho(\re)$, with
 $
 \rho(y)={\mathrm e}^{a
|y|^{4/3}}$,
where $a \in (0,2d)$  is a constant,
  and  has the discrete (point) spectrum
 \begin{equation}
\label{spec1}
  \sigma({\mathbf B}) = \big\{\lambda_l = -\mbox{$\frac
l4$}, \,\, l = 0,1,2,...\big\}.
\end{equation}
The corresponding eigenfunctions are given by  differentiating the
kernel,
 \be
 \label{psi1}
 \mbox{$
 \psi_l(y)= \frac {(-1)^l}{\sqrt{l!}} F^{(l)}(y), \quad l=0,1,2,... \,
 .
  $}
  \ee
The adjoint operator
 \be
 \label{ad1}
  \mbox{$
 {\bf B}^*=- D_y^4 - \frac 14 \, y D_y
  $}
  \ee
  has the same spectrum (\ref{spec1}) and
  eigenfunctions being generalized Hermite polynomials
  \begin{equation}
 \label{psi**1}
  \mbox{$
 \psi_l^*(y) = \frac 1{\sqrt{l !}}
\sum_{j=0}^{\lfloor-\lambda_l\rfloor} \frac {1}{j !}D^{4j}_y y^l,
\quad l=0,1,2,...\, ,
 $}
 \end{equation}
which form a complete and closed set  in $L^2_{\rho^*}(\re)$,
where $\rho^*(y)=\frac 1 {\rho(y)}$.
 The bi-orthonormality of eigenfunction sets holds:
 \be
 \label{ort1}
 \langle \psi_l, \psi^*_k \rangle = \d_{lk},
  \ee
  where $\langle \cdot, \cdot \rangle$ denotes the standard (dual)
  scalar product in $L^2(\re)$.

\section{Local asymptotic properties of solutions near interfaces}
  \label{SectLocR}

  We study the  oscillatory behaviour
  of solutions of (\ref{1.5}) close to  interfaces.

\subsection{Local properties of  
travelling waves (TWs)}


As customary,
 we use  TWs
     solutions,
 \be
 \label{t1}
 u(x,t)=f(y), \quad y=x- \l t,
 \ee
to describe propagation properties for
 (\ref{1.5}). Then $f(y)$ solves
 \be
 \label{t2}
 -\l f'= - (|f|^n f)^{(4)} +|f|^{p-1}f.
  \ee
  By a local analysis near the interface, at which $f=0$ (and, as we will see, also
  $f'=f''=f'''=0$ for all $n>0$),
  it
  is not difficult to see that the higher-order term $|f|^{p-1}f$
  on the right-hand side is negligible.
Therefore, near moving interfaces,
   we  consider the simpler equation
 \be
 \label{le2}
   (|f|^n f)'''= -f \quad \mbox{for} \quad y>0, \quad f(0)=0,
  \ee
 which is obtained on integration once.
 Here we set $\l=-1$  by scaling for propagating waves.
 It is convenient to use the natural change
  \be
  \label{le2N}
  F= |f|^n f \quad \Longrightarrow \quad F'''=- \big|F\big|^{-\frac
  n{n+1}}F.
   \ee
 We  describe
oscillatory solutions of changing sign of  (\ref{le2N}), with
zeros concentrating at the given interface point $y=0^+$.
 {\em Q.v.} the pioneering papers by Bernis and McLeod \cite{Bern88, BMc91}.

Actually, we are going to give a sharp description of the
behaviour of the solutions close to interfaces. By the scaling
invariance of (\ref{le2N}), we look for solutions of the form
 \be
 \label{le3}
  \mbox{$
 F(y) = y^\mu \varphi(s), \quad s= \ln y, \quad
 \mbox{where} \,\,\, \mu = \frac {3(n+1)}{n}> 3 \,\,\,\mbox{for} \,\,\, n>0,
  $}
  \ee
  where $\varphi(s)$ is the
  {\em oscillatory component}.
 Substituting (\ref{le3}) into (\ref{le2N}) yields
 the ODE
 \be
 \label{le4}
  P_3(\vp)= - |\vp |^{-\frac
  n{n+1}} \vp.
  \ee
   Linear differential operators $P_k$
 are given by the recursion
   (see
   \cite[p.~140]{GSVR})
    $$
    \begin{matrix}
    P_{k+1}(\varphi)= P'_k(\varphi) + (\mu-k) P_k(\varphi), \quad
    k \ge 0; \quad P_0(\varphi)=\varphi, \quad \mbox{so that}
    \smallskip\smallskip\\
P_1(\varphi)=\vp'+ \mu \vp, \quad P_2(\vp)= \vp''+(2\mu-1) \vp'+
\mu(\mu-1)\vp, \smallskip\smallskip \\
  P_3(\varphi)=\vp'''+ 3(\mu-1) \vp'' + (3 \mu^2- 6 \mu +2) \vp'
 + \mu(\mu-1)(\mu-2) \vp,  \,\,\, \mbox{etc}.
 \end{matrix}
 $$
In view of (\ref{le3}), we look for uniformly bounded global
solutions $\varphi(s)$ that are defined for all $s \ll -1$, i.e.,
can be extended to the interface at $y = 0^+$. The best candidates
for such global orbits of (\ref{le4}) are periodic solutions
$\varphi_*(s)$ that are defined for all $s \in \re$. Indeed, they
can describe suitable (and, possibly, generic) connections with
the interface at $s=-\infty$. Existence of such a periodic
solution $\vp_*(s)$ of (\ref{le4}) can be achieved by shooting
arguments; see details in \cite[\S~7]{Gl4}, where further
references  concerning periodic orbits of higher-order ODEs are
given. Uniqueness of $\vp_*$ remains an open mathematical problem,
though was always confirmed numerically.

\begin{figure}
\centering \subfigure[$n=0.75$]{
\includegraphics[scale=0.52]{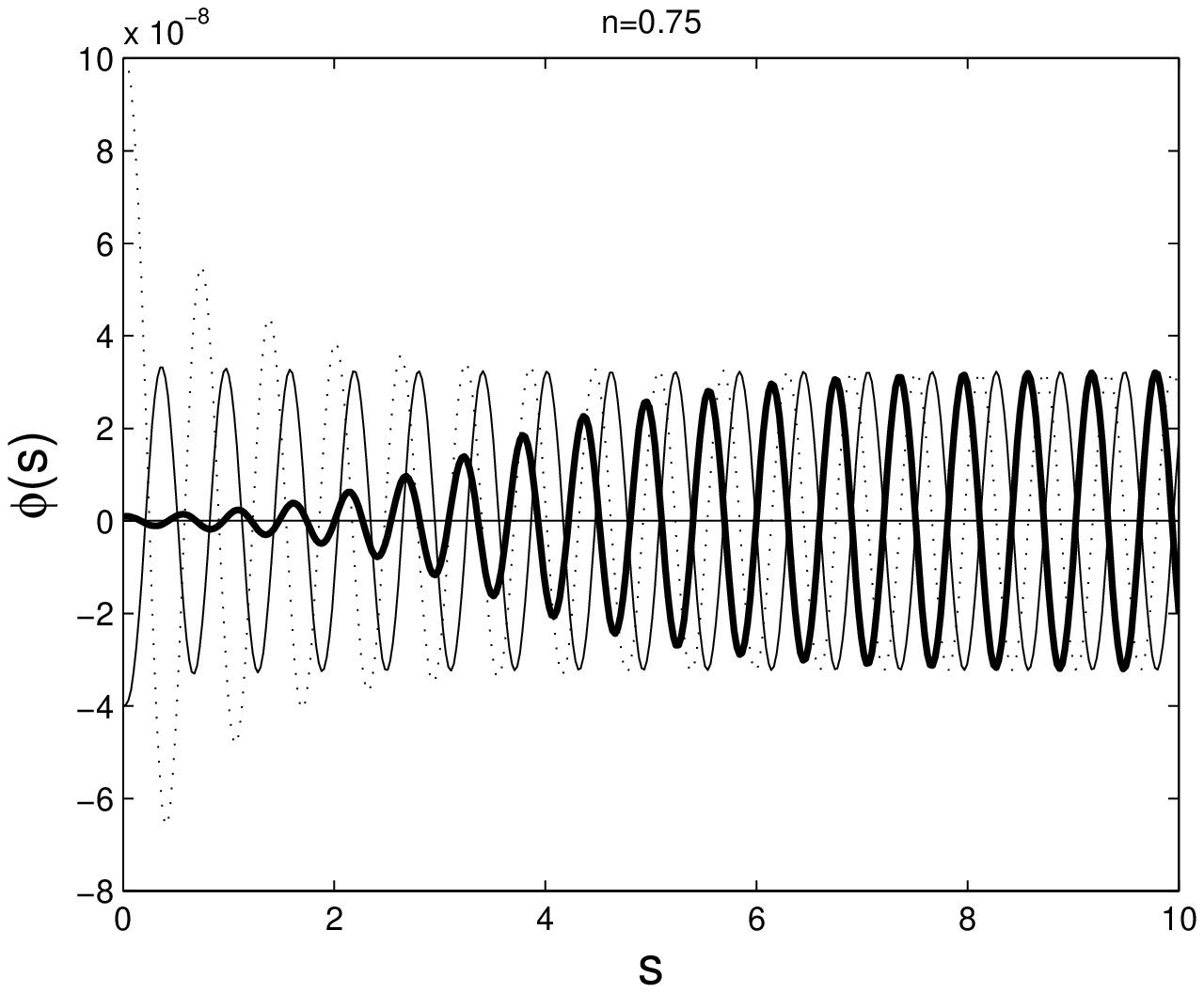}
} \subfigure[$n=1$]{
\includegraphics[scale=0.52]{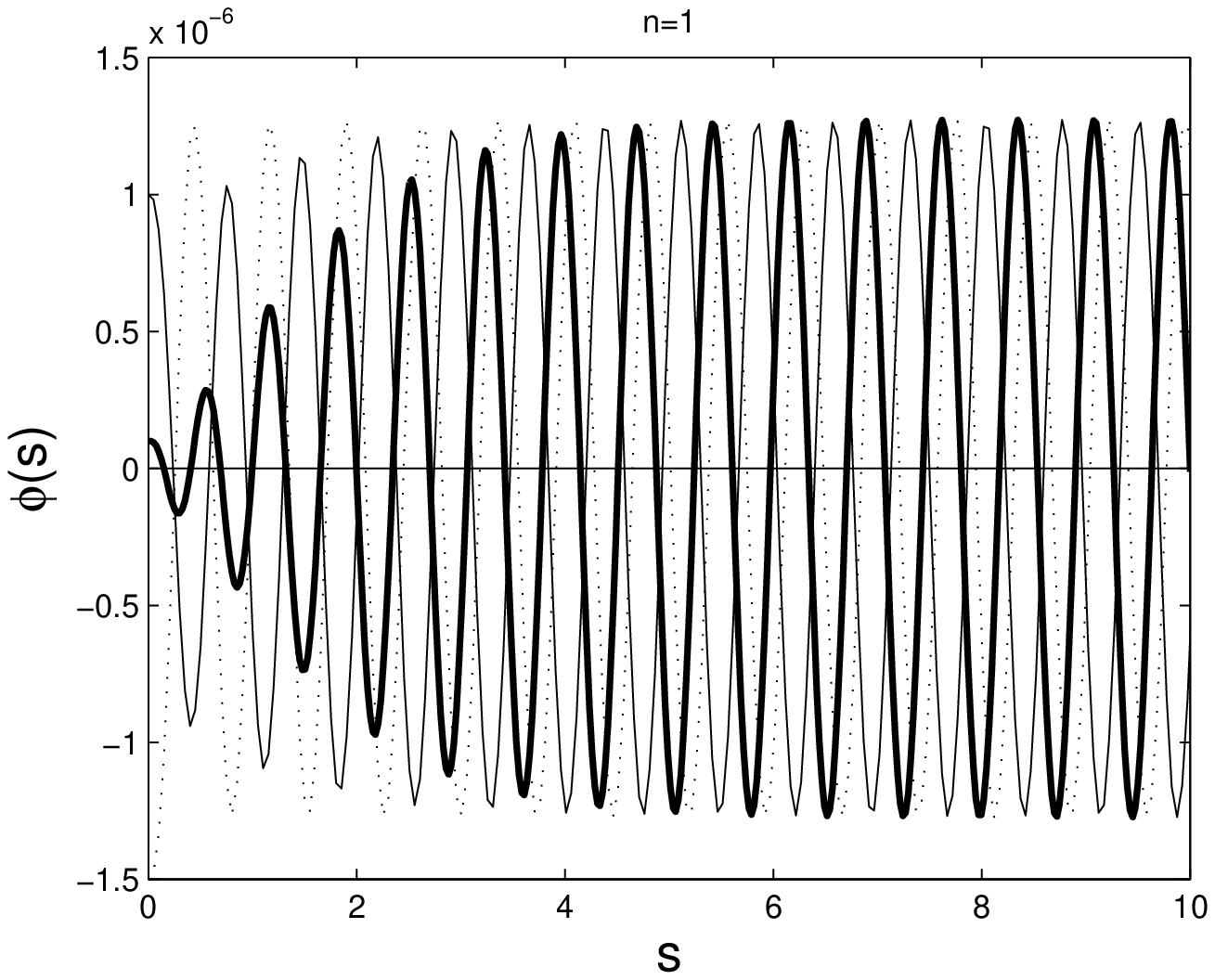}
}
 \subfigure[$n=2$]{
\includegraphics[scale=0.52]{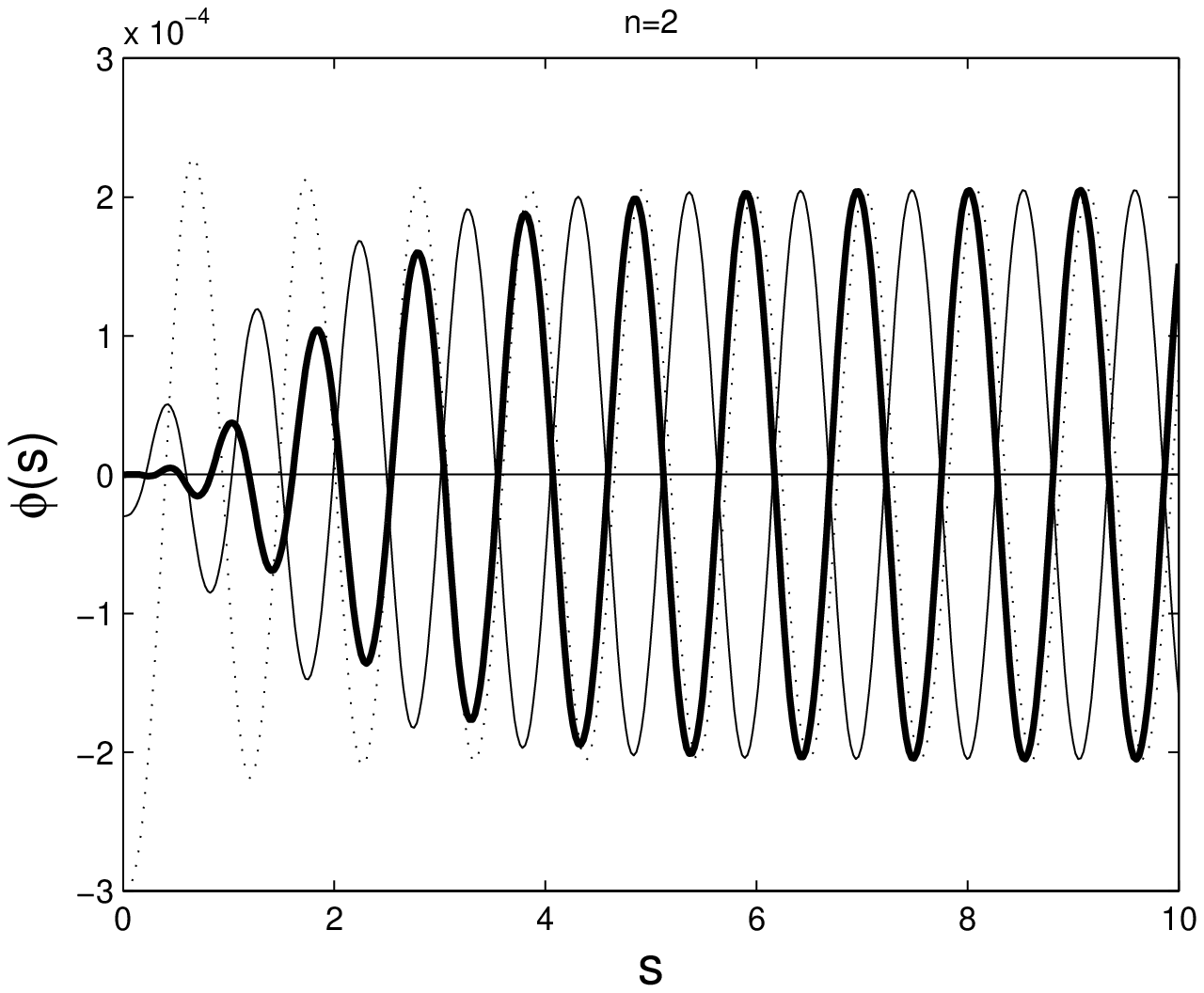}
} \subfigure[Large $n$]{
\includegraphics[scale=0.52]{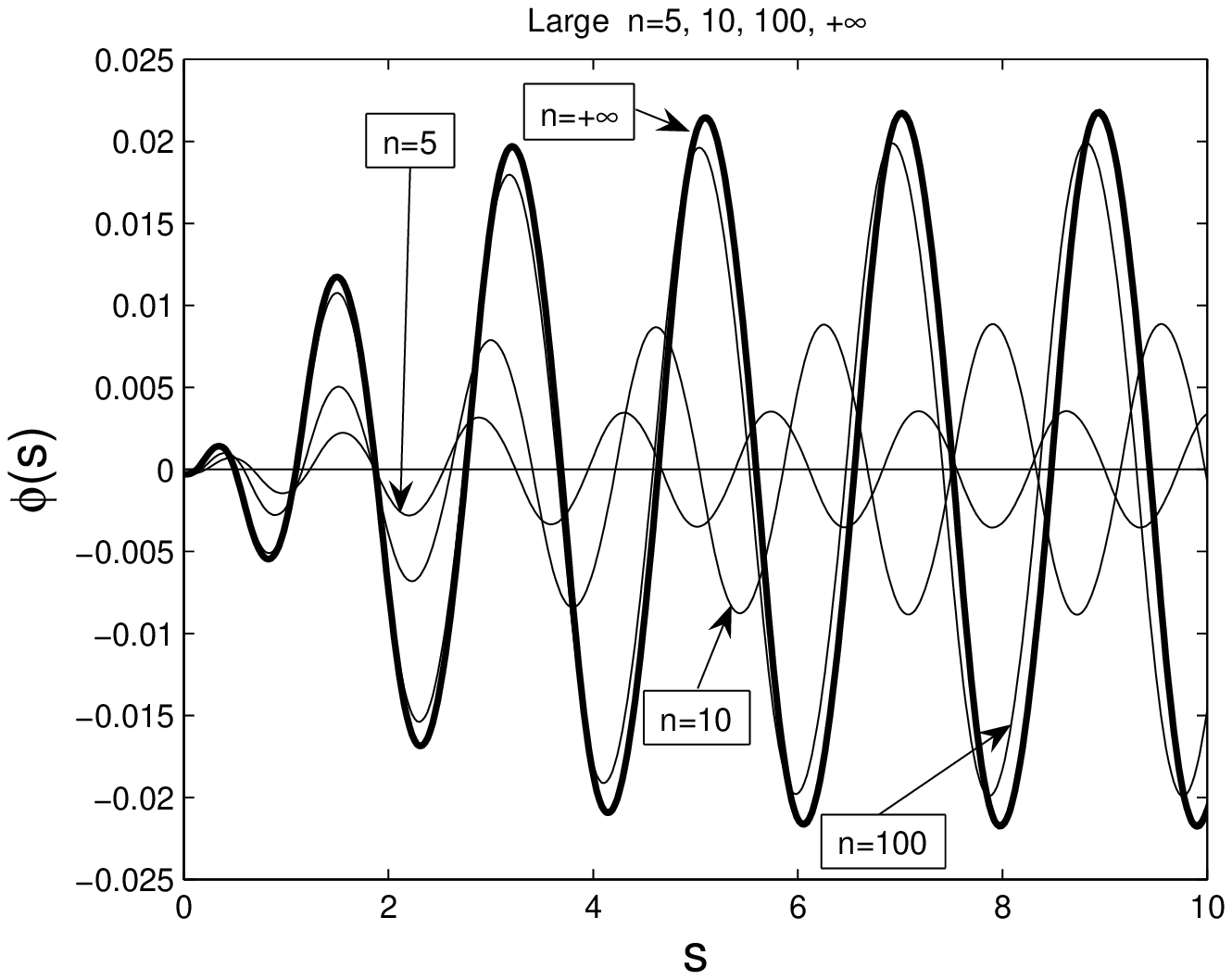}
}
 \vskip -.4cm
\caption{\rm\small Convergence to a stable periodic orbit
(\ref{le2N}) for various $n>0$.}
 \label{FOsc2}
\end{figure}

  Figure \ref{FOsc2}  shows fast convergence to such a unique (numerically) stable periodic
 solution of (\ref{le4}) for
 various $n > 0$. These numerics are obtained by using the ODE solver {\tt ode45}
  {\tt MatLab}
 with enhanced tolerances and a regularization in the singular term up to $10^{-12}$; see (\ref{4.1})
 and more comments in Section \ref{SectS}, where more advanced numerical techniques for solving boundary-value problems are employed.
 Different  curves therein correspond to different Cauchy data
$\varphi(0)$, $\varphi'(0)$,  $\varphi''(0)$  prescribed at $s=0$.
For $n$ smaller than $\frac 34$, the oscillatory component gets
extremely small, so an extra scaling is necessary, which is
explained in \cite[\S~7.3]{Gl4}. A more accurate passage to the
limit $n \to 0$ in (\ref{le4}) is done there in \S~7.6 and in
Appendix B.   In (d), we also present  the periodic solution for
$n=+\infty$ where (\ref{le4}) takes a simpler form,
 $$
 P_3(\varphi)=- {\rm sign} \, \varphi.
 $$


Finally, given the 
periodic $\varphi_*(s)$  of (\ref{le4}),
 as a natural
way to approach the interface point $y_0=0 $ according to
(\ref{le3}), we have that the ODE (\ref{le2})
 and, asymptotically,
 (\ref{t2}),
  generate at the singularity set $\{f=0\}$
 \be
 \label{as55}
\mbox{a 2D local asymptotic bundle with parameters $y_0$ and phase
shift in $s \mapsto s + s_0$}.
 \ee
This two-dimensional bundle will be shown to be sufficient for
matching with, typically, two boundary conditions that generate
blow-up patterns.


\subsection{Non-oscillatory case $\l >0$: 1D bundle of non-oscillatory asymptotics}

For $\l=+1$, we have the opposite sign in the ODE
 \be
 \label{le4N}
  P_3(\vp)=   \big|\vp \big|^{-\frac
  n{n+1}} \vp,
  \ee
  which admits two  equilibria
 \be
 \label{co1}
 \mbox{$
 \varphi_\pm= \pm \big[\mu
 (\mu-1)(\mu-2)\big]^{-\frac{n+1}n}.
  $}
  \ee
 Let us check the dimension of their stable manifolds as $s \to +
 \iy$ and, most importantly, $s \to - \iy$, which corresponds to
 approaching the interface. Note that, for (\ref{le4}) the trivial
 equilibrium $\varphi_0=0$ is evidently unstable in both
 directions with empty stable manifolds (in view of the
 non-Lipschitz nonlinearity at $\varphi=0$ on the right-hand
 side).

Thus, by  standard linearization, it follows that both are stable
as $s \to +\infty$: setting $\varphi = \varphi_+ +Y$ yields the
linear ODE
 \be
 \label{lin1}
 \mbox{$
 {\bf P} Y \equiv
\big[P_3- \frac 1{n+1} \,\mu(\mu-1)(\mu-2)I\big]Y=0.
  $}
  \ee
The characteristic equation by plugging $Y={\mathrm e}^{\l s}$
into (\ref{lin1}) takes the form
 \be
 \label{ch1NN}
 p(\l) \equiv \l^3 + 3(\mu-1)\l + (3 \mu^2-6 \mu+2) \l +
 3(\mu-1)(\mu-2)=0, \quad \mbox{where}
  \ee
 $$
  \mbox{$
 p(0)=3(\mu-1)(\mu-2)>0,  
 \quad p'(\l)=0 \,\,\, \mbox{at} \,\,\, \l_\pm= \frac
 13\big[-3(\mu-1) \pm \sqrt 3 \,  \big]<0.
  $}
  $$
It follows that all characteristic values of ${\bf P}$ satisfy
 \be
 \label{ch2}
  {\rm Re} \, \l_k < 0, \quad k=1,2,3.
 \ee

\begin{figure}
\centering \subfigure[Stability of equilibria]{
\includegraphics[scale=0.52]{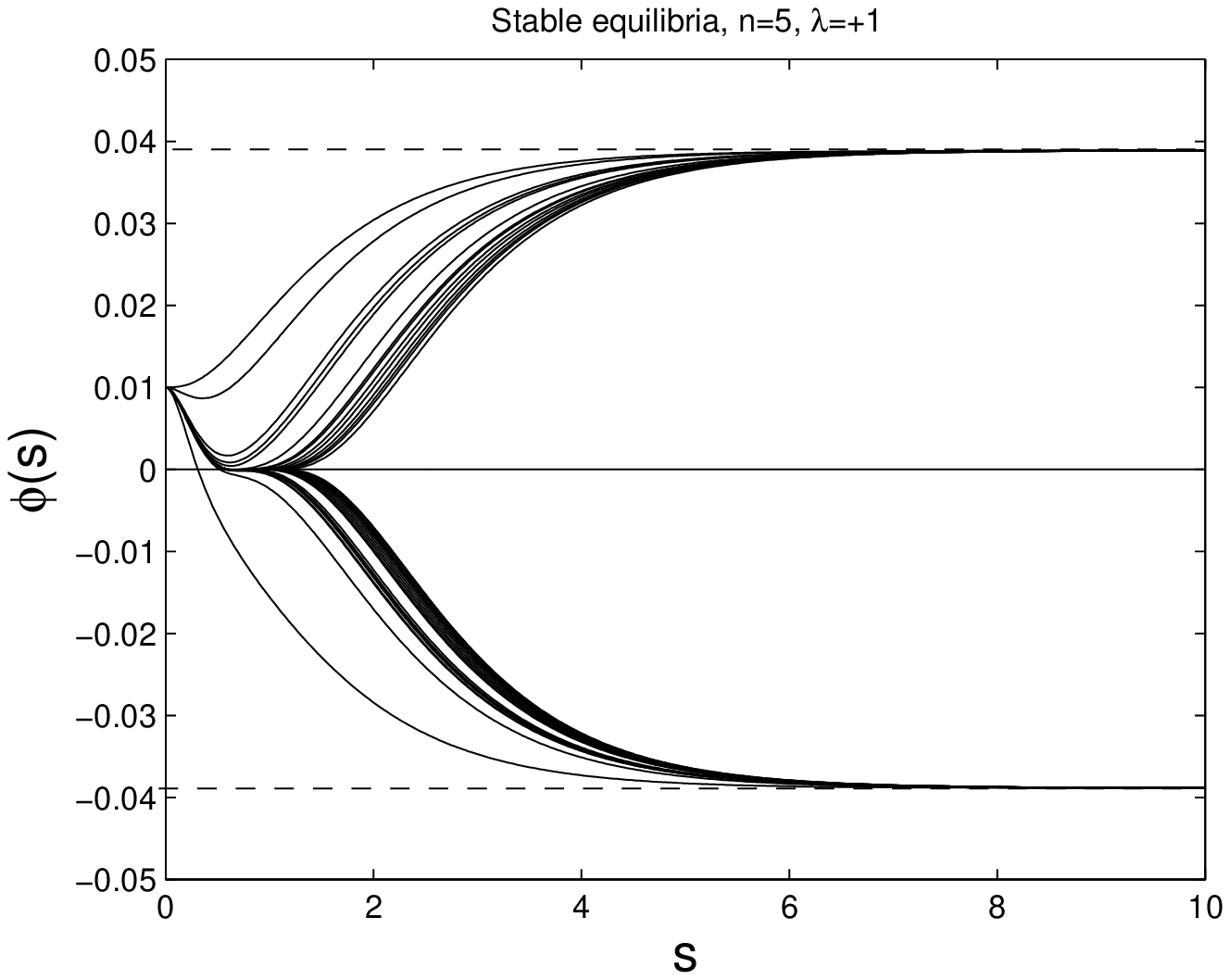}
} \subfigure[Behaviour in between]{
\includegraphics[scale=0.52]{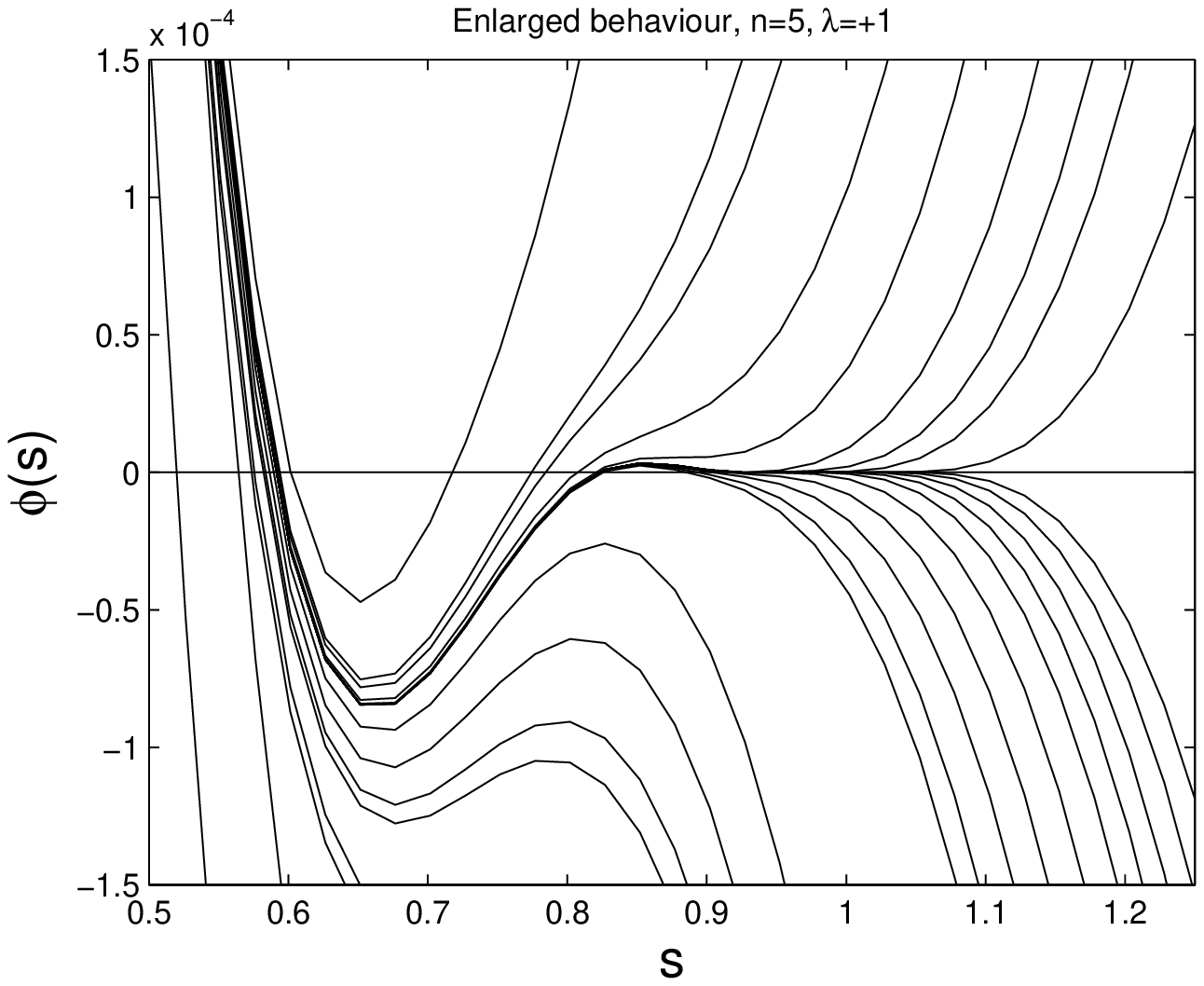}
}
 \vskip -.4cm
\caption{\rm\small Non-oscillatory behaviour for the ODE
(\ref{le4N}) for $n=5$; stability of equilibria (\ref{co1}) (a),
and enlarged non-periodic behaviour in between, (b).}
 \label{FOsc22}
\end{figure}

 Figure
\ref{FOsc22}(a) also confirms that as $s \to + \infty$ the
equilibria (\ref{co1}) are stable.
  In (b), which gives the enlarged image of the
behaviour from (a) close to $\vp=0$,  we observe a changing sign
orbit, which is not a periodic one. Therefore, this cannot be
extended as a bounded solution up to the interface at $s=-\infty$.
There are many other similar  ODEs, where such behaviour between
two equilibria is periodic; cf. \cite[p.~143]{GSVR}.

In a whole, theses results confirm that for $\l>0$, the TWs {\em
are not oscillatory at interfaces}, so that, formally, this {\em
backward} propagation can be performed by positive solutions. On
the other hand, (\ref{ch2}) establishes that  the stable manifold
of equilibria $\varphi_\pm$
 as $s \to -\infty$ (i.e., towards the interface) is empty, so
 that unlike (\ref{as55}), for $\l = +1$,
  \be
  \label{as11}
   \exists \,\,
   \mbox{a 1D local asymptotic bundle with parameter $y_0$}.
 \ee
One can see from (\ref{le3}) that these positive asymptotics are
 given by (for fixed $\l=+1$ which has been scaled out as a non-essential parameter)
 \be
 \label{as111}
 f(y)= (y-y_0)^{\frac 3n} [\mu(\mu-1)(\mu-2)]^{-\frac 1n} (1+o(1))
  \quad \mbox{as}
 \quad y \to y_0^+.
  \ee
 This 1D bundle is not enough for construction of typical global connections
 via two boundary conditions, so  that the backward
 propagation is either not possible at all for almost all (a.a.) initial data in the Cauchy problem
 for (\ref{p1}), or it is
 not performed by TWs.

In other words, this shows that, for the PME$-4$ (\ref{p1}), the
forward propagation of interfaces is common. Concerning the
backward one, though  existing for specific initial data, it is
not plausible in general. Recall that for the standard PME
 \be
 \label{p21}
 u_t=\big(|u|^n u \big)_{xx},
  \ee
the backward propagation of non-negative solutions is completely
prohibited by the Maximum Principle (by the straightforward
comparison from below with slowly moving small TWs). For
higher-order equations, these barrier techniques via comparison
fail, but anyway, via  the reduction of the bundle dimension in
(\ref{as11}), we justify a similar phenomenon, which is now true
not for {\em all}, but for {\em a.a.} initial data. It is worth
mentioning that the PME (\ref{p21}) also admits oscillatory
solutions near interfaces (see references and comments in
\cite[p.~30]{GalGeom}). By Sturm's First Theorem on non-increase
of the number of zeros, such behaviour is not generic for
second-order parabolic PDEs in the sense that such solutions
cannot appear from data with any finite number of sign changes.


\section{Blow-up similarity solutions: problem setting and preliminaries}
 \label{Sect2}

\subsection{Posing the similarity problem}
The parabolic PDE (\ref{1.5}) admits the following similarity
solutions with blow-up as $t \to T^-$:
 \be
    \label{RVarsE}
 \mbox{$
 u_S(x,t) = (T-t)^{-\frac 1{p-1}}  f(y),
\quad y = {x}/(T-t)^{\b},  \quad \mbox{with} \quad \b =
\frac{p-(n+1)}{4(p-1)}.
 $}
 \ee
The rescaled  profile $f(y)$ solves a quasilinear fourth-order ODE
 \be
  \label{f11E}
 \mbox{$
  {\bf A}(f) \equiv - (|f|^n f)^{(4)}  - \b \, y f' -
\frac 1{p-1} \, f + |f|^{p-1}f=0 \quad \mbox{in} \quad \re.
 $}
 \ee
Equation (\ref{f11E}) has the constant equilibria
  \be
 \label{g1N}
 \pm f_*(p)=  \pm(p-1)^{-\frac 1{p-1}} \to \pm \iy \quad \mbox{as} \quad p \to 1^+.
  \ee
 For even solutions $f(y)$,
 symmetry conditions at the origin $y=0$ are posed,
   \be
\label{BCs}
 f'(0) = 0 \quad \mbox{and} \quad f'''(0) = 0 \quad (\mbox{if} \,\,\, f(0) \not =
 0),
\ee
 while for odd ones we impose the anti-symmetry ones
  \be
\label{as9}
 f(0) = 0 \quad \mbox{and} \quad f''(0) = 0.
\ee
 A natural setting for the Cauchy problem assumes that, for $p \in
 (1,n+1]$,
   \be
 \label{com1}
 f(y) \,\, \mbox{is sufficiently smooth and compactly supported}.
  \ee
The actual regularity of $f(y)$ close to interfaces was determined
in the previous section by the local asymptotic analysis.
For $p>n+1$, the solutions are not compactly supported.


The ODE (\ref{f11E}) is a difficult fourth-order nonlinear
equation, so that, in general, existence and multiplicity results
for the above boundary-value problems reduce to complicated, at
least, 2D shooting problems (recall (\ref{as55}) and two boundary
conditions). This principally differs (\ref{f11E}) from its
second-order counterparts associated with the combustion model
(\ref{sec1}), where  the most of the results admit either a clear
phase-plane interpretation or reduce to 1D shooting in the
presence of the Maximum Principle; see typical results and
references in \cite[Ch.~4]{SGKM}.
 Bearing in mind that, for the PME$-6$ with source (\ref{661}),
the  parameter space of shooting will be 3D, we cannot rely on
these geometric ideas, and have to use a principally different
approach for reliable detecting multiple blow-up patterns. This
will be variational approach to the case $p=n+1$ plus
$p$-branching approach for $p>n+1$ and $p<n+1$.

\subsection{Blow-up self-similar profiles: preliminaries}
\label{SectEx}

We  next study global
 compactly supported solutions of the ODE
(\ref{f11E}). For $p \le n+1$,
 the local interface analysis from Section \ref{SectLocR}
  applies to 
  (\ref{f11E}). Indeed, close to the interface point
 $y=y_0>0$  of the similarity profile
$f(y)$,
 the ODE (\ref{f11E}) contains the same leading terms as in
 (\ref{le2}) and other linear two are negligible as $y \to y_0^-$.


It is key that, taking into account the local result (\ref{as55})
 and bearing in mind the two boundary  conditions in (\ref{BCs}) or (\ref{as9}), we
may expect that
 \be
 \label{co22}
 \mbox{there exists not more than a countable set $\{f_k\}$ of
 solutions}.
  \ee
This speculations assume a certain ``analyticity" hypothesis
concerning the dependence on parameters in the degenerate ODE
(\ref{f11E}), which is plausible but not easy to prove.
    Actually, this means that, relative to the parameter $p>1$,
we can expect at most a countable set of $p$-branches of
solutions. To begin with, this is true for the linear case $n=0$
and $p=1$.

\subsection{Countable set of similarity solutions for
$n=0$, $p=1$, and  homotopic connection}

These questions are discussed in detail in \cite[\S~3.4]{GalpLap},
so we briefly comment that the linear PDE (\ref{bh1}) has a {\em
countable} set of patterns given by eigenfunctions (\ref{psi1}):
 \be
 \label{u11}
 \mbox{$
 u_l(x,t) = {\mathrm e}^{-t}\, t^{-\frac {1+l}
4} \,\psi_l \big(\frac
 x{t^{1/4}}\big), \quad l=0,1,2,... \,
  $}
  \ee
Our intention is to see  how the blow-up similarity patterns
(\ref{RVarsE}) can be deformed as $n \to 0$ and  $p \to 1$ to
those in (\ref{u11}).
Then (\ref{u11}) suggests that,
for $n>0$, 
there exists a countable number of ``branches" $\{f_l(y;n,p)\}$,
which
``bifurcate" from the point $\{n=0, \,  p=1\}$.  
It is possible to  observe the space-time structures in
(\ref{u11}) in nonlinear blow-up analysis, where similar countable
sets of patterns will be shown to exist; see \cite{GW2}.
 We then say that the above two (linear for $n=0$, $p=1$ and nonlinear for $n>0$,
$p>1$) asymptotic problems admit a continuous {\em homotopic}
connection as $n \to 0$, $p \to 1$. See \cite{GalpLap} for extra
details.




\section{Regional blow-up similarity profiles for $p=n+1$}
\label{SectS}

As in \cite[\S~4]{GalpLap}, the case $p=n+1$ is simpler since
$\b=0$, so (\ref{RVarsE}) is a solution in separable variables
 $$
 u(x,t)=(T-t)^{-\frac 1n} f(x).
 $$
 Then
  $f=f(y)$ (we continue to use $y$ as the independent spatial
  variable for application to $p \not = n+1$) solves an autonomous fourth-order ODE of the form
 \be
  \label{f11ES}
 \mbox{$
  {\bf A}(f) \equiv - (|f|^n f)^{(4)}   -
\frac 1{n} \, f + |f|^{n}f=0 \quad \mbox{in} \quad \re.
 $}
 \ee
 We again use the same change as in (\ref{le2N}),
 \be
  \label{1}
  |f|^n f= n^{-\frac{n+1}n} F \quad \Longrightarrow \quad F^{(4)}=F- \big|F\big|^{-\frac
  n{n+1}}F \quad \mbox{in} \quad \re.
   \ee

This equation is much simpler than that in \cite{GalpLap}, and was
studied in \cite{GPHom}.
We briefly present
the main conclusions that will play an important role for the
cases $p>n+1$ and $p<n+1$ later on.

\subsection{Variational setting}

The operators in (\ref{1}) are potential, so the problem admits a
variational setting, so the solutions can be obtained as critical
points of a $C^1$ functional   of the form
 \be
 \label{V1}
  \mbox{$
 {E}(F)= - \frac 12  \int (F'')^2 \, {\mathrm d}y + \frac 12 \int
 F^2  \, {\mathrm d}y -\frac{1}{\nu} \, \int |F|^{\nu} \, {\mathrm d}y, \quad \mbox{where} \,\,\, \nu=\frac {n+2}{n+1}\in
 (1,2).
  $}
  \ee
  We are interested in  critical points in $W^2_2(\re)
  \cap L^2(\re) \cap L^{\nu}(\re)$.
Especially, we are interested in localized compactly supported
solutions, so we choose a sufficiently
  large interval $B_R=(-R,R)$ and consider the variational problem for (\ref{V1})
  in $W_{2,0}^2(B_R)$, where we assume Dirichlet boundary
  conditions at the end points $\partial B_R=\{\pm R\}$.
  By Sobolev embedding theorem, $W_{2,0}^2(B_R)
     \subset L^2(B_R)$ and in $L^{p+1}(B_R)$ compactly for any $p \ge 1$.
   Continuity of any bounded solution $F(y)$ is guaranteed
 by Sobolev embedding $H^2(\re) \subset C(\re)$.
We also need the following:

\begin{proposition}
 \label{Pr.CS}
  Let $F$ be a continuous weak solution of the
  equation $(\ref{1})$ such that
   \be
   \label{ffy1}
   F(y) \to 0 \quad \mbox{as} \quad y \to \infty.
    \ee
Then $F$ is compactly supported in $\re$.
 \end{proposition}

\smallskip

\noi {\em Proof.}
 Consider the corresponding parabolic
equation with the same elliptic operator,
 \be
 \label{Par1}
  \mbox{$
 w_t = - w_{xxxx} +w   -
 | w|^{p-1} w \quad \mbox{in} \,\,\, \re_+
 \times\re,
 \quad \mbox{where} \quad p = \frac 1{n+1} \in (0,1),
  $}
  \ee
  with initial data $F(y)$. Setting $w = {\mathrm e}^t \hat w$
yields the equation
  \be
  \label{inst1}
  \mbox{$
 \hat w_t = - \hat w_{xxxx} -{\mathrm e}^{-\frac n{n+1} \, t}\,
 | \hat w |^{p-1}\hat w,
  $}
  \ee
where the operator is monotone in $L^2(\re)$. Therefore, the CP
with initial data $F$ has a unique weak solution,
\cite[Ch.~2]{LIO}.
 Thus, (\ref{Par1}) has the unique solution
 $w(y,t) \equiv F(y)$.
 In the presence of the  singular absorption $-|u|^{p-1}u$, with $p < 1$,
 there occurs the phenomenon of {\em instantaneous
compactification} or {\em shrinking of the support} of the
solution for any data satisfying (\ref{ffy1}) (or even for more
general data in $L^p$-spaces).
  Such  phenomena for quasilinear
  absorption-diffusion equations for $p<1$
   have been known
  since  the 1970s.
  By energy estimates, similar
  results were proved for a number of
  quasilinear higher-order parabolic equations with non-Lipschitz
  absorption terms, \cite{Bern01, Shi2}.
 By the instantaneous compactification,  the multiplier ${\mathrm e}^{-\frac n{n+1} \,
 t}$ in the absorption in (\ref{inst1})  changes nothing.
   $\qed$

\smallskip

Thus,
 to revealing  compactly supported
patterns $F(y)$, we have to pose the problem in bounded
sufficiently large intervals.

 \subsection{L--S theory and direct
 application of fibering method}

Unlike \cite[\S~4]{GalpLap}, this application is more standard.
Namely, we apply classic
 Lusternik--Schnirel'man (L--S) theory of calculus of variations
\cite[\S~57]{KrasZ} in the form of the fibering method
  \cite{Poh0, PohFM}.
  Then the number of critical points of the
functional (\ref{V1}) depends on the {\em category} (or
Krasnosel'skii's  {\em genus} \cite[\S~57]{KrasZ}) of the
functional subset, on which fibering is taking place. The critical
points of ${E(F)}$ are  obtained by the {\em spherical fibering}
in the form
 \be
 \label{f1}
 F= r(v) v \quad (r \ge 0),
  \ee
  where $r(v)$ is a scalar functional, and $v$ belongs to the subset
   \be
   \label{f2}
    \mbox{$
    {\mathcal H}_0=\bigl\{v \in W_{2,0}^2(B_R): \,\,\,H_0(v)
     \equiv  -  \int (v'')^2 \, {\mathrm d}y +  \int
 v^2 \, {\mathrm d}y =1\bigr\}.
    $}
    \ee
The new functional
 \be
 \label{f3}
  \mbox{$
H(r,v)= \frac 12 \, r^2 - \frac 1{\nu}\, r^{\nu} \int |v|^{\nu} \,
{\mathrm d}y
 $}
  \ee
 has the absolute minimum point, where
 \be
 \label{f31}
  \mbox{$
 H'_r \equiv r-  r^{\nu-1} \int |v|^{\nu} \, {\mathrm d}y =0
  \,\,\Longrightarrow \,\,
   r_0(v)=\bigl(\int |v|^{\nu} \, {\mathrm d}y\bigr)^{\frac 1{2-\nu}},
   $}
   \ee 
  at which
   $
   H(r_0(v),v)=- \frac {2-\nu}{2\nu} \,  r_0^2(v).
    $
Introducing
 \be
 \label{f4}
 \mbox{$
 \tilde H(v) = \bigl[ - \frac {2\nu}{2-\nu}H(r_0(v),v)
 \bigr]^{\frac {2-\nu} 2} \equiv \int |v|^{\nu} \, {\mathrm d}y,
  $}
  \ee
  yields an even, non-negative, convex,
 and uniformly differentiable functional, to which
    L--S  theory applies, \cite[\S~57]{KrasZ};
      see also \cite[p.~353]{Deim}.
   Following \cite{PohFM}, searching
  for critical points of $\tilde H$ in  ${\mathcal H}_0$,
  one needs to estimate the category $\rho$
  of the set ${\mathcal H}_0$.
The details on this notation and basic results can be  found in
Berger \cite[p.~378]{Berger}. Notice that the Morse index $q$ of
the quadratic form $Q$ in Theorem 6.7.9 therein is  the dimension
of the space where the corresponding form is negatively definite.
This includes all the multiplicities of eigenfunctions involved in
the corresponding subspace.



For application, it is convenient to recall that utilizing
Berger's version \cite[p.~368]{Berger} of this minimax analysis of
L--S category theory \cite[p.~387]{KrasZ},  the critical values
$\{c_k\}$ and the corresponding critical points $\{v_k\}$ are
given by
 \be
 \label{ck1}
  \mbox{$
 c_k = \inf_{{\mathcal F} \in {\mathcal M}_k} \,\, \sup_{v \in {\mathcal
 F}} \,\, \tilde H(v),
  $}
  \ee
where  ${\mathcal F} \subset {\mathcal H}_0$ are  closed sets,
 and
 ${\mathcal M}_k$ denotes the set of all subsets of the form
  $
  B S^{k-1}
\subset {\mathcal H}_0,
 $
 where $S^{k-1}$ is a suitable sufficiently
smooth $(k-1)$-dimensional manifold (say, sphere) in ${\mathcal
H}_0$ and $B$ is an odd continuous map.
 Then each member of ${\mathcal M}_k$ is of  genus at least $k$
 (available in ${\mathcal H}_0$).
   It is also important to remind that the
definition of genus \cite[p.~385]{KrasZ} assumes  that
$\rho({\mathcal F})=1$, if no {\em component} of ${\mathcal F}
\cup {\mathcal F}^*$, where
 $
 {\mathcal F}^*=\{v: \,\, -v \in {\mathcal F}\},
 $
 is the {\em reflection} of ${\mathcal F}$ relative to 0,
 contains a pair of antipodal points $v$ and $v^*=-v$.
 Furthermore, $\rho({\mathcal F})=n$ if each compact subset of
${\mathcal F}$ can be covered by, minimum, $n$ sets of genus one.

According to (\ref{ck1}),
 $
 c_1 \le c_2 \le ... \le c_{l_0},
 $
 where $l_0=l_0(R)$ is the category of ${\mathcal H}_0$ (see an estimate
 below) satisfying
  \be
  \label{l01}
  l_0(R) \to + \infty \quad \mbox{as} \quad R \to \infty.
  \ee
  Roughly speaking,
since the dimension of the sets ${\mathcal F}$ involved in the
construction of ${\mathcal M}_k$ increases with $k$, this
guarantees that the critical points delivering critical values
(\ref{ck1}) are all different.

  It follows from (\ref{f2}) that the category
$l_0=\rho({\mathcal H}_0)$  of the set ${\mathcal H}_0$ is equal
to the number (with multiplicities) of the eigenvalues $\l_k>-1$
of the linear bi-harmonic operator,
 \be
 \label{f55}
 - w^{(4)}= \l_k \psi, \quad \psi \in W^2_{2,0}(B_R).
  \ee
  Since the dependence of the spectrum on $R$ is, obviously,
   \be
   \label{f56}
   \l_k(R)= R^{-4} \l_k(1), \quad k=0,1,2,... \, ,
    \ee
 the category $\rho({\mathcal H}_0)$ can be arbitrarily large for
$R \gg 1$,  (\ref{l01}) holds, and we obtain:

\begin{proposition}
 \label{Pr.MM}
The ODE problem $(\ref{1})$ has at least a countable set of
different solutions denoted by $\{F_l, \, l \ge 0\}$, each one
obtained as a critical point of the functional $(\ref{V1})$ in
$W^2_{2,0}(B_R)$ with sufficiently large $R>0$.
 \end{proposition}



 \subsection{First basic pattern and local structure of zeros}

\begin{figure}
 \centering
\includegraphics[scale=0.8]{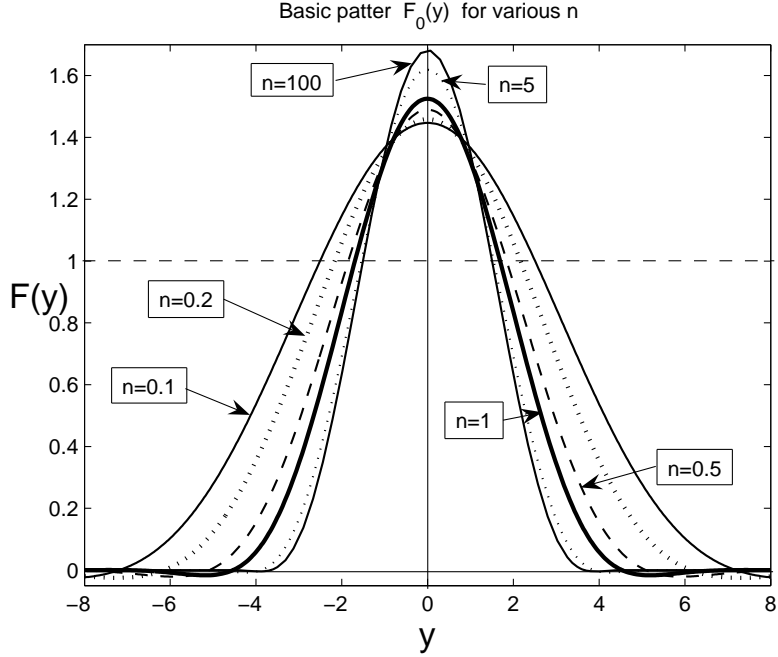}  
 \vskip -.4cm
\caption{\rm\small The first solution $F_0(y)$  of  (\ref{1}) for
various $n$.}
   \vskip -.3cm
 \label{G1}
\end{figure}

 We next present numerical results concerning existence
 and multiplicity of solutions for equation (\ref{1}).
 In Figure
\ref{G1}, we show the first basic 
pattern for (\ref{1}) called the  $F_0(y)$ for various $n \in [
0.1,100]$. Note that (\ref{1}) admits a natural passage to the
limit $n \to +\infty$ that gives the ODE
 \be
 \label{sign1}
 F^{(4)}= F- {\rm sign} \, F \equiv \left\{
 \begin{matrix}
 F-1 \quad\,\,\, \mbox{for} \quad F  \ge 0, \ssk\\
 F+1 \quad \mbox{for} \quad F <0.
  \end{matrix}
 \right.
  \ee
A unique oscillatory solution of (\ref{sign1}) can be treated by
an algebraic approach; cf. \cite[\S~7.4]{Gl4}. The solutions for
$n=1000$ or $n=+\infty$ practically do not differ from the last
profile for $n=100$ in Figure \ref{G1}.

\begin{figure}
\centering \subfigure[scale $10^{-3}$]{
\includegraphics[scale=0.52]{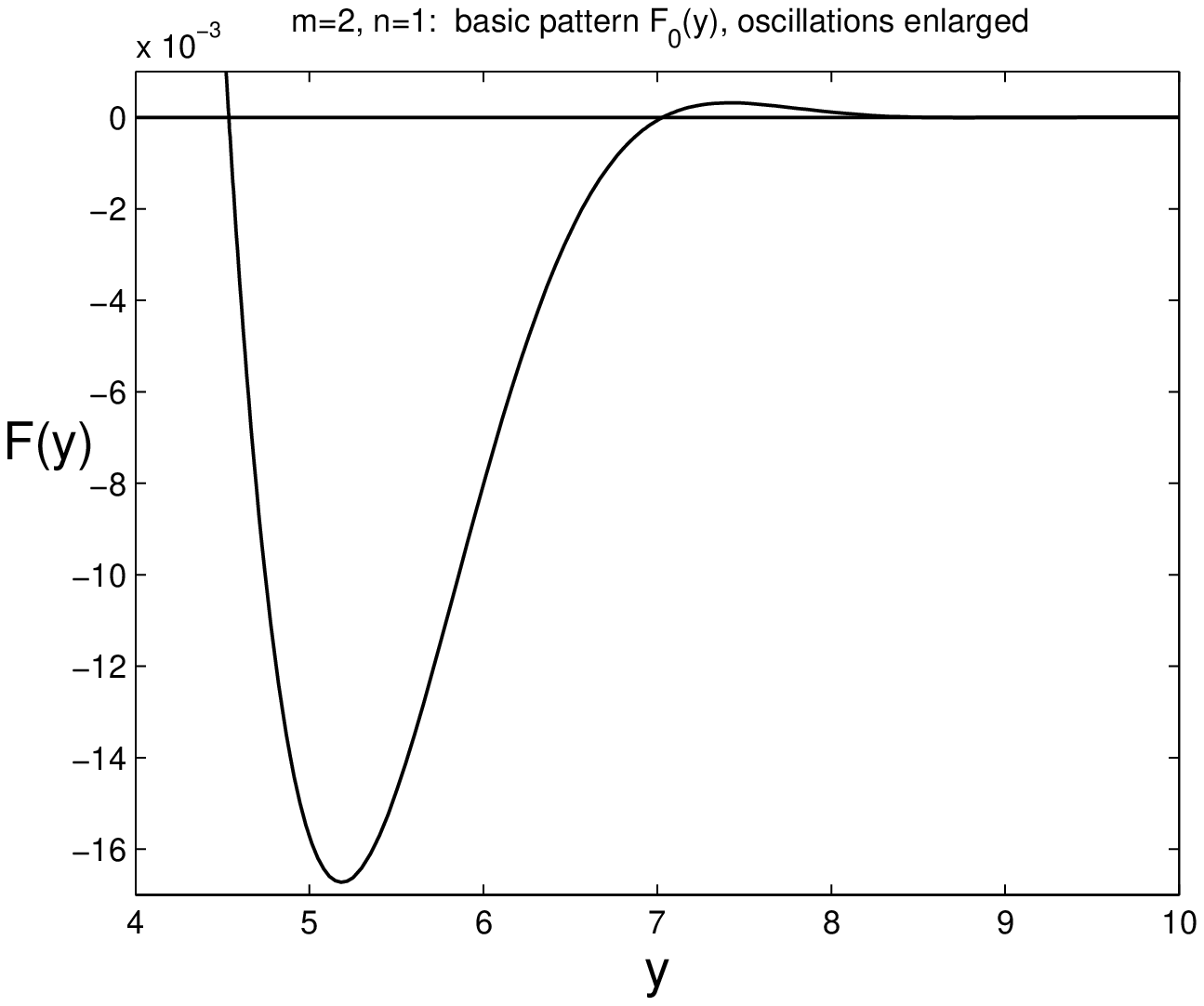}
} \subfigure[scale $10^{-4}$]{
\includegraphics[scale=0.52]{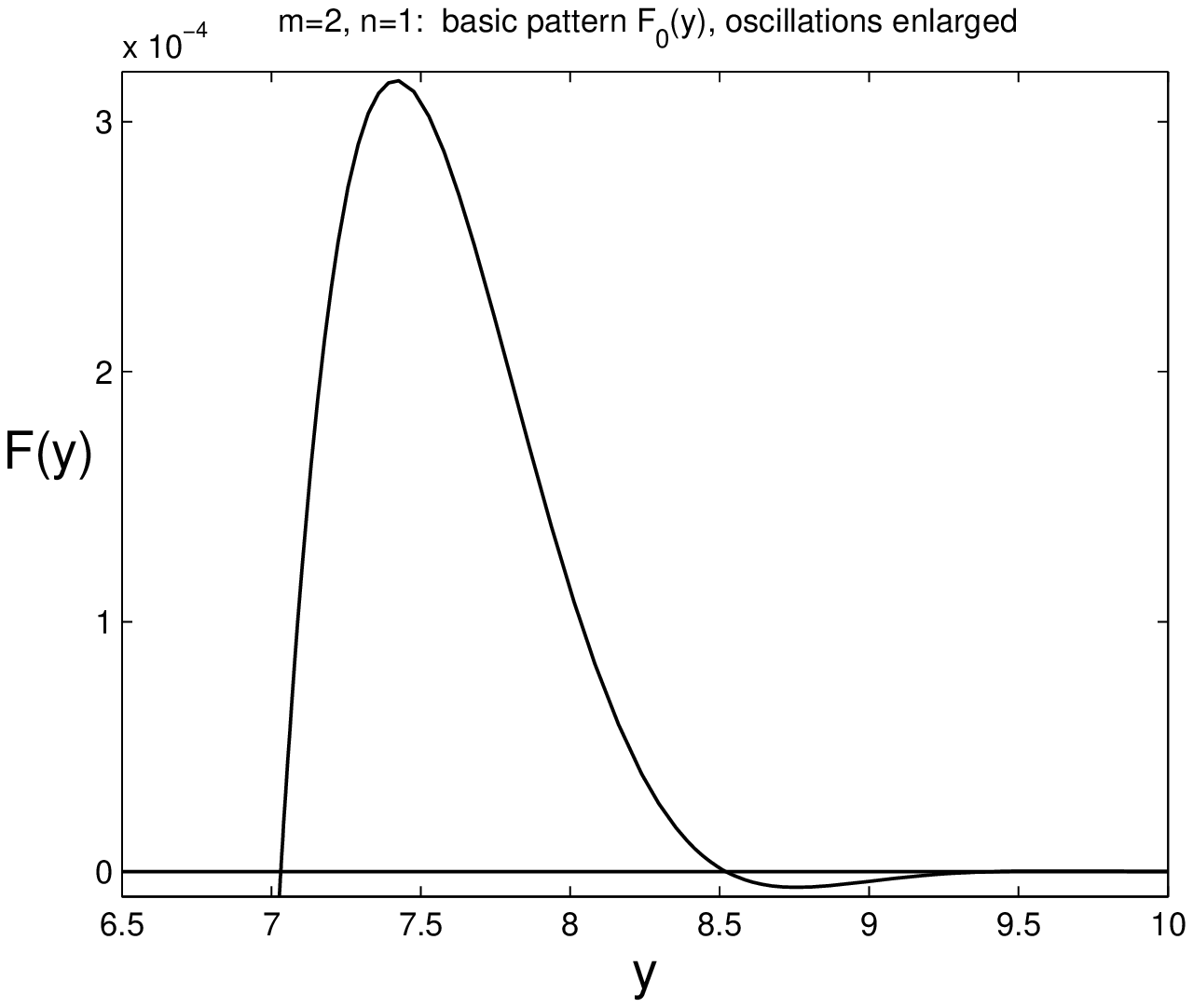}               
}
 \vskip -.2cm
\caption{\rm\small  Enlarged zero structure of the profile
$F_0(y)$ for $n=1$ in  the linear scale.}
 \label{G2}
\end{figure}

 These profiles are  constructed by {\tt MatLab}
 by using a natural regularization in the singular term in (\ref{1}),
 \be
 \label{4.1}
  \mbox{$
 F^{(4)}=F- 
 \big(\e^2+ F^2 \big)^{-\frac n{2(n+1)}} F \quad \mbox{in} \quad
 \re.
 $}
 \ee
 Here,  the regularization parameter $\e$
 and both absolute and relative tolerances in the {\tt bvp4c}
solver have been enhanced and took the values
 \be
 \label{zz0}
  \e= 10^{-10} \quad  \mbox{and} \quad  {\rm Tols}=10^{-10}.
 \ee
 This allows us also to reveal the refined local structure of
 multiple zeros at the interfaces. Figure \ref{G2} for $n=1$ shows how
the zero structure repeats itself in a self-similar manner from
one zero to another in the usual linear scale.

\begin{figure}
 \centering
\includegraphics[scale=0.7]{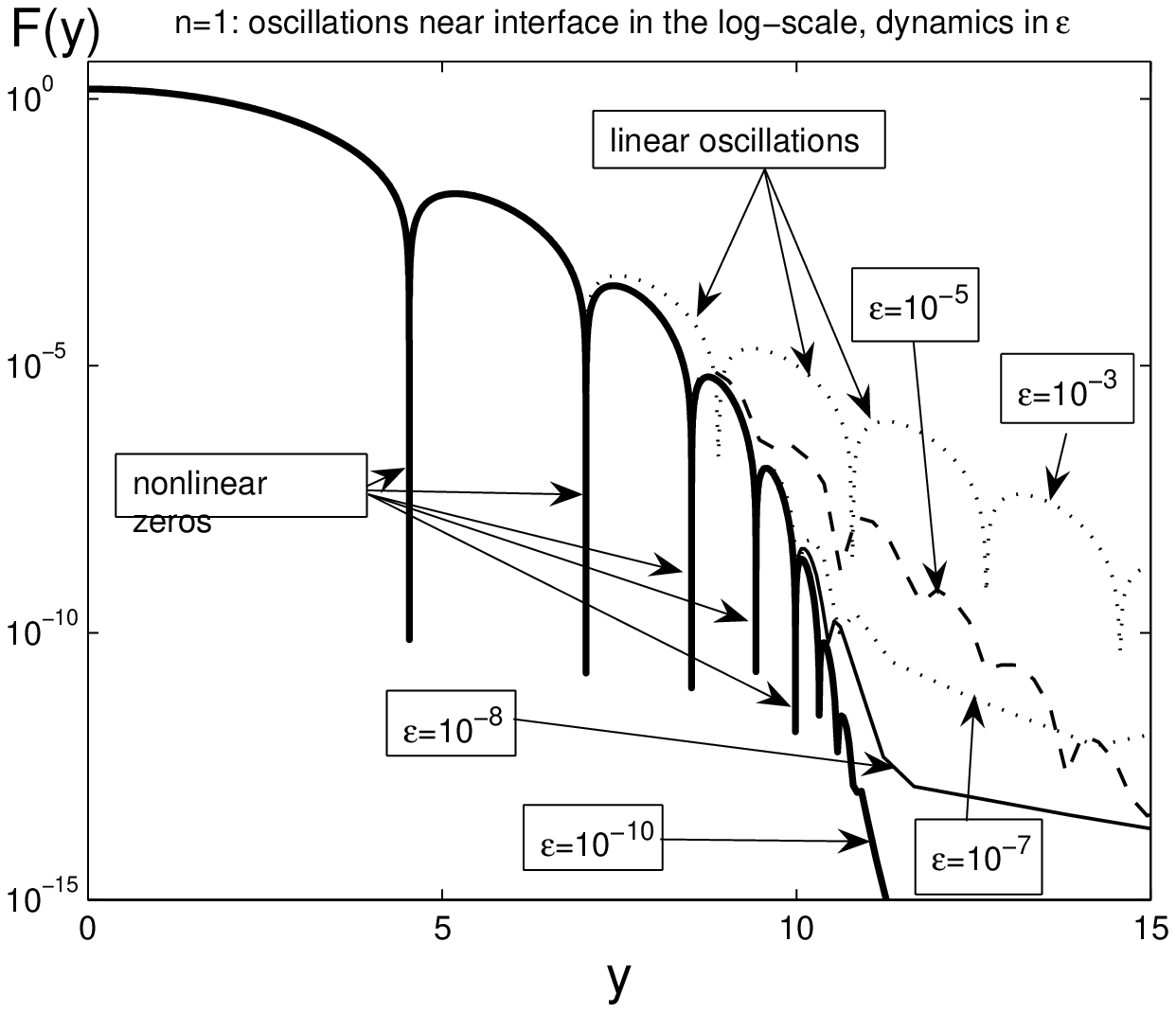} 
 \vskip -.4cm
\caption{\rm\small Behaviour of $F_0(y)$  for
 $n=1$ in the log-scale near the interface: formation of five nonlinear zeros as $\e={\rm Tols}$
 decrease from $10^{-3}$ to $10^{-7}$.}
   \vskip -.3cm
 \label{ZZ1}
\end{figure}

 In Figure \ref{ZZ1}, we
present the results for $n=1$ that show the oscillatory structure
in the log-scale such as (\ref{le3}) with $y \mapsto y_0-y$ and
 \be
 \label{zz1}
  \mbox{$
 \mu=  \frac{4(n+1)}n \big|_{n=1}=8 \quad \Longrightarrow
 \quad \ln |F(y)|= 8 \ln(y_0-y) + \ln |\varphi_*(\ln(y_0-y))|+...
  $}
  \ee
  (see  also (\ref{b3}) below).
 This figure shows an ``$\e$-dynamic" formation of at least six-seven nonlinear zeros, when we
 decrease $\e={\rm Tols}$ from $10^{-3}$ (just two first nonlinear
 zeros; the rest correspond to the linear ones in (\ref{zz2}))
 to $10^{-10}$ (6-7 zeros are nonlinear). Observe a ``concave"
 shape of the graph for the last $\e=10^{-10}$,
 which is consistent with the log-shape in (\ref{zz1})
(we claim that even a trace of the multiplier 8 can be
distinguished by 6 zeros).

 Further decrease  of $\e$ and Tols leads to quick divergence
 of {\tt bvp4c} solver, which has the limit minimal Tols
 $10^{-13}$, which is not that helpful in comparison with
 $10^{-10}$ achieved in Figure \ref{ZZ1}.
 The right-hand interface is then estimated as (it is quite a
challenge  to detect numerically the free-boundary point more
accurately)
 $$
 y_0 \approx 12.
  $$
 The ``nonlinear area" ends at $y \sim 11$ (due to the
 requires accuracy and $\e$-regularization), and next we observe
 the ``linearized area" where (\ref{4.1}) implies an exponential
 behaviour for $y \gg 1$ governed by the linearized ODE
   \be
   \label{zz2}
   \mbox{$
 F^{(4)}=- 
 \e^{-\frac n{n+1}} F+...
   \quad \Longrightarrow \quad
   F(y) \sim {\mathrm e}^{- c_1 \e^{-  \frac n{4(n+1)}}} \,
    \cos\big(c_2  \e^{- \frac n{4(n+1)}}y+c_3 \big),
     $}
  \ee
  where $c_k$ are constants and $c_1>0$.
 From $y \sim 12.5$ ($F \sim 10^{-12}$), the numerics are not
 reliable at all.

\subsection{The first basic pattern $F_0(y)$:  on the behaviour as $n \to 0$}

\begin{figure}
\centering \subfigure[profiles $F_0(y)$]{
\includegraphics[scale=0.52]{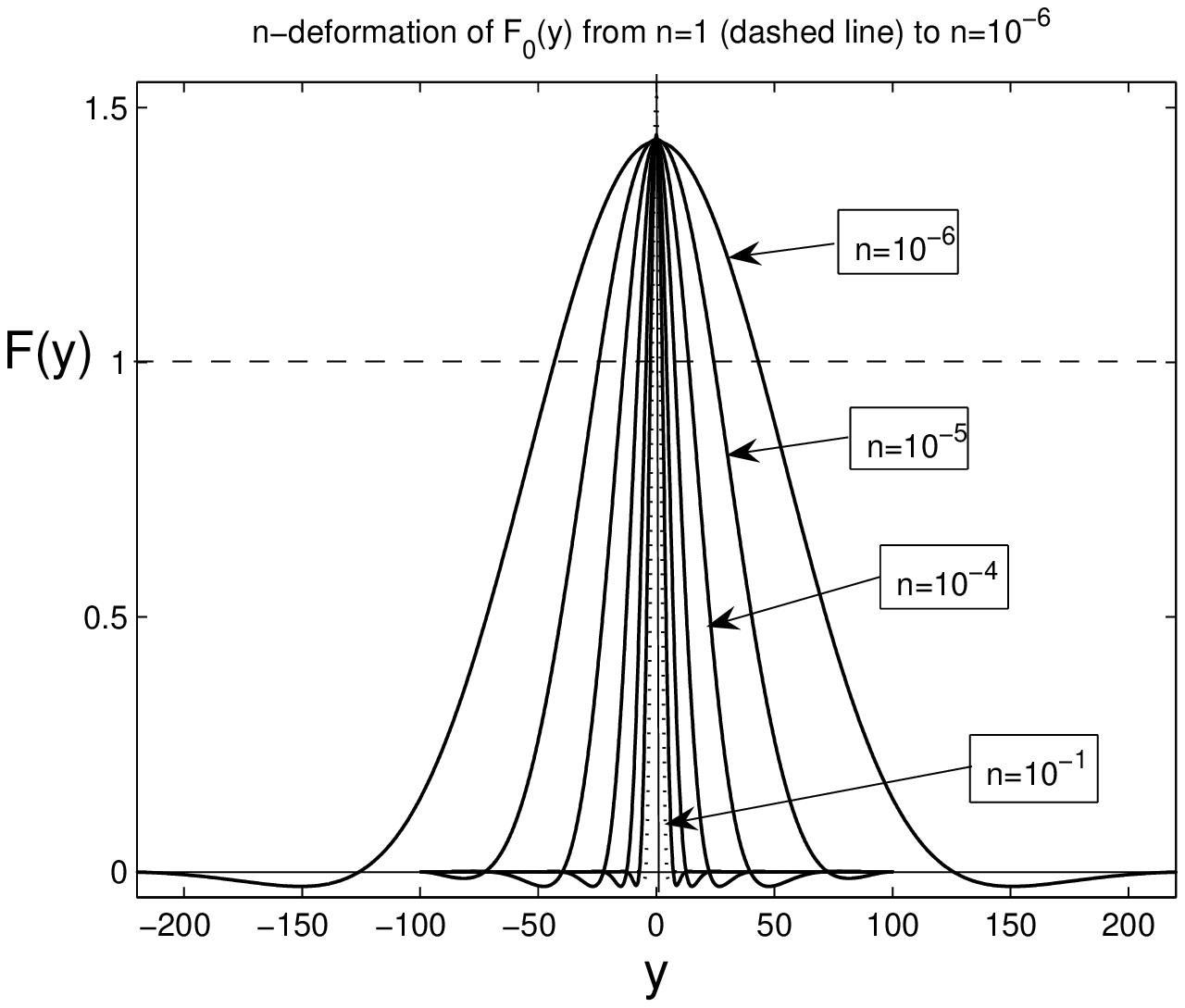}
} \subfigure[enlarged for $y \approx 0$]{
\includegraphics[scale=0.52]{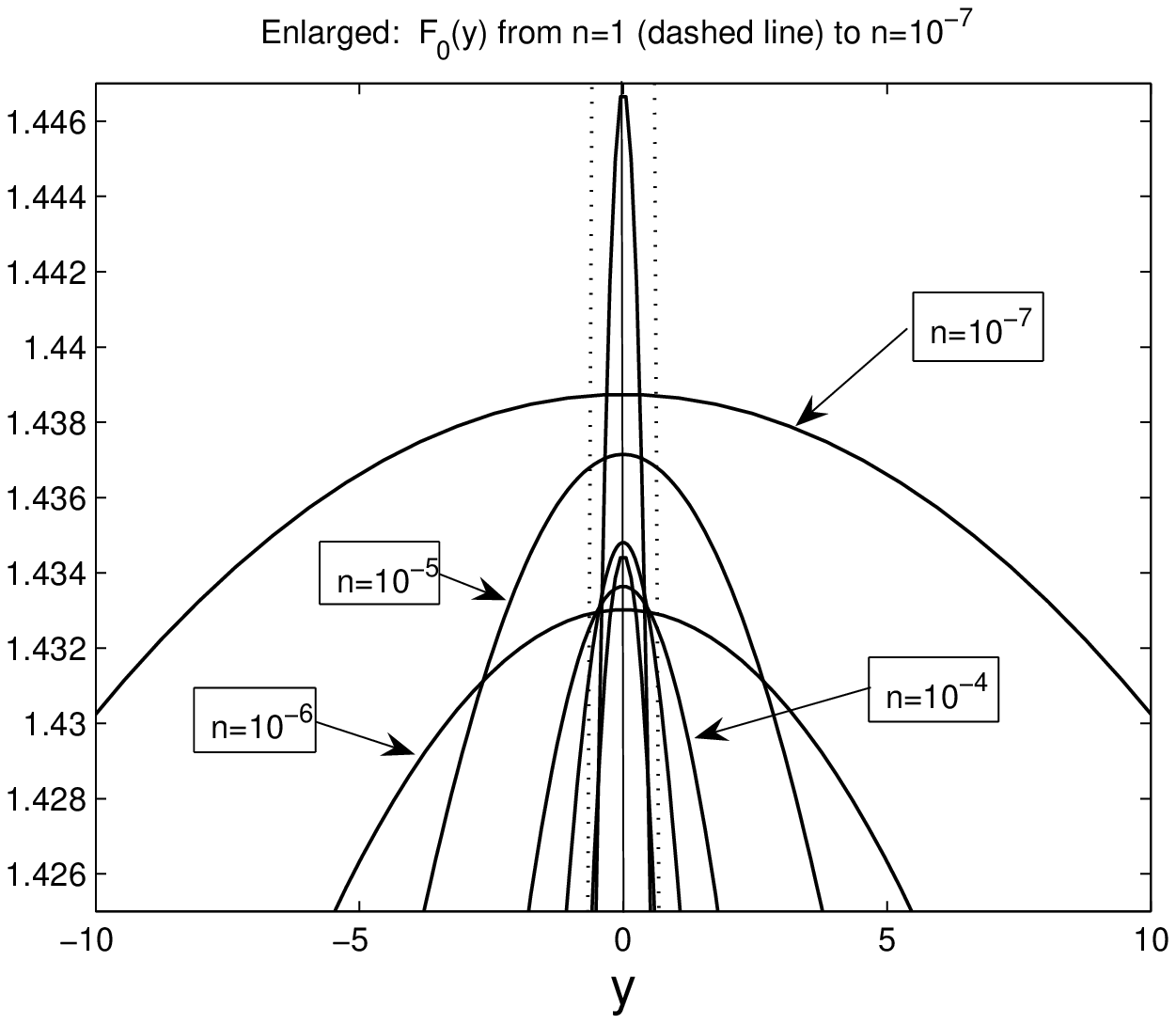}               
}
 \vskip -.2cm
\caption{\rm\small  Deformation of the similarity solutions
$F_0(y)$ of the ODE (\ref{1}) for small $n>0$.}
 \label{Fn01}
\end{figure}

      In order to
understand the continuous connection with the linear PDE
(\ref{bh1}) and other homotopy results, one needs to pass to the
limit $n \to 0^+$ in the ODE (\ref{1}). Figure \ref{Fn01} shows
the behaviour of the first basic profile $F_0(y)$ for small
$n=10^{-1}$, $10^{-2}$, ... , $10^{-7}$ (a), while (b) explains
the behaviour of the maximum value $F_0(0)$ and around.

It follows from equation (\ref{1}) that, for $n \approx 0^+$,
close to the origin and uniformly on bounded intervals in $y$, the
smooth solution $F_0(y)$ approaches an even function satisfying
the linear homogeneous equation
 \be
 \label{b1}
 F^{(4)}=0 \quad \Longrightarrow \quad
 F_0(y) \approx C_0(n) - A_0(n) y^2,
  \ee
  so that $F_0(y)$ essentially depends on small $n$.
Here both unknown coefficients $C_0(n) \to C_0(0)$ and $A_0(n) \to
0$ are slightly oscillatory as $n \to 0$ (in view of the
non-monotone  behaviour in (\ref{b3}) near the interface), with
 \be
 \label{b2}
 C_0(0) \approx 1.435... \, .
  \ee
The behaviour (\ref{b1}) is then to be matched with the
asymptotics near the interface at some $y=y_0^-(n) \gg 1$ (cf.
(\ref{le3})) that is governed by the non-Lipschitz term in
(\ref{1}),
 \be
 \label{b3}
 F^{(4)}= - \big|F \big|^{-\frac n{n+1}} F+...
 \quad \Longrightarrow \quad F(y)=\big(y_0-y\big)^{\frac {4(n+1)} n} \varphi_*\big(\ln
 (y_0-y)\big)+...\, ,
  \ee
  where $\varphi_*$ is the corresponding periodic solution of the
  ODE for the oscillatory component. A dimensional
  analysis implies that the interface position can be estimated as
  follows:
   \be
   \label{b4}
   y_0(n) \sim n^{-\frac 34} \quad \mbox{as}
   \quad n \to 0^+;
    \ee
    see \cite[\S~9.6]{Gl4} for details. First rough matching of (\ref{b1})
    and (\ref{b3}), (\ref{b4}) yields
    $$
    A_0(n) \sim n^{\frac 32}.
     $$
A more accurate matching is difficult and typically  assumes
appearing logarithmic, $\ln n$, and other terms including
oscillatory ones; see   \cite[\S~7]{Gl4}.

  \subsection{Basic countable family}
 \label{Sect54}

\begin{figure}
\centering \subfigure[ $F_0(y)$ ]{
\includegraphics[scale=0.52]{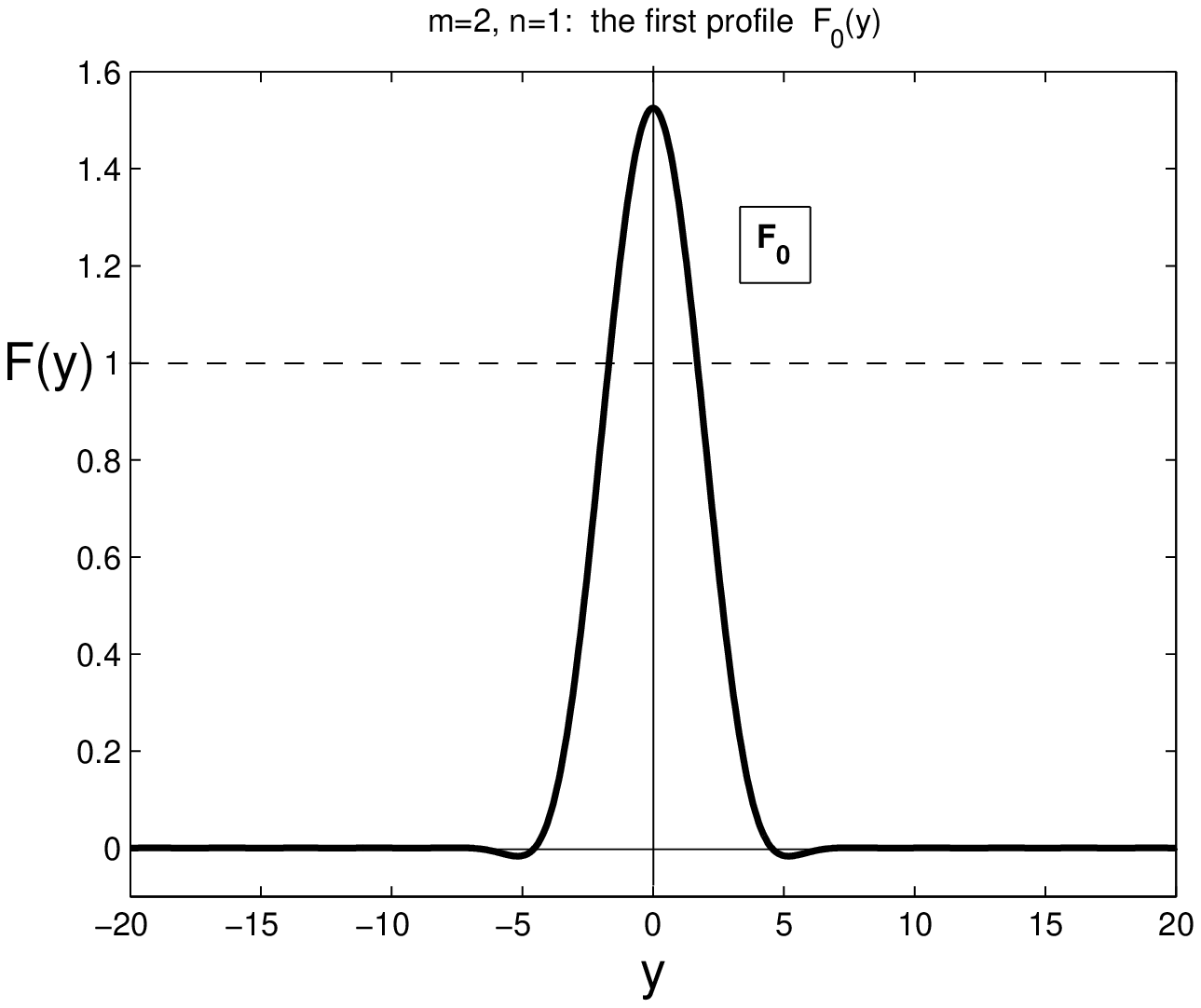}
} \subfigure[$F_1(y)$]{
\includegraphics[scale=0.52]{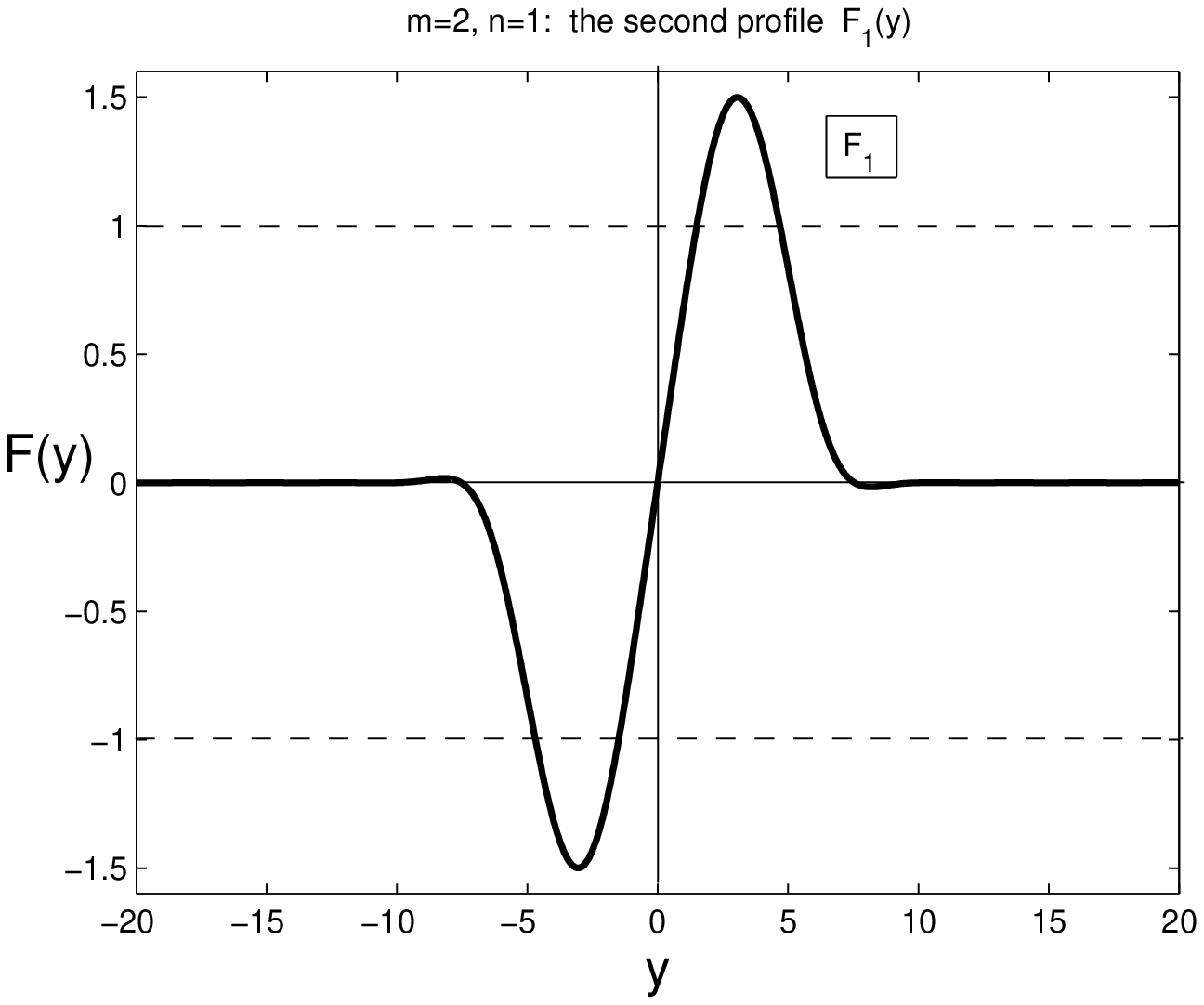}
} \subfigure[$F_2(y)$]{
\includegraphics[scale=0.52]{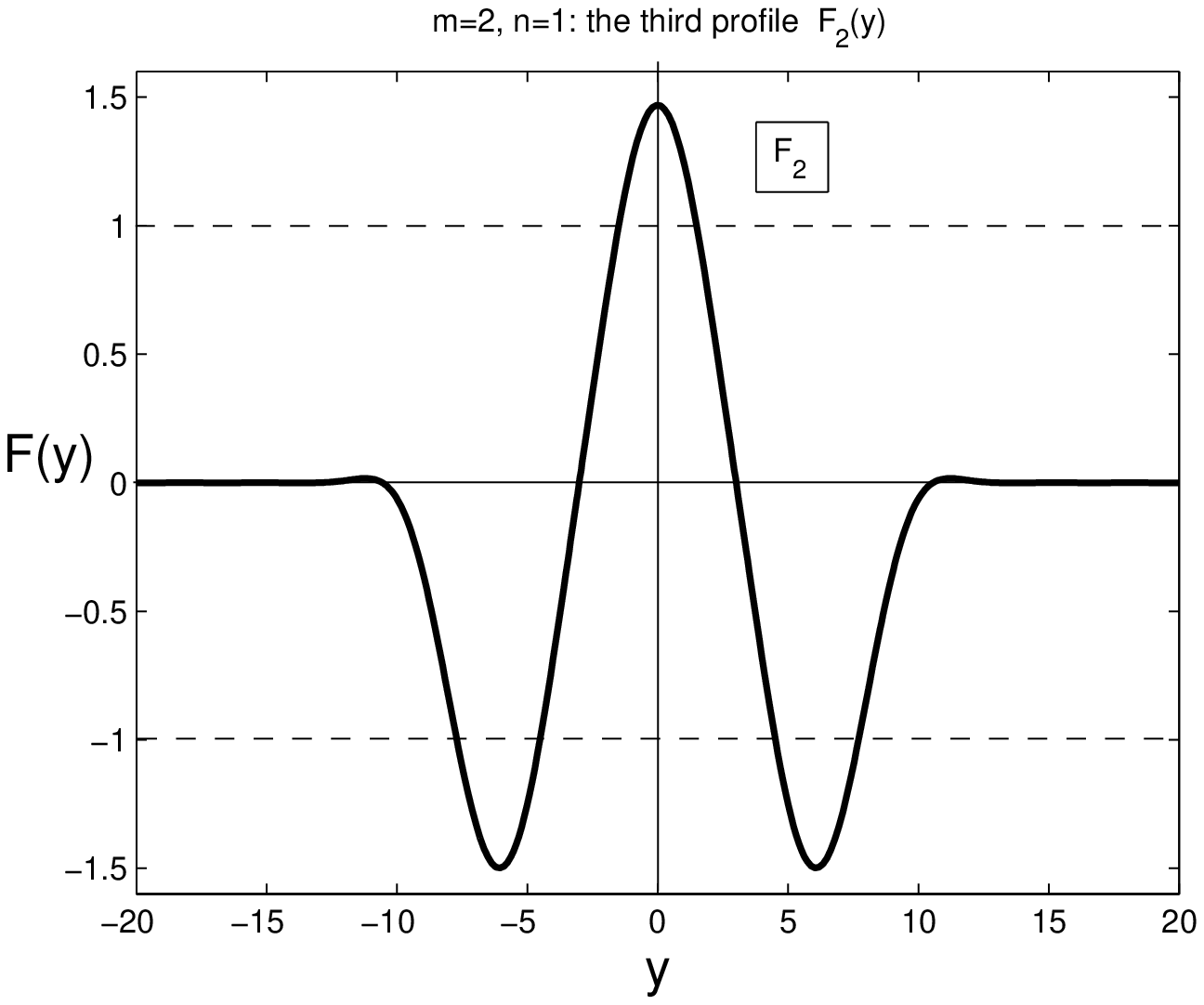}
} \subfigure[$F_3(y)$]{
\includegraphics[scale=0.52]{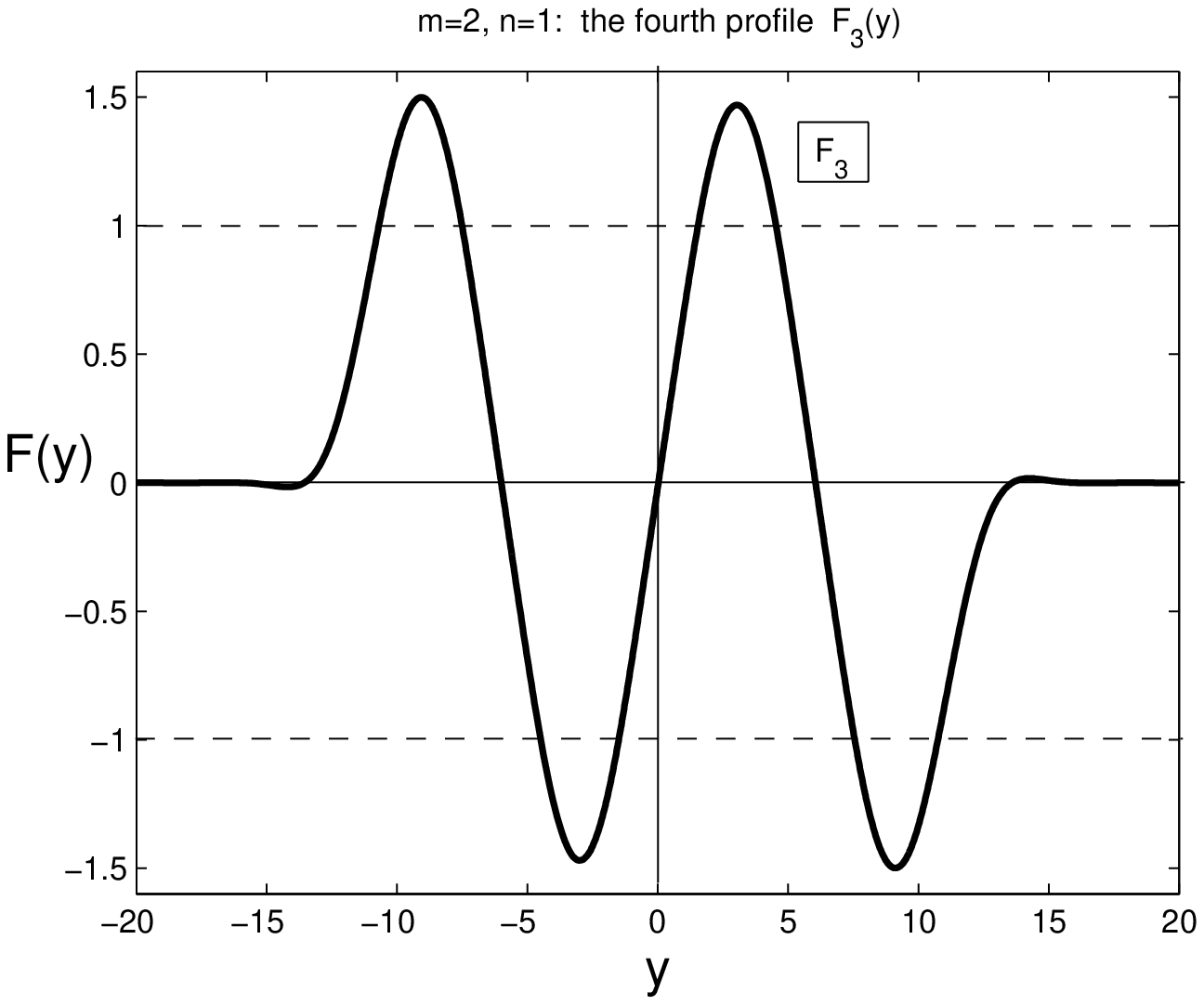}
} \subfigure[$F_4(y)$]{
\includegraphics[scale=0.52]{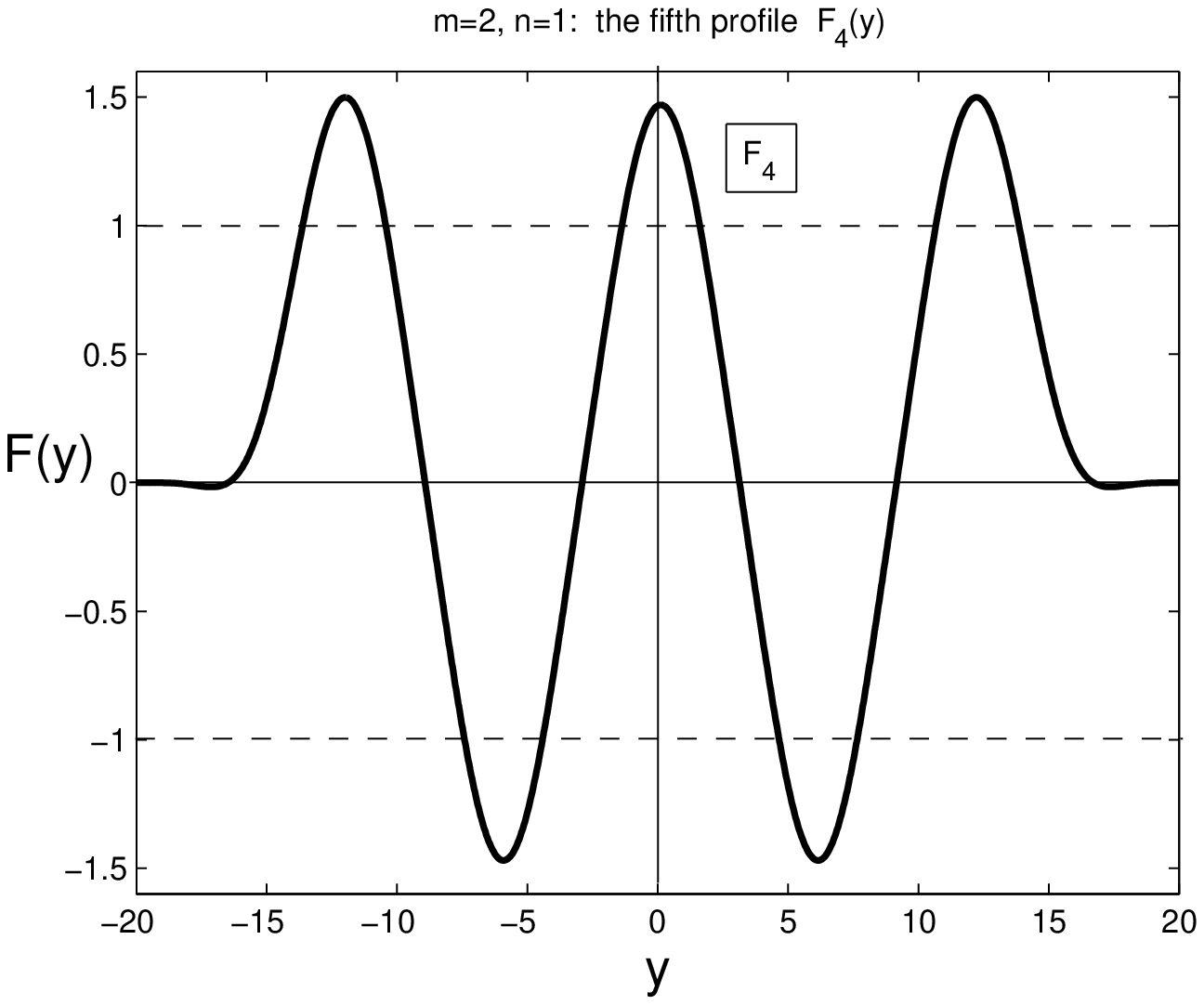}
} \subfigure[$F_5(y)$]{
\includegraphics[scale=0.52]{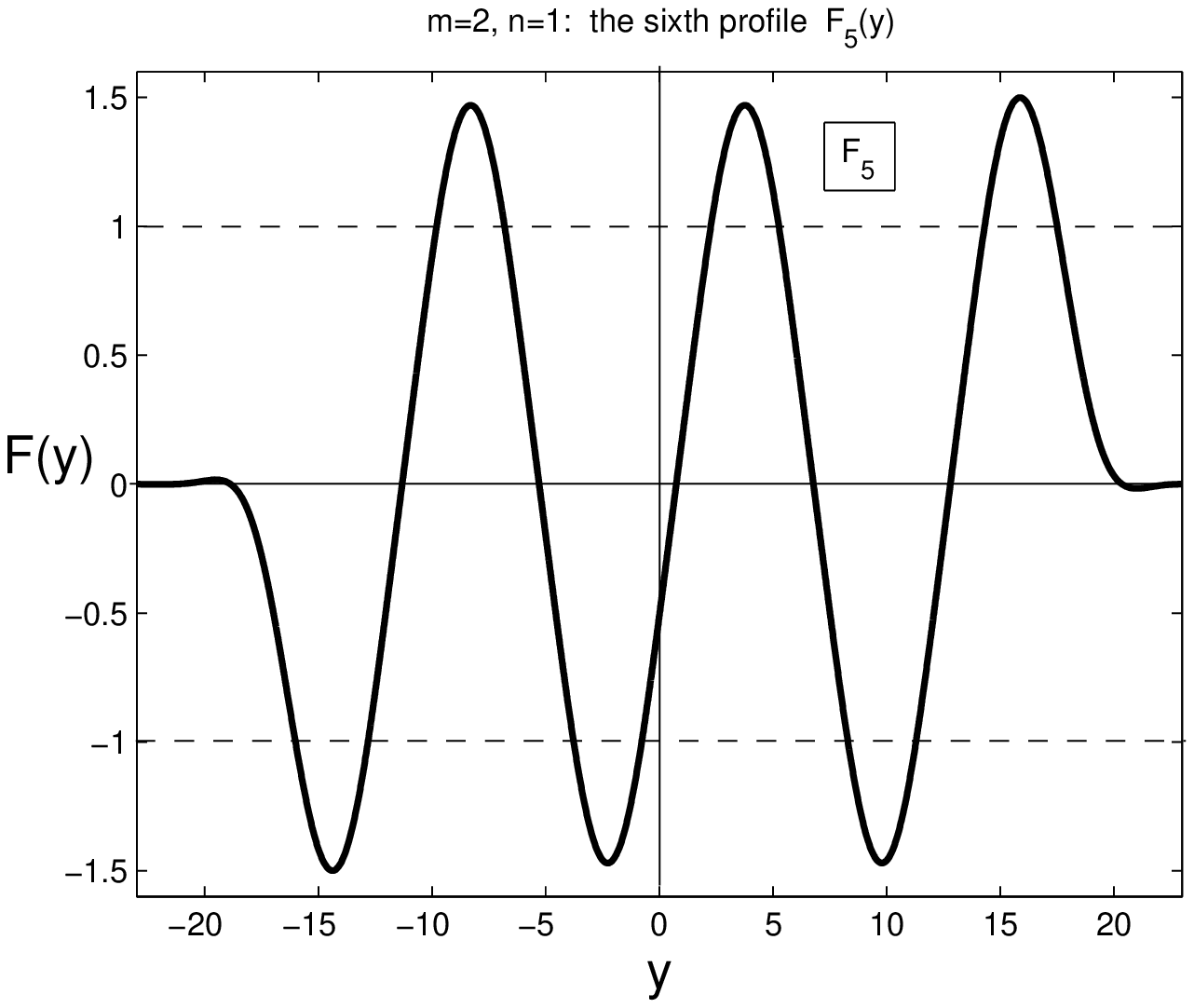}
}
 \vskip -.2cm
\caption{\rm\small The first six patterns of the basic family
$\{F_l\}$ of the ODE (\ref{1}) for $n=1$.}
 \label{G4}
\end{figure}

  Figure \ref{G4} shows\footnote{Amazingly, these profiles
   look almost identical to those for the PDE (\ref{pLap1}) \cite[\S~4]{GalpLap},
   corresponding
   to a completely different diffusion operator and much harder and distinct L--S theory.
   These similarities underline some quite obscure and involved common evolution
    properties of blow-up for different PDEs.}
    the basic family denoted by
 $
 \{F_l, \, l=0,1,2,...\}
  $
   of solutions of (\ref{1}) for $n=1$. This
   is
   associated with the application of  L--S and fibering theory,
   \cite{GPHom}.
Each profile $F_l(y)$ has  $l+1$ ``dominant" extrema and $l$
``transversal" zeros; see \cite[\S~5]{GPHom} and
\cite[\S~4]{GHUni} for further details. We claim that
 $$
\mbox{all the internal zeros of $F_l(y)$ are {transversal}}
 $$
(excluding the oscillatory end points of the support), and this is
a fundamental conclusion of our analysis, which deserves more
rigorous treatment.
  In other words, each profile $F_l$ is
approximately obtained by a simple ``interaction" (gluing
together) of $l+1$ copies of the first pattern $\pm F_0$ taking
with necessary signs; see further comments below.

Actually, if we forget for a moment about the complicated
oscillatory structure of solutions near interfaces, where an
infinite number of extrema and zeros occur, the dominant geometry
of profiles in Figure \ref{G4}  approximately obeys Sturm's
classic zero set property, which is true rigorously for
 the second-order ODE only,
 \be
 \label{4.4}
  \mbox{$
 F''=-F + \big|F \big|^{-\frac n{n+1}}F \quad \mbox{in} \quad \re.
 $}
  \ee
For (\ref{4.4}), the basic family $\{F_l\}$ is constructed by
direct gluing together simple patterns $\pm F_0$ given explicitly;
see \cite[p.~168]{GSVR}. Therefore, each $F_l$ consists of  $l+1$
patterns  (with signs $\pm F_0$), so that Sturm's property is
clearly true.



 \subsection{Countable family of $\{F_0,F_0\}$-gluing and extensions}

This procedure is similar to that in \cite[\S~4.7]{GalpLap}, but
we again recall that these are completely different variational
problems.
 Further patterns to be introduced do not
exhibit as clear a ``dominated" Sturm property and are associated
with a double fibering technique where both Cartesian and
spherical representations are involved; see details
\cite[\S~6]{GPHom}. Let us present basic explanations.

\begin{figure}
 \centering
\includegraphics[scale=0.65]{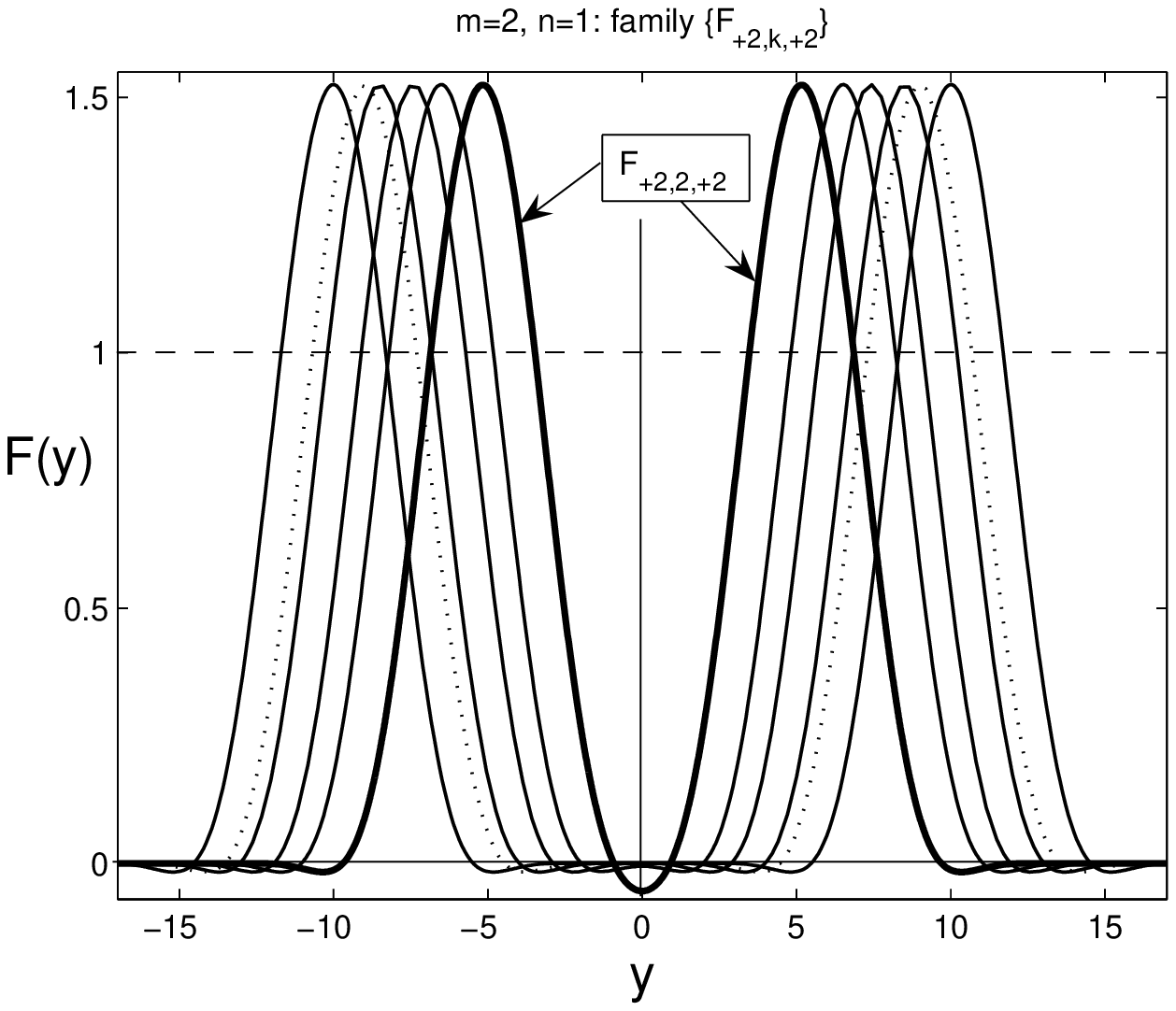}
 \vskip -.4cm
\caption{\rm\small First  patterns from the family
$\{F_{+2,k,+2}\}$ of the  $\{F_0,F_0\}$-gluing;  $n=1$.}
   \vskip -.3cm
 \label{G6}
\end{figure}

The nonlinear interaction of the two first
 patterns $F_0(y)$ leads to a new family of profiles.
In Figure \ref{G6} for $n=1$, we show the first profiles from this
family denoted by $\{F_{+2,k,+2}\}$, where in each function
$F_{+2,k,+2}$ the multiindex $\s=\{+2,k,+2\}$ means, from left to
right, +2 intersections with the equilibrium +1, then next $k$
intersections with zero, and final +2 stands again for 2
intersections with +1. Later on, we will use such a multiindex
notation to classify other  patterns obtained.

Since $F_0(y)$ is infinitely oscillatory at the interfaces, the
family $\{F_{+2,k,+2}\}$ is expected to be countable, so such
functions exist for any even $k=0,2,4,...\, $. Therefore,
$k=+\infty$ leads to the non-interacting pair with no overlapping
of supports,
  \be
 \label{F0+}
 F_0(y+ y_0) + F_0(y-y_0), \quad \mbox{where} \,\,\,\,\,{\rm supp}\,
 F_0(y) = [-y_0,y_0].
  \ee
There exist various triple $\{F_0,F_0,F_0\}$ and any multiple
interactions $\{F_0,...,F_0\}$ of arbitrarily large $k$ single
profiles, with a variety of distributions of zeros between any
pair of neighbours.

 \subsection{Countable family of  $\{-F_0,F_0\}$-gluing and extensions}

\begin{figure}
 \centering
\includegraphics[scale=0.65]{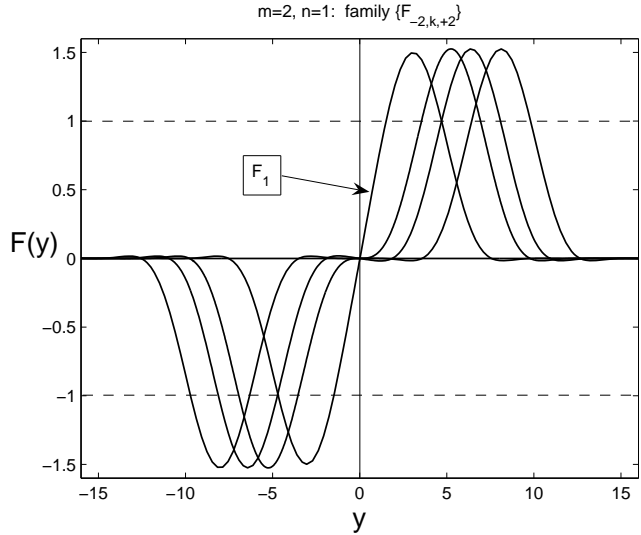}
 \vskip -.4cm
\caption{\rm\small First four patterns from the family
$\{F_{-2,k,+2}\}$ of the  $\{-F_0,F_0\}$-gluing;  $n=1$.}
   \vskip -.3cm
 \label{G7}
\end{figure}

Consider next the interaction of $-F_0(y)$ with
  $F_0(y)$.
In Figure \ref{G7}, for $n=1$, we show the first profiles from
this
family denoted by $\{F_{-2,k,+2}\}$, where, for  
 the
multiindex  $\s=\{-2,k,+2\}$,
  the first number $-2$ means 2
intersections with the equilibrium $-1$, etc.  The enlarged zero
structure shows that the first two profiles belong to the same
class
 $
 F_{-2,1,2},
 $
i.e., both have a single zero for $y \approx 0$.
  The last solution shown is $F_{-2,5,+2}$.
 This family $\{F_{-2,k,+2}\}$ is expected to be countable,
with profiles  existing for any odd $k=1,3,5,...$, where the pair
for $k=+\infty$ is non-interacting,
 $- F_0(y+ y_0) + F_0(y-y_0)$.
There can be constructed families of an arbitrary number of
interactions such as $\{\pm F_0, \pm F_0,..., \pm F_0\}$
consisting of any $k \ge 2$ members.

\subsection{Periodic solutions in $\re$}

\begin{figure}
 \centering
\includegraphics[scale=0.60]{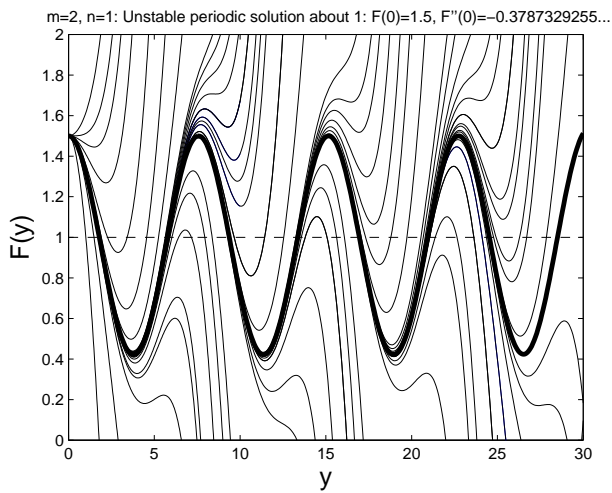}
 \vskip -.4cm
\caption{\rm\small An example of a periodic solution of the ODE
(\ref{1}) for $n=1$.}
   \vskip -.3cm
 \label{GP1}
\end{figure}

We are going to introduce new types of patterns, which exhibit a
different geometric shape. We then
 need first to
describe non-compactly supported periodic solutions in $\re$. As a
variational problem, equation (\ref{1}) admits an infinite number
of periodic solutions; see e.g. \cite[Ch.~8]{MitPoh}.  Figure
\ref{GP1} for $n=1$ presents a special  unstable (as $y \to +
\infty$) periodic solution obtained by shooting from the origin
with conditions
 $
 F(0)=1.5$, $F'(0)=F'''(0)=0$, and  $F''(0)=-0.3787329255... \,
 .
 $
 The periodic orbit $F_*(y)$ with the value
 \be
 \label{**1}
 F_*(0) \approx 1.535...
 \ee
 will be shown to play a key role in the
 construction of other families of compactly supported
patterns. Namely, this variety of solutions of (\ref{1}) having
oscillations about equilibria $\pm 1$ are  close to $\pm F_*(y)$
there.

\subsection{Family $\{F_{+2k}\}$}
 \label{S5.9}

\begin{figure}
\centering \subfigure[ $F_{+2k}(y)$ ]{
\includegraphics[scale=0.52]{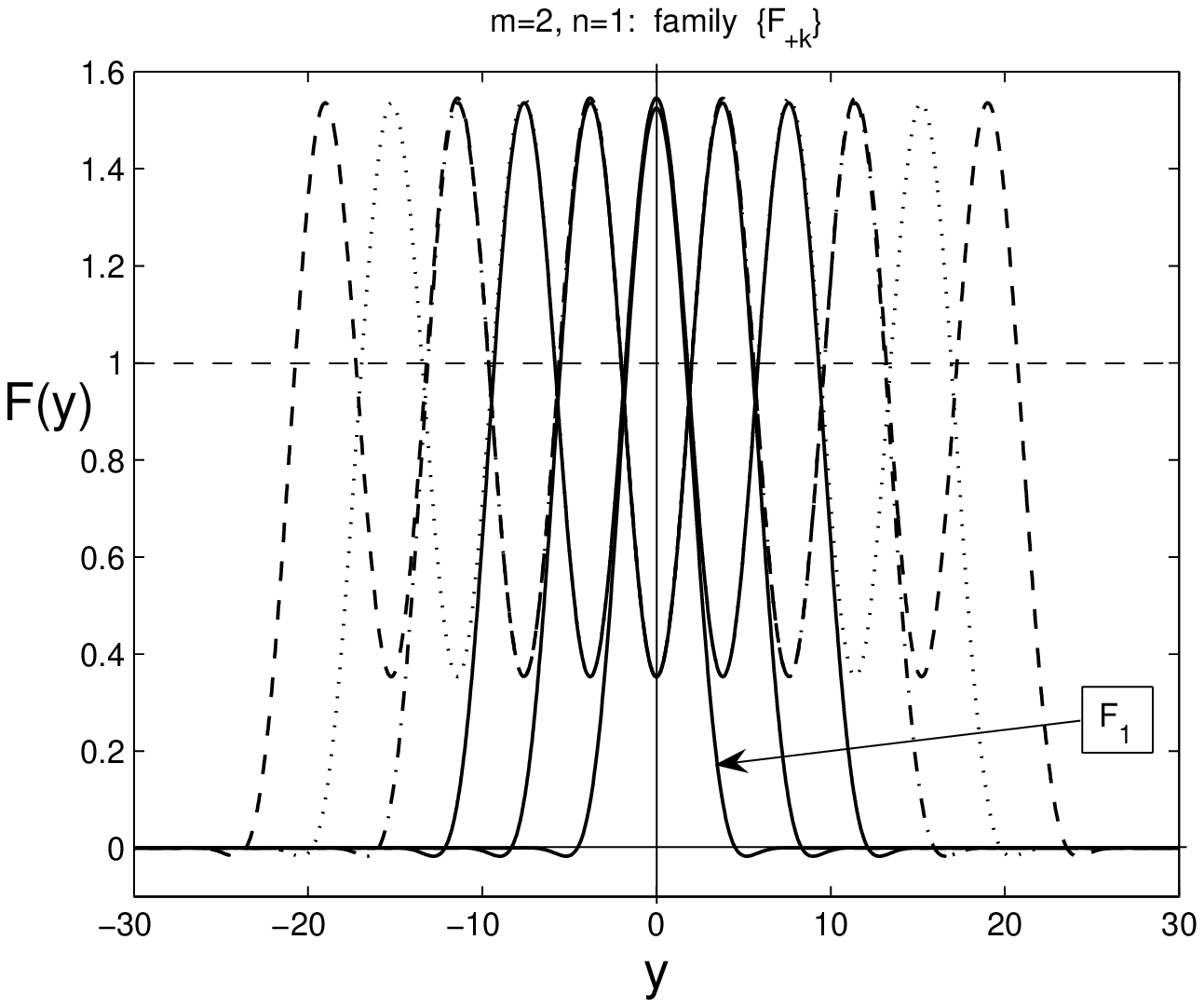}
} \subfigure[$F_{+k,l,-m,l,+k}$]{
\includegraphics[scale=0.52]{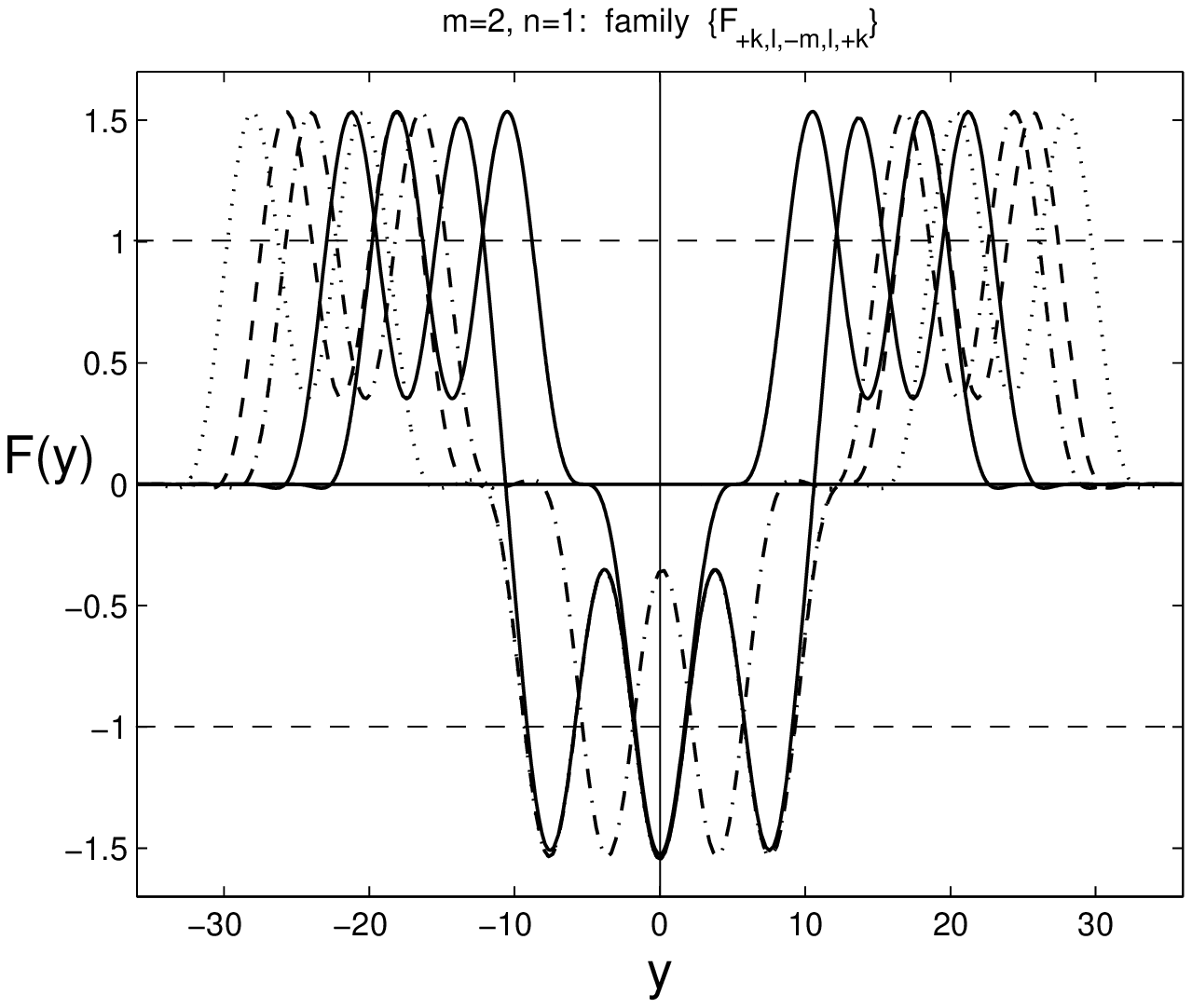}
}
 \vskip -.2cm
\caption{\rm\small Two families of solutions  of (\ref{1}) for
$n=1$; $F_{+2k}(y)$ (a) and $F_{+k,l,-m,l,+k}$ (b).}
 \label{G8}
\end{figure}

The patterns $F_{+2k}$ for $k \ge 1$ have $2k$ intersections with
the single equilibrium +1  only and have a clear ``almost"
periodic structure of oscillations about it; see Figure
\ref{G8}(a). The number of intersections with $F=+1$ denoted by
$+2k$ is an extra characterization of the Strum index to such a
pattern.
In this notation,
 $
 F_{+2}=F_0.
  $

\subsection{Complicated and chaotic patterns}

\begin{figure}
 \centering \subfigure[$\s=\{+2,2,+4,2,+2,1,-4\}$]{
\includegraphics[scale=0.52]{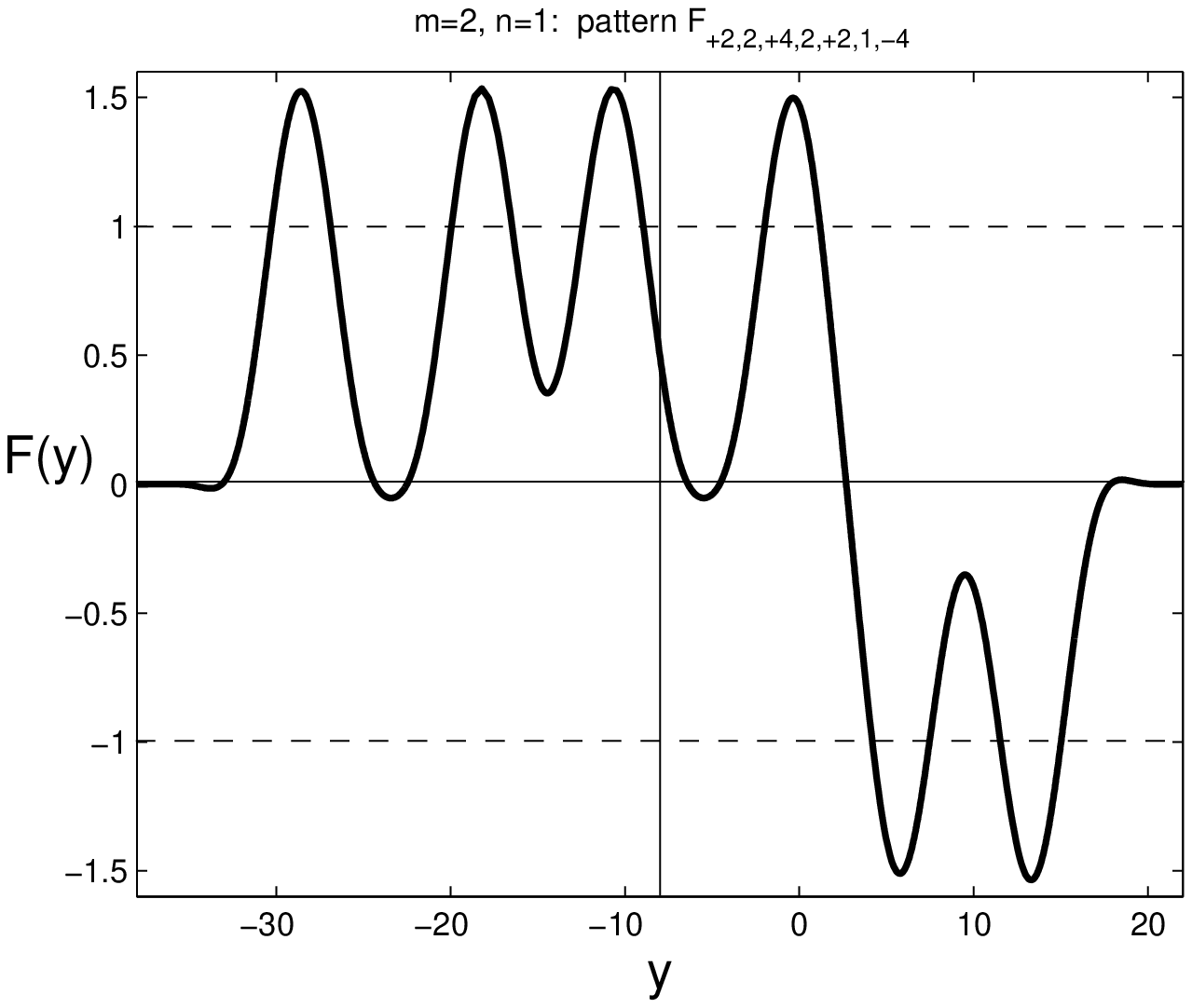}
} \centering \subfigure[$\s=\{+6,3,-4,2,-6\}$]{
\includegraphics[scale=0.52]{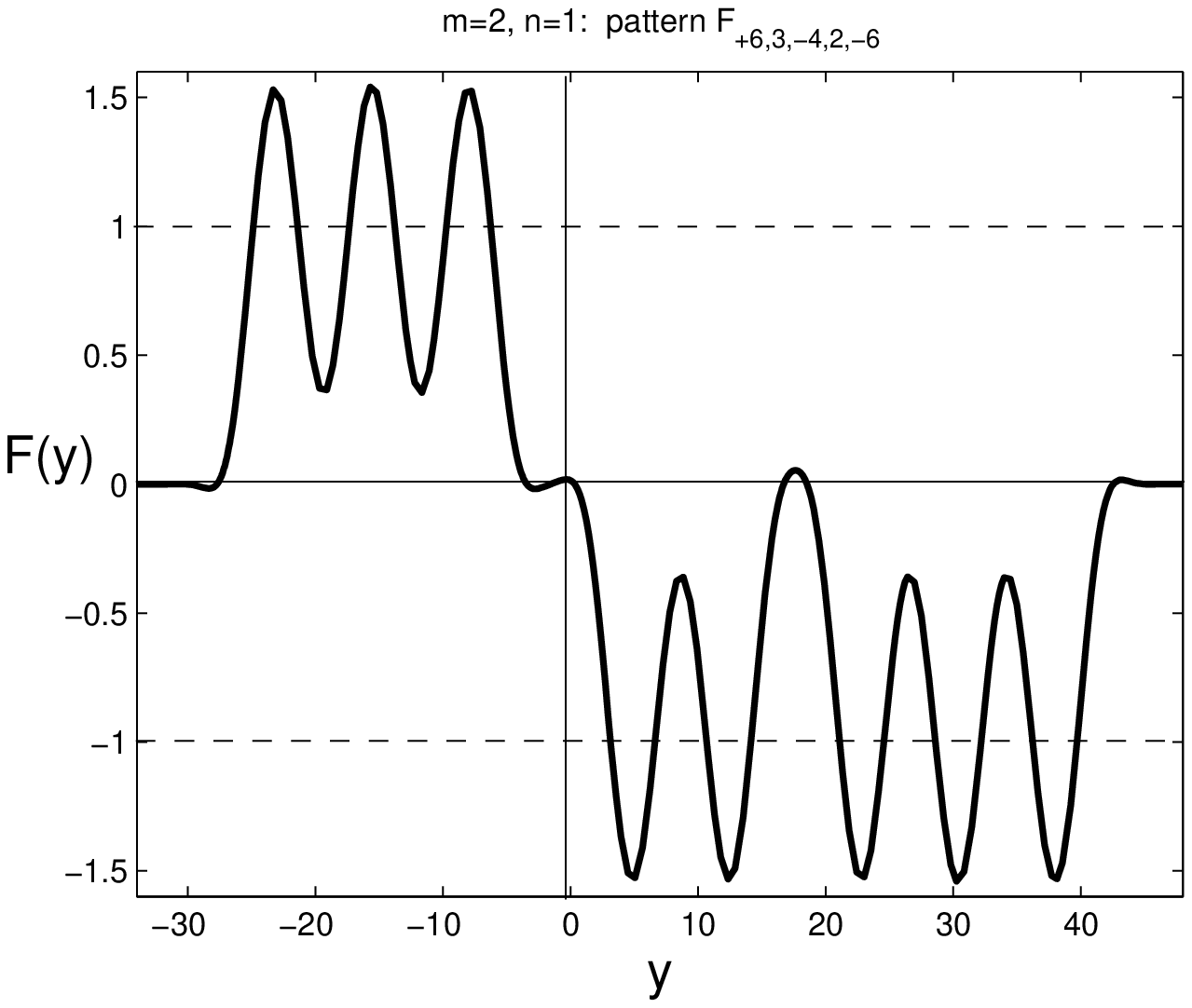}
}
 \vskip -.2cm
\caption{\rm\small  Two patterns for  (\ref{1}) for $n=1$.}
 \label{FC1}
\end{figure}

\begin{figure}
 \centering
\includegraphics[scale=0.65]{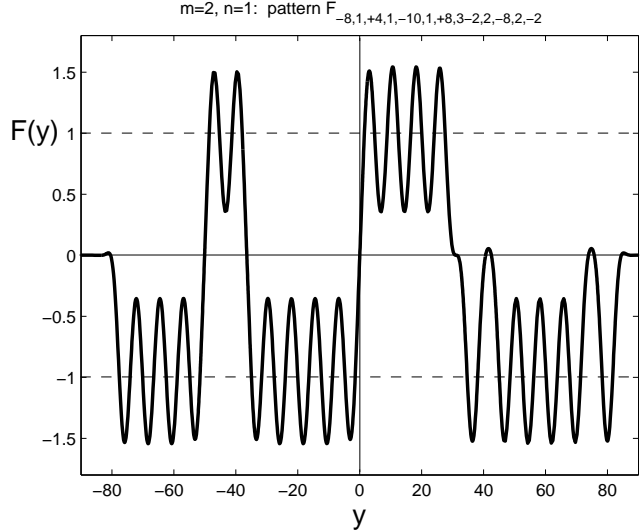}
 \vskip -.4cm
\caption{\rm\small A complicated pattern $F_{\s}(y)$ for (\ref{1})
for $n=1$.}
   \vskip -.3cm
 \label{FC2}
\end{figure}

On the basis of our previous experience of dealing with various
patterns, we will classify other solutions
(possibly, a class of patterns)
by introducing  multiindices  of any length
 \be
 \label{mm1}
 \s=\{\pm \s_1, \s_2, \pm \s_3, \s_4,..., \pm \s_l\}.
  \ee
 Figure \ref{G8}(b)  shows several profiles from the
family with the index
 $
  \s=\{+k,l,-m,l,+k\}.
  $
  In Figure \ref{FC1}, we show further two different patterns, while in
  Figure \ref{FC2}, a single most complicated pattern is
  presented, for which
   \be
   \label{ch1}
   \s=\{-8,1,+4,1,-10,1,+8,1,3,-2,2,-8,2,2,-2\}.
    \ee
   All computations are performed for $n=1$ as usual, and numerics
   were well converged with sufficient accuracy and regularization
 at least  $\sim 10^{-4}-10^{-5}$ and better.
 This shows that the multiindex (\ref{mm1}) can be arbitrary, i.e., can  take
any finite part of any non-periodic fraction. Though we do not
insist that, for a given $\s$, the profile $F_\s(y)$ is unique, we
have seen that the uniqueness fails. Note that the homotopy
approach  \cite{KKVV00, VV02} does not apply to ODEs such as
(\ref{1}) with infinite oscillatory properties.

 Actually, this means
{\em chaotic features}
 of the whole family of solutions $\{F_\s\}$.
 These chaotic types of behaviour are known for other fourth-order ODEs with
coercive operators,  \cite[p.~198]{PelTroy}.

\section{Single point blow-up for $p>n+1$}
 \label{SectLS}

We now consider the ODE (\ref{f11E}) in the case $p>n+1$, which,
in view of the spatial rescaled variable $y$ in (\ref{RVarsE}),
corresponds to single point blow-up. It is crucial that
(\ref{f11E}) for $p \not = n+1$ {\em is not variational}.
Therefore, the solutions of (\ref{f11E}) can be traced out by a
complicated shooting and matching procedures, which are still not
completely justified. For practical reasons, we will use a
continuation in parameters  approach, which allows us to predict
solutions by using those in the variational case $p=n+1$. Recall
that for such ODEs, using a standard inverse function theorem is
not straightforward at all since the differential operator in
(\ref{f11E}) is degenerate and singular. Nevertheless,
 we will suggest to use Schauder's fixed point theorem
 and arrive at some
 convincing conclusions concerning the solvability and the
multiplicity of solutions (the so-called $p$-branches of
solutions).

\subsection{Asymptotics at infinity and first numerical results}

 Recall that, according to our local analysis close to interfaces,
such a behaviour for the TWs given by (\ref{le2}), (\ref{le2N})
asymptotically coincides with that  for the blow-up similarity
solutions; cf. ODEs (\ref{1}) and (\ref{b3}). We can continue to
use this convenient analogy in the present case, since the ODE
(\ref{f11E}) provides us with  similar asymptotics at the
interfaces.
 Obviously, the behaviour at finite interface
$y_0>0$ now corresponds to $\l=+1$,
 and therefore the {\em one-dimensional} bundle (\ref{as11})
is not sufficient to shoot two boundary conditions in (\ref{BCs})
or (\ref{as9}). This explains why we need another asymptotic
expansion as $y \to + \infty$. Thus, compactly supported
similarity profiles are very unlikely (though may exist for some
very special parameters values involved).

Thus, unlike the previous case of regional blow-up for $p=n+1$, in
the present case,
 in order to get a guaranteed successful shooting,
 we must use another full {\em two-dimensional} bundle of  non-compactly supported
  solutions of  (\ref{f11E})
 with the following behaviour:
 \be
 \label{as1}
 f(y) \sim (C_0 y^\g +...) + (C_1 \, {\mathrm e}^{-b_0
 y^\nu}+...) \quad \mbox{as\,\, $y \to +\infty$}.
 \ee
 Here $C_0 \not = 0$ and $C_1 \in \re$ are arbitrary constants
 and
 $$
 \mbox{$
  \g = - \frac 4{p-(n+1)}<0, \quad
 \nu= \frac{4(p-1)}{3[p-(n+1)]}>0, \quad
   b_0=\frac 1 \nu \big( \frac \b{n+1} \, C_0^{-n}\big)^{\frac
   13}>0.
   $}
   $$
Roughly speaking, the first bracket in (\ref{as1}) represents an
``analytic" part of the expansion (e.g., for integer $p$,  in
terms of multiples of  $y^\g$ it can be represented as an analytic
series; the proof of convergence is difficult), while the second
braces gives the essentially ``non-analytic" part.
 Such a structure in (\ref{as1}) is typical for saddle-node
 equilibria, \cite[p.~311]{Perko}. Of course, for the fourth-order
 ODE (\ref{f11E}), this expansion does not admit a simple
 phase-plane interpretation (though the algebraic  origin of the expansion is
 clear). Justification of (\ref{as1}) needs  technical
 applications of fixed point theorems in weighted spaces of
 continuous functions defined for $y \gg 1$; see a typical example
  in \cite[p.~29]{SGKM} and related references.

 One can see passing to the limit $t \to T^-$  in (\ref{RVarsE})
that the asymptotic behaviour (\ref{as1}) gives the following {\em
final-time profile} of this single point blow-up for even profiles
$f=f(|y|)$:
 \be
 \label{ft1}
 u_S(x,T^-)= C_0 |x|^{-\frac 4{p-(n+1)}} < \infty \quad \mbox{for
 all \, $x \not = 0$.}
  \ee

Returning to the asymptotic expansion, we conclude that
(\ref{as1}) represents
 \be
 \label{as2}
 \mbox{a 2D asymptotic bundle.}
  \ee
Hence, the bundle (\ref{as1}) is well suitable for matching with
also two symmetry conditions at the origin (\ref{BCs}), so we
expect  not more than a countable set of solutions. For first
patterns, we keep the same notation as in Section \ref{SectS}) for
$p=n+1$. Moreover, one can expect that these  profiles can be
continuously deformed to each other as $n \to 0^+$.

\begin{figure}
 \centering
\includegraphics[scale=0.65]{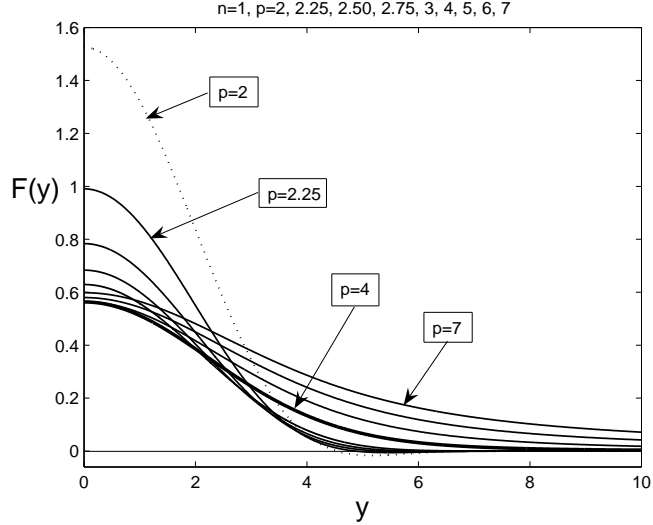}
 \vskip -.4cm
\caption{\rm\small Single point blow-up patterns $F_0(y)$ of
(\ref{f11E}) for $n=1$ and various $p \ge 2$.}
   \vskip -.3cm
 \label{FLS1}
\end{figure}

In Figure \ref{FLS1}, we present the first pattern $F_0(y)$ for
$n=1$ with  $p =2$ (the dotted line for comparison), 2.25, 2.5,
2.75, 3, 4, 5, 6, 7. This  shows
 that, for larger $p$, the
profiles get the positive asymptotic behaviour (\ref{as1}) with
$C_0>0$, and become strictly positive in $\re$.

\begin{figure}
 \centering
\includegraphics[scale=0.65]{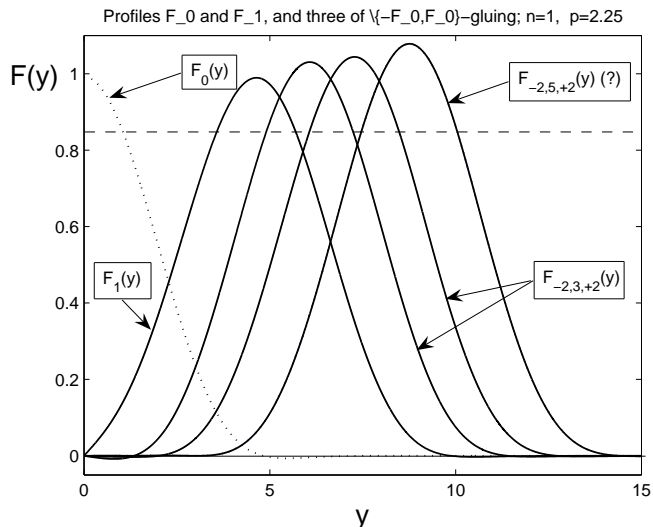}    
 \vskip -.4cm
\caption{\rm\small Single point blow-up patterns $F_0(y)$ and
$F_1(y)$ of (\ref{f11E}) and three profiles from the family
$\{F_{-2,2k+1,+2}\}$ of the $\{-F_0,F_0\}$-gluing; $n=1$, $p=
2.25$.}
   \vskip -.3cm
 \label{FLS2N}
\end{figure}

 Figure \ref{FLS2N}
 shows
   the first  $F_0(y)$ and the second $F_1(y)$ (dipole-like) patterns for
$n=1$ and  $p =2.25$. In addition, we show therein three profiles
from the corresponding (possibly, finite for $p>n+1$) family of
$\{-F_0,F_0\}$-gluing. This family can be viewed as the continuous
$p$-extension of the family shown in Figure \ref{G7} in the
variational case $p=n+1$. Observe that we have two profiles
$F_{-2,3,+2}(y)$, and in the last one $F_{-2,5,+2}(y)$ the
structure of zeros close to $y=0$ was not identified completely
clear. Numerically, this is very difficult and convergence is very
slow with the maximal number 20 000 of mesh points used in {\tt
bvp4c} solver with Tols $\sim 3\cdot 10^{-3}$ achieved.

 These numerics also confirm that the
profiles have non-oscillatory asymptotic behaviour as in
(\ref{as1}), though these can change sign a few (couple of) times,
thus inheriting this finite oscillation property from the
oscillatory one of the type  (\ref{le3}) for $p=n+1=2$
(sufficiently close to the current $p=2.25$).

\subsection{On branching of solutions from variational critical points: degree theory
and Schauder's theorem}

 We again reduce (\ref{f11E}) to the semilinear problem,
 \be
 \label{2}
  \mbox{$
 F=|f|^n f \,\, \Longrightarrow \,\,
  -F^{(4)} - \b(1-\a)|F|^{-\a}F'y- \frac 1{p-1} \, |F|^{-\a}F+
  |F|^{p(1-\a)-1}F=0,
   $}
   \ee
   where  $\a= \frac n{n+1}$.
  For $p=n+1$, i.e., for $\b=0$, this ODE has been studied in Section
\ref{SectS}. Setting $\e=p-(n+1)$
 and
 writing (\ref{2}) as
  \be
  \label{2.11}
   \mbox{$
 {\bf G}(F) \equiv -F^{(4)} - \frac 1{p-1} \, |F|^{-\a}F+
  |F|^{p(1-\a)-1}F= \e \, \frac{1-\a}{4(p-1)}\, |F|^{-\a}F'y,
  $}
   \ee
 we have that the non-autonomous term on the right-hand side becomes
 an asymptotically small perturbation for $|\e| \ll 1$ (in the
 ``weak" integral sense for the equivalent integral equation, so the singularity
 $|F|^{-\a}F'$ is supposed to be eliminated via integration by parts; see (\ref{ss1})).
  On the
 other hand, the operator on the left is variational and hence the
 non-perturbed problem for $\e=0$ admits families of solutions
 described in Section \ref{SectS} (see more details in
 \cite{GPHom}). Classic perturbation and branching
 theory \cite{VainbergTr, KrasZ} suggests that, under natural hypotheses, the variational
 problem for $\e=0$ can generate a countable family of $p$-branches,
 which can be extended for some sufficiently small $|\e|>0$.
The analysis of bifurcation, branching, and continuous extensions
is then performed for the equivalent to (\ref{2.11}) integral
equations with Hammerstein compact operators; see typical examples
in \cite{BGW1, GHUni, GW2}, where similar perturbation problems
for blow-up and global patterns were investigated. Note that these
$p$-branches of solutions are fully extensible and can end up
either at a singularity point or at another bifurcation values.

Thus, for convenience, we write (\ref{2.11}) as follows:
 \be
  \label{ss1}
   \begin{matrix}
 ({\bf B}^*-I) F= h(F,\e) \equiv -\frac 14 \, y F'
 + \e \, \frac{1-\a}{4(p-1)}\, |F|^{-\a}F'y + \frac 1{p-1} \,
 |F|^{-\a}F \qquad
  \smallskip\smallskip\\
-
  |F|^{p(1-\a)-1}F \quad \Longrightarrow \quad
  F= {\bf A}(F,\e) \equiv({\bf B}^*-I)^{-1}h(F,\e) ,\qquad
 \end{matrix}
  \ee
 where $({\bf B}^*-I)^{-1}$ is a compact linear operator in
 $L^2_{\rho^*}(\re)$, \cite{Eg4}.
 For $\e=0$, (\ref{ss1}) gives an integral equation with a variational
 operator, which admits the same critical values and points as the differential one studied
 before.

An efficient way to prove  branching of solutions from $p=n+1$,
which applies for the  lack of differentiability and  regularity
of nonlinearities involved, is using degree-index theory.
This establishes branching from an isolated
solution, say, the first one $F_0$ for simplicity, from the
branching point $\e=0$ provided its index (rotation of the vector
field $I-{\bf A}'(F_0,0)$) satisfies \cite[p.~353]{KrasZ}
 \be
 \label{ind1}
\g= {\rm ind}(F_0,I-{\bf A}'(F_0,0)) \not = 0.
  \ee
Unfortunately, for the equation (\ref{ss1}), this approach hardly
applies, since  the operator contains nonlinearities $|F|^{-\a}F$,
which are not differentiable at 0, so it is not completely clear
how to treat the spectral properties of
the self-adjoint operator ${\bf G}'(F_0)$ for setting in the whole
space $\re$ (not on a bounded interval corresponding to Dirichlet
problems).
  Besides, the  Fr\'echet differentiability of ${\bf A}$
 in $F$ for $\e \not = 0$ is non-existent.
  Note that standard
alternatives of bifurcation-branching theory without
differentiability hypotheses assume sufficient regularity of the
perturbations; cf. \cite[Thm.~28.1]{Deim}. Justification of
branching phenomena in the present problem needs further deeper
analysis and more tricky ``functional topology" involved.


 Continuing to describe  spectral features of ${\bf G}'(F_0)$,
 we note that
 $$
  \l=0 \quad (\mbox{with the eigenfunction} \,\,\, \psi_0 \sim
  F_0')
  $$
  is  an eigenvalue that corresponds to the translational invariance
  of the original PDE (\ref{1.5}) with the infinitesimal generator
  $D_y$. Therefore, $\l=1$ is an eigenvalue of ${\bf A}'(F_0,0)$ (this derivative exists),
  and this leads the {\em critical case}, where computing
  of the index is more difficult and is  performed as in
  \cite[\S~24]{KrasZ}.
   It is more important that
 $$
  \l=1 \quad (\mbox{with the eigenfunction} \,\,\, \psi_1 \sim
  F_0)
  $$
 is also an eigenvalue of ${\bf G}'(F_0)$ and this is associated
 with the generator $D_t$ of the group of
  translations in $t$.

Thus, the index condition (\ref{ind1}) needs special additional
treatment, in particular, associated with spectral properties of
the linearized operator ${\bf G}'(F)$, and this is an open
problem. In this connection, we conjecture
 that $\e=0$ is a {\em point of changing index}
  for the basic
  family $\{F_l\}$,
   and that \cite[Thm.~56.2]{KrasZ} applies to generate
a countable set of  continuous $\e$-curves from any basic pattern
$F_l(y)$
constructed in Section \ref{Sect54}.

\ssk

Continuing branching approach,  for non-differentiable
nonlinearities as in (\ref{ss1}), as an alternative (and less
effective relative to branches detected) approach, we have that
Schauder's Theorem can be applied to get solutions of (\ref{2.11})
for small $\e>0$ and to trace out $p$-branches of the suitable
profiles.
Let us be more precise about this extension of $p$-branches from
the variational critical points at $p=n+1$. Thus, we consider
(\ref{ss1}) as an integral equation in $\re$ using the well-known
spectral properties of the operator ${\bf B}^*$ in Section
\ref{SectExi} with compact resolvent; see \cite{GW2, GHUni} for
similar reductions.
  We then need to apply
 Schauder's fixed point
Theorem \cite[p.~90]{Berger} in the framework of a weighted
$L^p_{\rho^*}$-metric, in which the integral Hammerstein-type
operators involved are naturally compact; see \cite[\S~17]{KrasZ}.
One can see that the right-hand side in (\ref{ss1}) is continuous
in this topology at $\e=0$, so that this gives at least one
solution which is close to $F_0$ for $\e \approx 0$.

This somehow
 settles existence of at least one solution of
(\ref{2}) for small $|\e|>0$ in a 
convex neighbourhood of solutions $F_l$ for $p=n+1$ (i.e.,
$\e=0$). In this framework, uniqueness  becomes a very difficult
problem since the integral operators are not contractive in this
metric. But we then obtain at least a single continuous
$p_l$-branch emerging from $p=n+1$. The global behaviour of
$p$-branches is a hard problem to be tackled next.

\subsection{Bifurcation from constant equilibria}

 Here we  study other bifurcation phenomena
 in this problem: the  $p$-bifurcation
  from the constant
 equilibrium
  $F(y) \equiv F_*$ for the ODE (\ref{2}) (or (\ref{ss1}),
 where difficulties concerning non-differentiable
   nonlinearities do not occur.
   This leads to other non-basic profiles $F$.
To this end, we  use the known spectral properties of the adjoint
operator (\ref{ad1}). Set
 \be
 \label{20}
 F(y)= F_* + Y(y), \quad \mbox{where}
 \quad F_*= f_*^{n+1}
  \ee
 (i.e., $\pm F_*$ are equilibria for (\ref{2})),
 and
 write down (\ref{2}) in the
following form:
 \be
 \label{3}
  \mbox{$
 -Y^{(4)} - \b(1-\a)|F_*|^{-\a}Y'y- \frac 1{p-1} \,(1-\a)|F_*|^{-\a}Y+
  G(y,Y,Y')=0,
   $}
\ee
 where $G(\cdot)$ is the nonlinear part of the operator to be treated
 as a perturbation in the space $L^2_{\rho^*}$.
We next set
 \be
 \label{4}
  \mbox{$
 y=a z, \quad \mbox{where}
  \quad a^4 \b (1-\a)|F_*|^{-\a}= \frac 14.
   $}
   \ee
Then (\ref{3}) reduces to
 \be
 \label{30}
  \mbox{$
  {\bf B}^*  Y \equiv
 -Y^{(4)} - \frac 14 \, Y'y =  \frac 1{p-(n+1)}\, Y
 -a^4 G(z,Y,Y').
   $}
 \ee

We finally apply classic bifurcation-branching theory \cite{KrasZ,
VainbergTr} for the equivalent integral operator with compact
operators in $L^2_{\rho^*}(\re)$ constructed similar to
(\ref{ss1}); see also extra details in \cite{GW1, BGW1, GHUni}. It
follows from the linearized operator in (\ref{30}) that
bifurcation in $p$ can occur at the following points:
 \be
 \label{6}
  \mbox{$
 \frac 1{p-(n+1)} \in \s({\bf B}^*)= \big\{- \frac l4\big\}
 \quad \Longrightarrow \quad p=p_l= n+1 - \frac{4}{l} \quad \big(l >  \frac 4 n \big).
  $}
  \ee
Recall that all the eigenvalues $\l_l=- \frac l4$ of ${\bf B}^*$
are simple, and hence, under typical assumptions, correspond to
bifurcation points; see e.g., \cite[p.~381]{Deim}. The rigorous
justification in the framework of functional setting in  weighted
$L^2_{\rho^*}$-spaces and, especially, asymptotic expansions as
$\e=p-p_l \to 0$ become rather tricky and involved; we do not do
this here and we refer to typical examples in \cite[\S~5]{BGW1}.
 This creates a countable set of
$p$-bifurcation branches to be described next.
The adjoint polynomials $\psi_l^*(y)$ do not change sign (see
(\ref{psi**1})), so that the patterns obtained for $p \approx
p_l^+$ are expected to have a minimal number (sometimes, none) of
intersections with the constant equilibrium $F_*$. We do not know
for sure how to classify these patterns and expect that as $p \to
n+1$ they are converted into the non-basic profiles $F_{+2k}$ from
Section \ref{S5.9} with large $k=k(l) > \frac 4n$ (an open
problem).


\ssk


Thus, in general,   there exists a countable family of
$p$-branches, and some of them, being extended to $p=1^-$, $n=0$,
are expected, after necessary scaling, to match with the countable
set of linear profiles given in (\ref{u11}); cf.
\cite[\S~6.1]{GW2}.
 A reliable rigorous identification of these $p$-branches for
 the integral equation of the type (\ref{ss1})
  is a
 difficult analytical, as well as numerical problem.



\begin{figure}
 \centering \subfigure[$p_0$-branch]{
\includegraphics[scale=0.52]{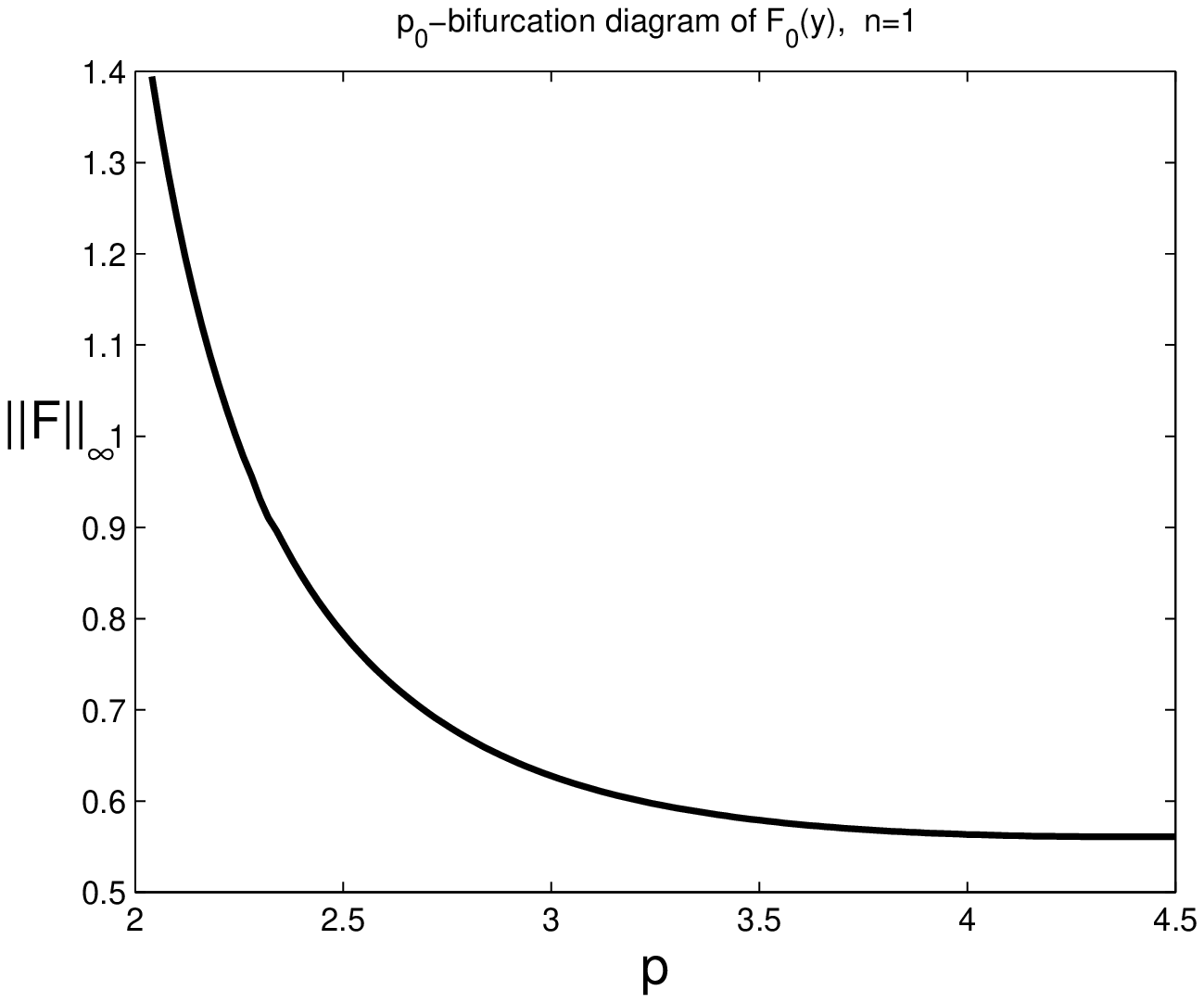}
} \centering \subfigure[$F_0$ profiles]{
\includegraphics[scale=0.52]{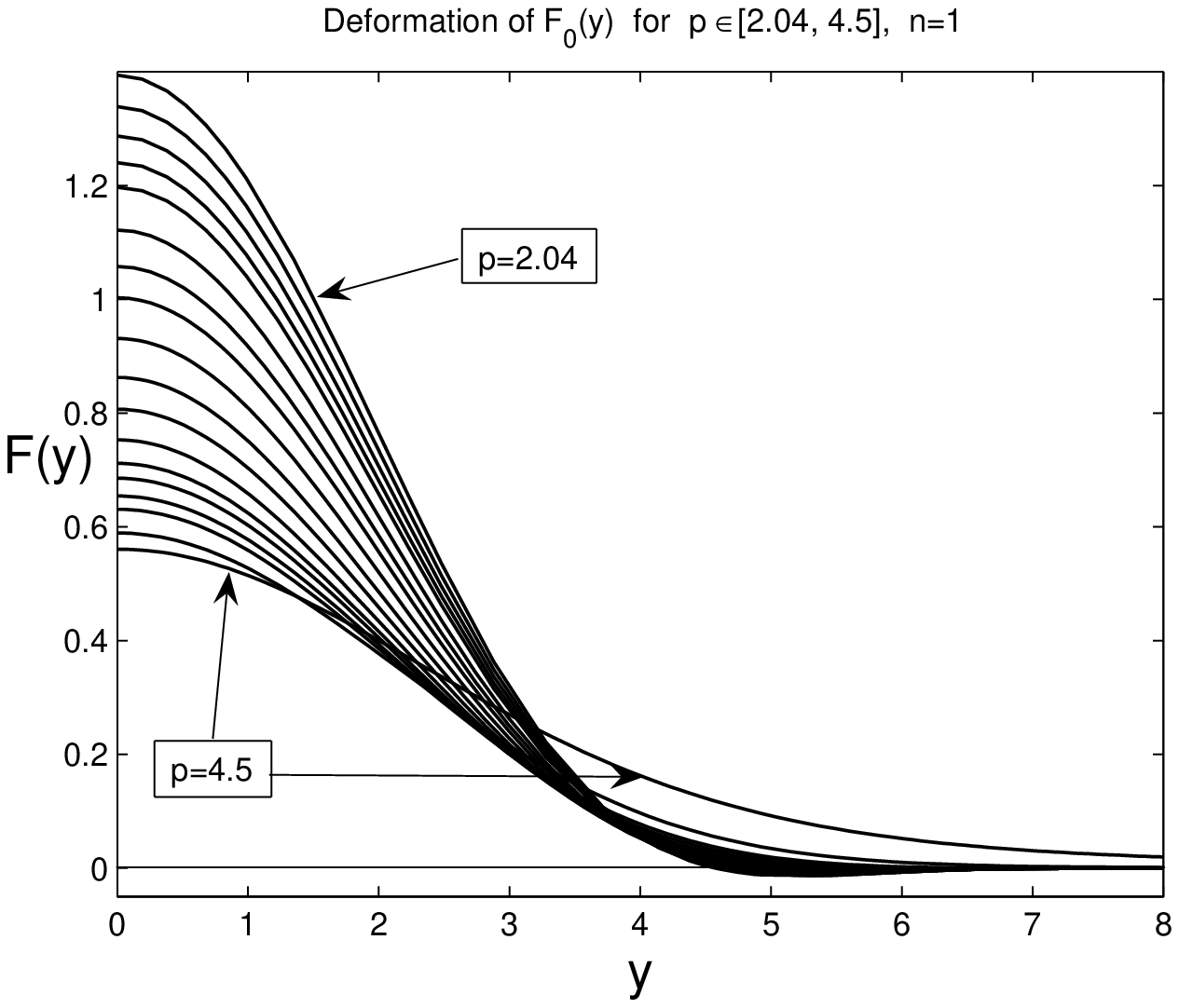}
}
 \vskip -.2cm
\caption{\rm\small  The first $p_0$-branch of solutions $F_0(y)$
of equation (\ref{2}) for $n=1$ (a); corresponding deformation of
$F_0$ (b).}
 \label{Fp0}
\end{figure}

Since (\ref{6}) makes no sense for $l=0$ (and other small $l$'s),
classic bifurcation theory \cite[p.~401]{Deim} suggests that the
basic $p_0$-branch of $F_0(y)$ of the simplest shape (as well as
other first ones with $l \le \frac 4n$) is assumed to exist for
all $p>1$. In other words, this branch cannot appear at
bifurcation points such as those indicated in (\ref{6}). In Figure
\ref{Fp0}, we show the first $p$-branch of $F_0$ (a) and the
deformation of the profiles $F_0(y)$, (b). We expect that this
$p_0$ branch is composed from stable solutions and hence
represents the generic asymptotic blow-up behaviour for the
parabolic PDE (\ref{1.5}).

\begin{figure}
 \centering \subfigure[$p_1$-branch]{
\includegraphics[scale=0.52]{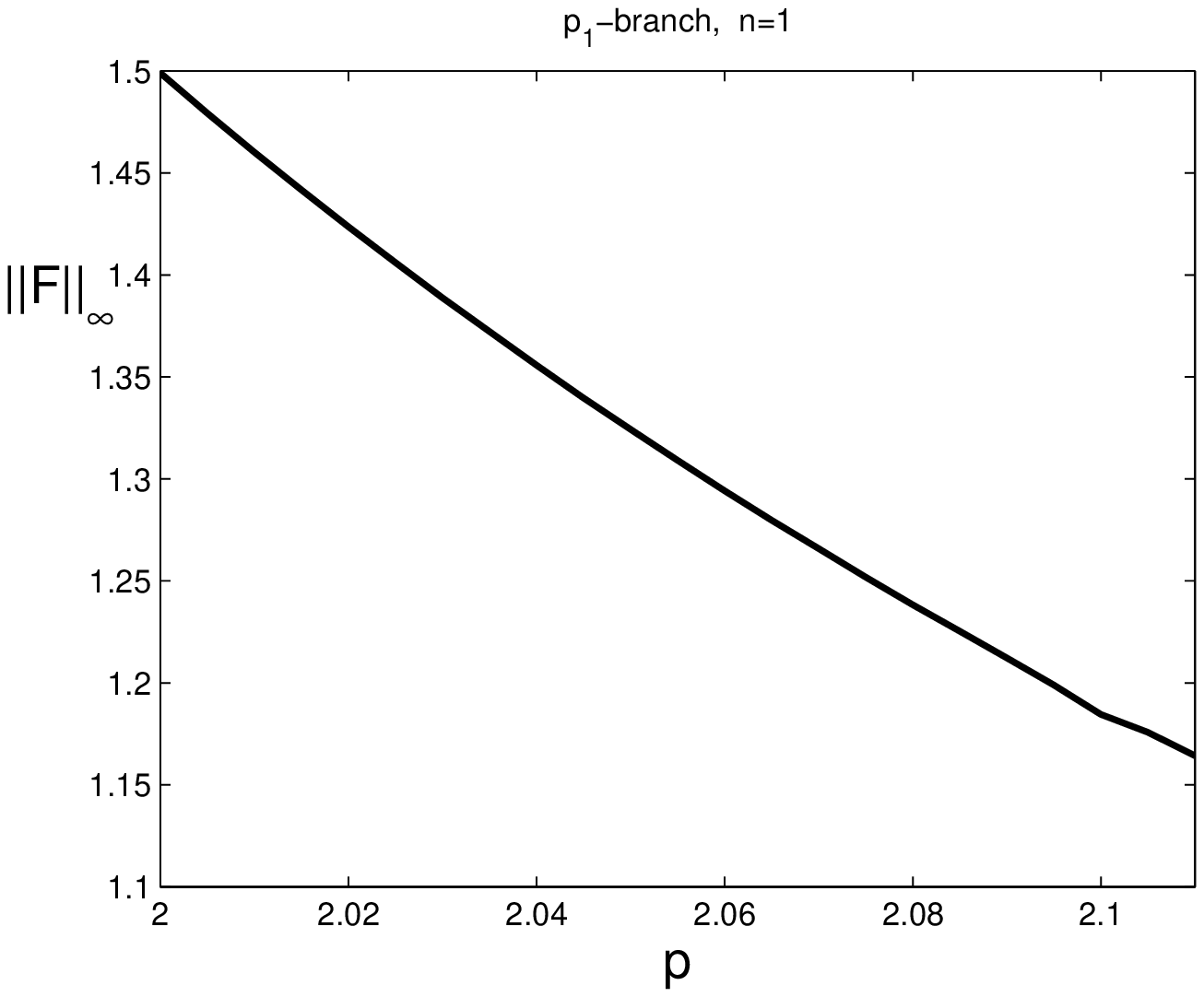}
} \centering \subfigure[$F_1$ profiles]{
\includegraphics[scale=0.52]{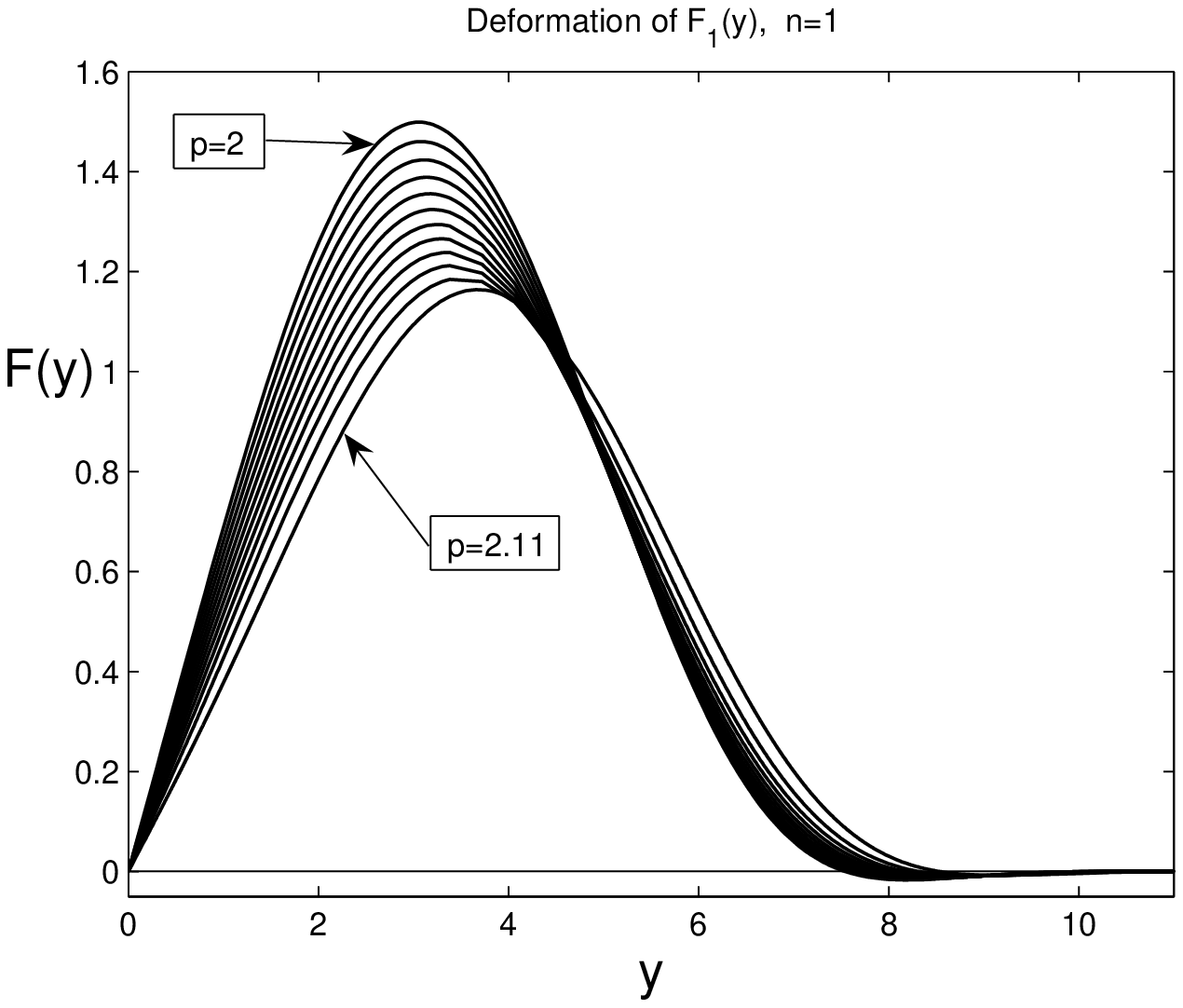}
}
 \vskip -.2cm
\caption{\rm\small  The second $p_1$-branch of solutions $F_1(y)$
of equation (\ref{2}) for $n=1$ (a); corresponding deformation of
$F_1$ (b).}
 \label{Fp1}
\end{figure}

A small part of the next $p_1$-branch of dipole-like profiles
$F_1(y)$ for $n=1$ and $p \in [2,2.11]$ is shown in Figure
\ref{Fp1}(a), where (b) demonstrates the corresponding deformation
with $p$ of $F_1(y)$. Further extension of this branch beyond
$p=2.17$ leads to strong instabilities where the profiles suddenly
jump to different   shapes (which possibly belong to other
$p$-branches nearby having the geometric structure as in Figure
 \ref{G7}; we did not construct such neighbouring   branches).


The $p$-branches can connect various profiles, with rather obscure
understanding of possible geometry of such branches and their
saddle-node bifurcations (turning) points.
 For $p=n+1$, the questions on connections with respect to regularization parameters
  are addressed in \cite[\S~7]{GPHom} posing problems of
  homotopy classification of patterns in variational problems and
  approximate ``Sturm index" of solutions.

\begin{figure}
 \centering
\includegraphics[scale=0.70]{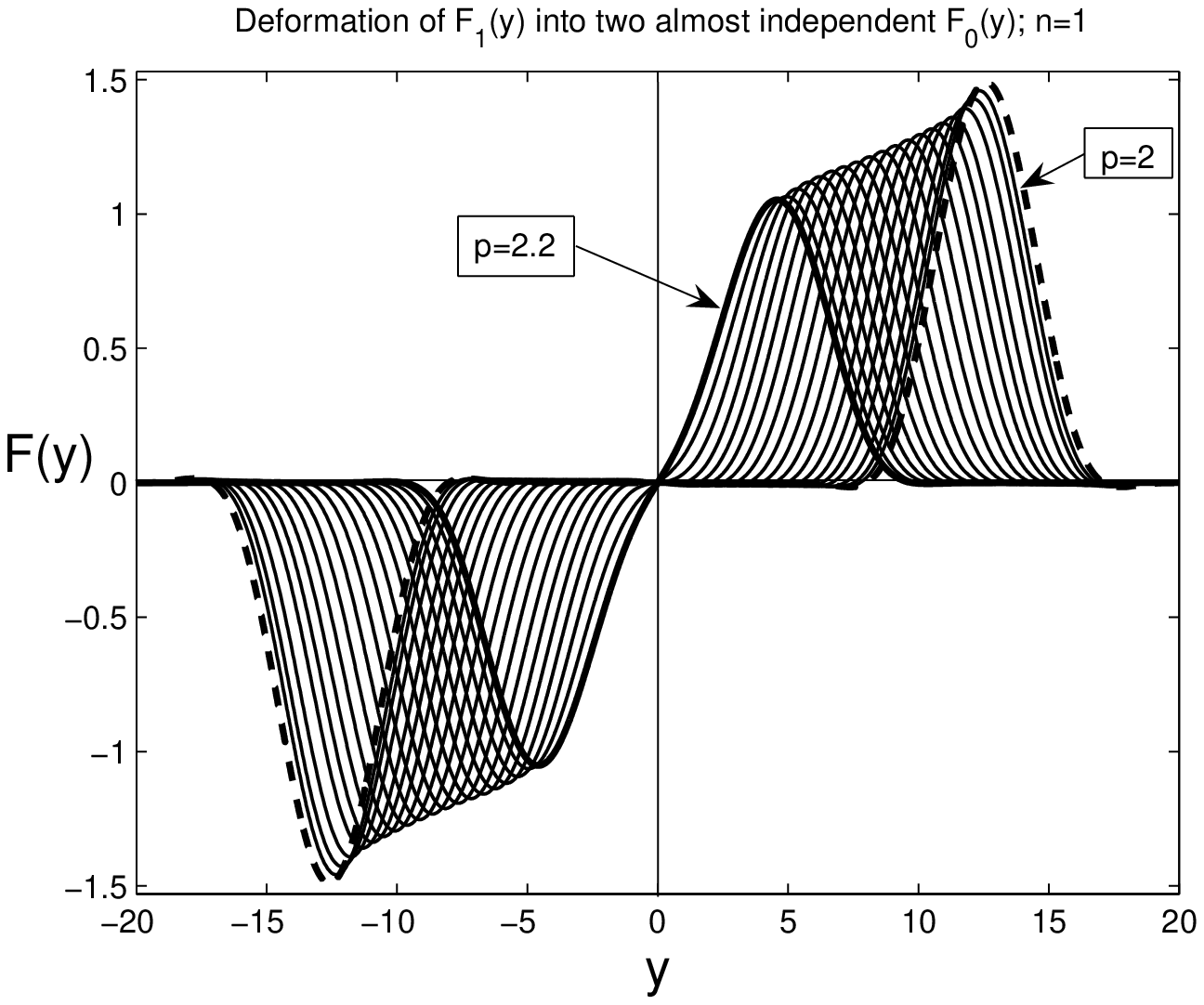}    
 \vskip -.4cm
\caption{\rm\small Deformation of $F_1(y)$ for $p \in [2,\, 2.2]$;
$n=1$.}
   \vskip -.3cm
 \label{FNN1SS}
\end{figure}

 For instance, Figure
\ref{FNN1SS} for $n=1$ shows the connection of the profile
$F_1(y)$ for $p=2.2$ (cf. the previous figure) with two ``almost
independent" mutually shifted  profiles $\pm F_0$ for $p=2$. Next
Figure \ref{FNN1SSN}, in the enlarged form, explains formation of
the zero set of profiles and shows in (b)
 that the eventual structure for $p=2$
actually belongs for a member of the family $\{-F_0,F_0\}$-gluing,
i.e., it is the profile $F_{-2,7,+2}(y)$, with exactly seven
transversal zeros between $\pm F_0(y)$ structures.

\begin{figure}
 \centering \subfigure[enlargement $10^{-4}$]{
\includegraphics[scale=0.52]{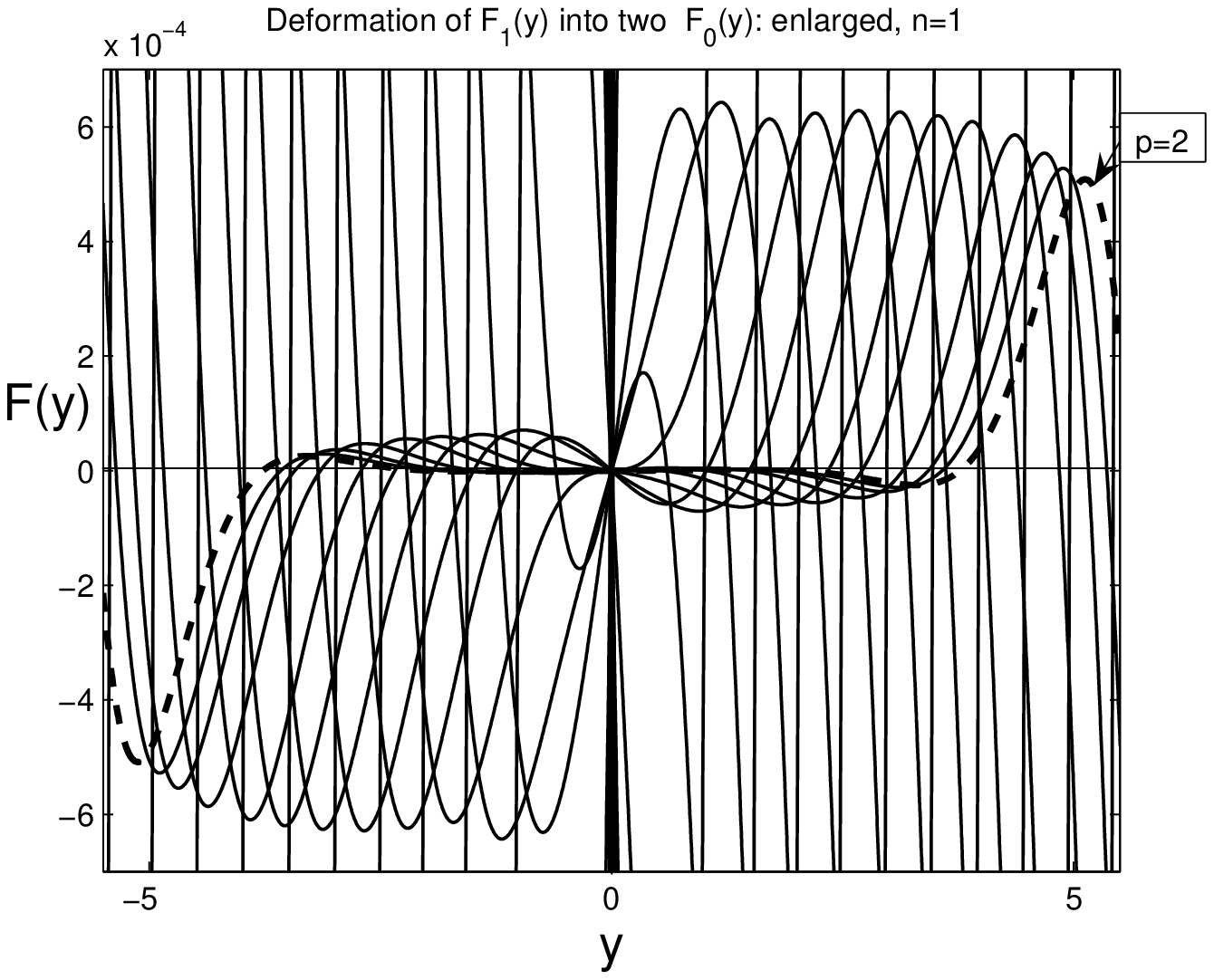}    
} \centering \subfigure[enlargement $10^{-6}$]{
\includegraphics[scale=0.52]{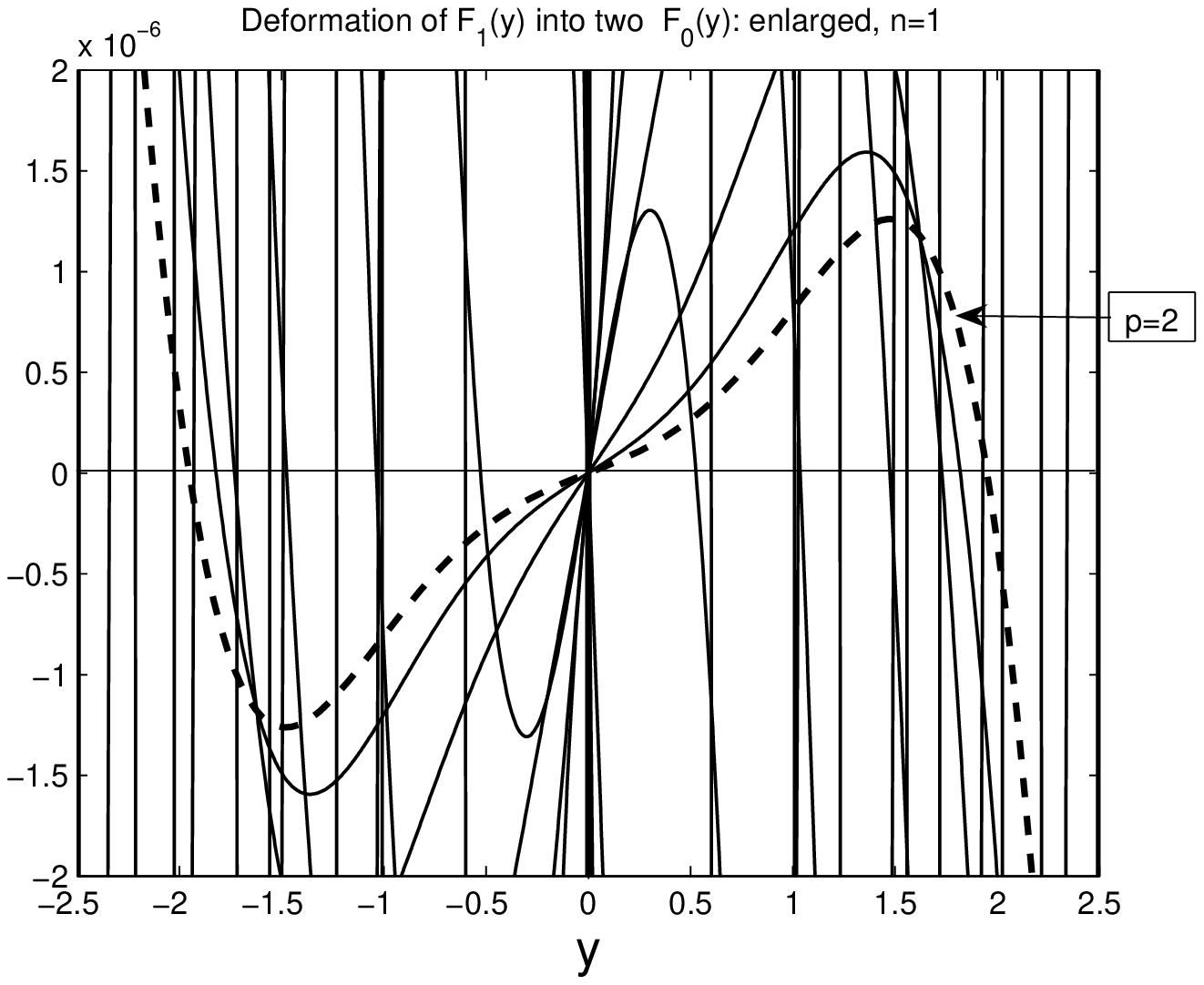}    
}
 \vskip -.2cm
\caption{\rm\small   Enlarged zero set for profiles in Figure
\ref{FNN1SS}; (b) shows that for $p=2$, it is $F_{+2,7,-2}(y)$.}
 \label{FNN1SSN}
\end{figure}

\begin{figure}
 \centering \subfigure[$p_2$-branch]{
\includegraphics[scale=0.52]{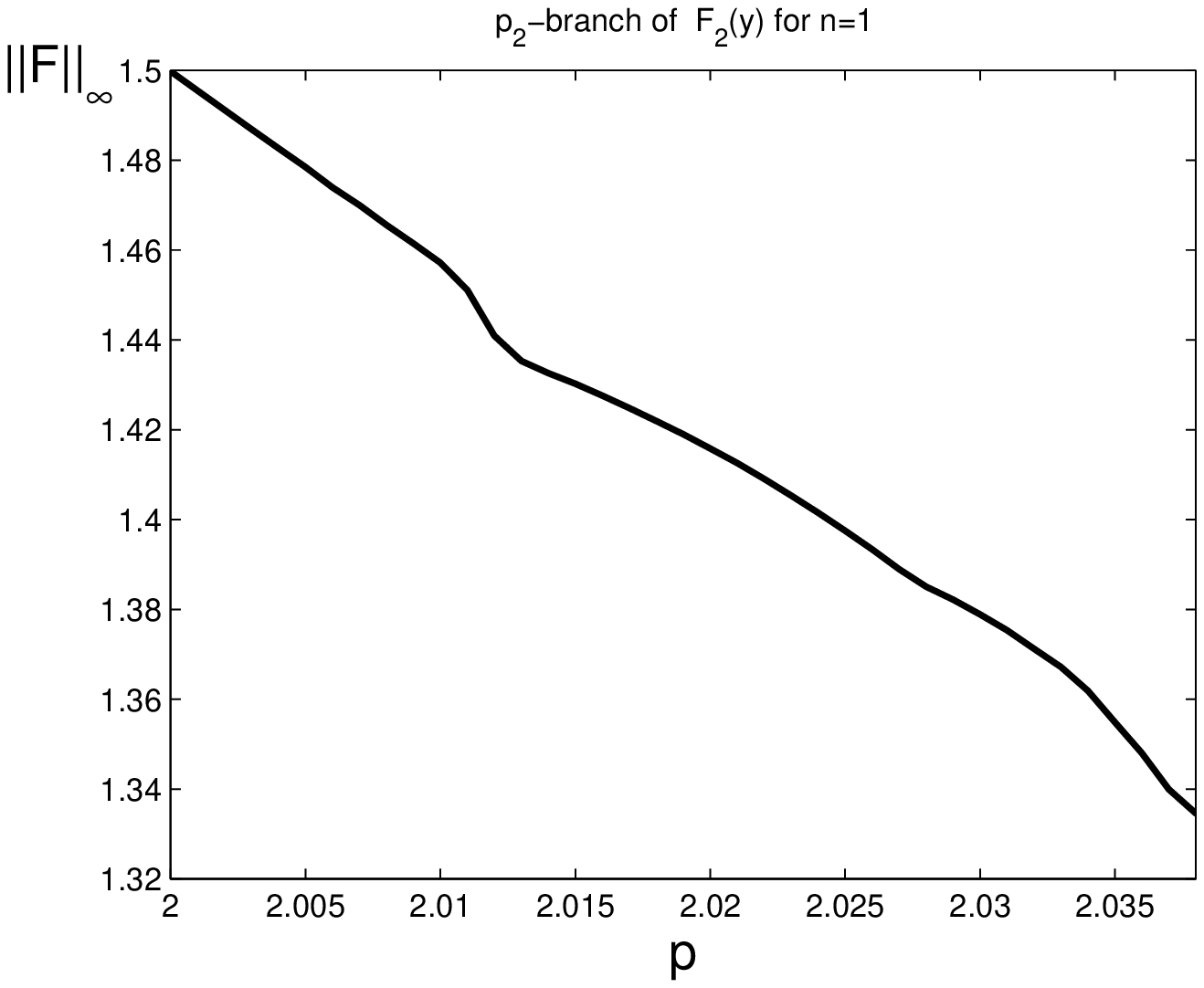}
} \centering \subfigure[$F_2$ profiles]{
\includegraphics[scale=0.52]{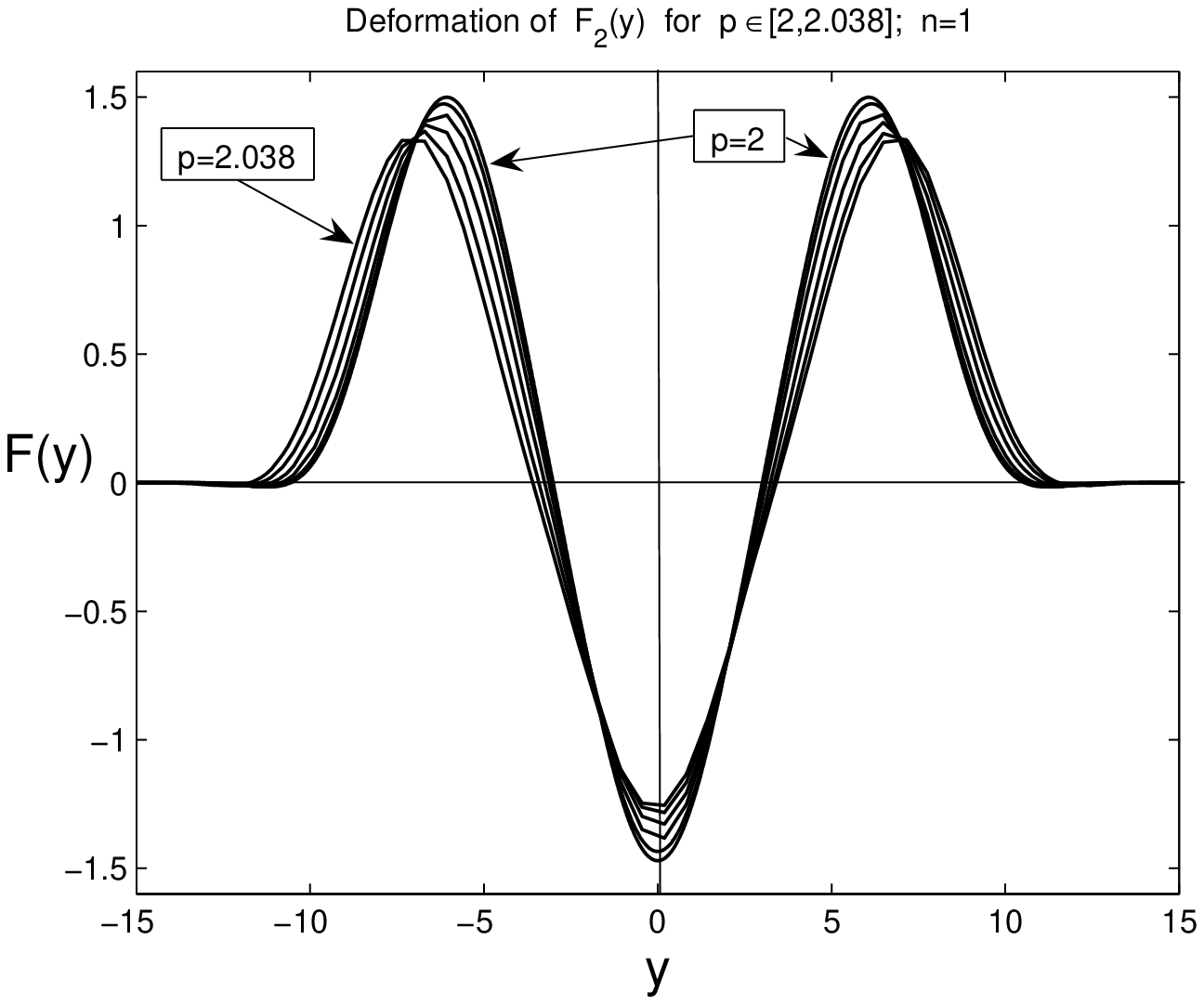}
}
 \vskip -.2cm
\caption{\rm\small  The third $p_2$-branch of solutions $F_2\equiv
F_{+2,1,-2,1,+2}$ of equation (\ref{2}) for $n=1$ (a);
corresponding deformation of $F_2$ (b).}
 \label{Fp2}
\end{figure}

The principle fact that higher-order $p$-branches of the basic
family $\{F_l\}$ (see Section \ref{Sect54}) can be originated from
$p=n+1$ is illustrated in Figure \ref{Fp2} where we show the
$p_2$-branch of profiles $F_2(y)$ (a) (for $p=n+1$, this profiles
is given in Figure \ref{G4}(c)) and the deformation of the
profiles (b).

\begin{figure}
 \centering \subfigure[$p$-branch]{
\includegraphics[scale=0.52]{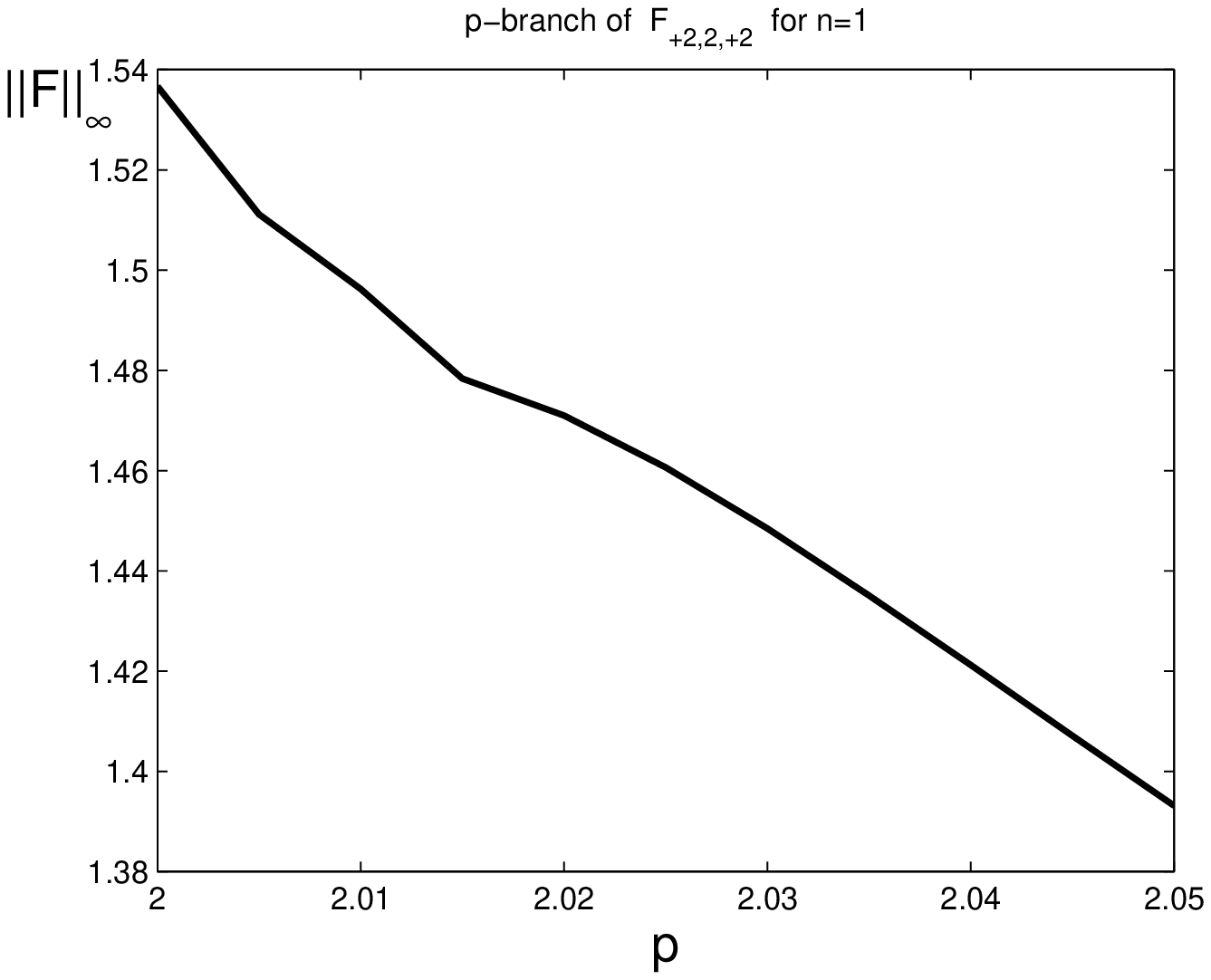}
} \centering \subfigure[$F_{+2,2,+2}$ profiles]{
\includegraphics[scale=0.52]{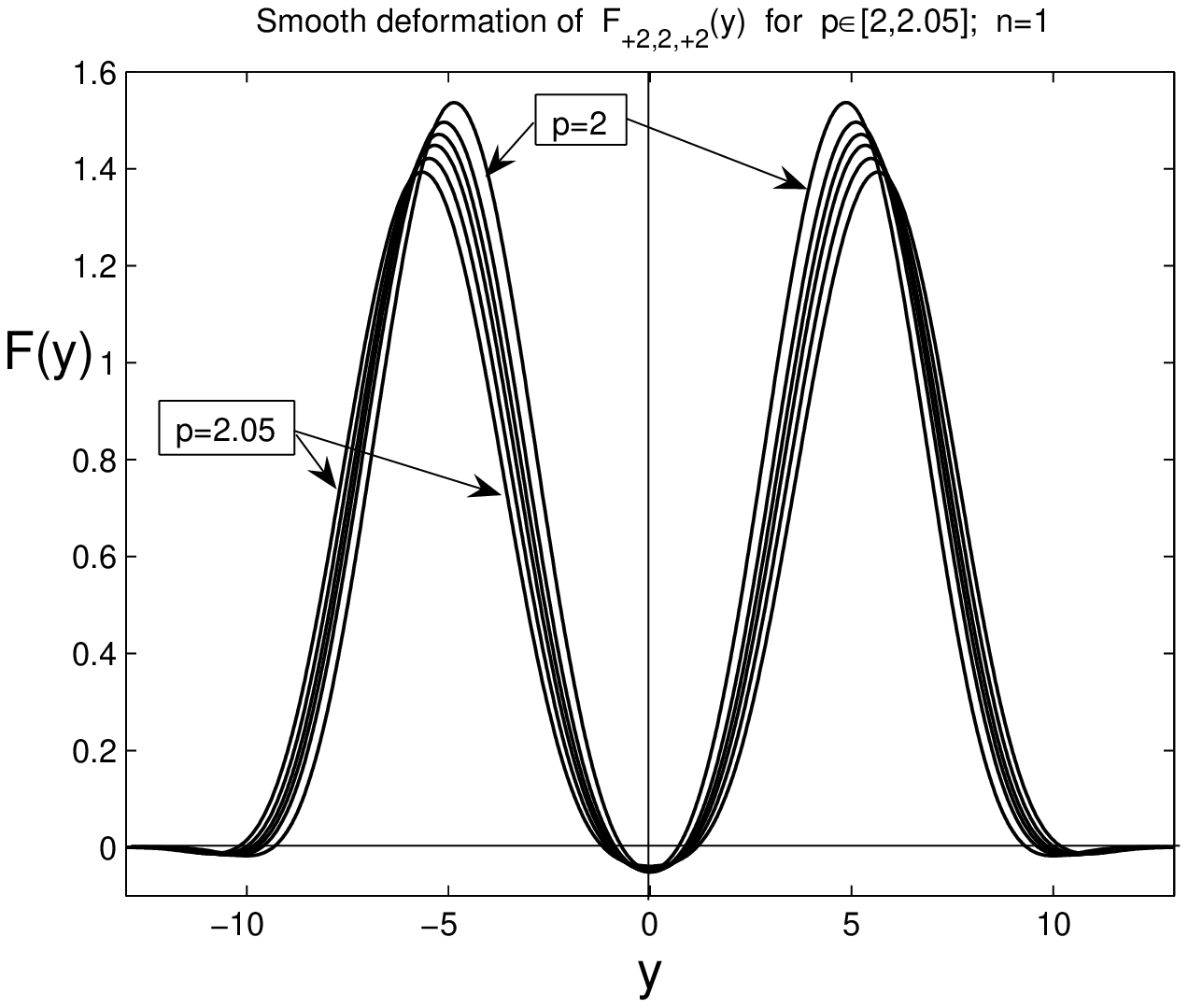}
}
 \vskip -.2cm
\caption{\rm\small  The  $p$-branch of solutions $ F_{+2,2,+2}$ of
equation (\ref{2}) for $n=1$ (a); corresponding deformation of
$F_{+2,2,+2}$ (b).}
 \label{Fp22}
\end{figure}

Another principal issue  is justified in Figure \ref{Fp22} that
shows the $p$-branch of the profile $F_{+2,2,+2}$(y); cf. Figure
\ref{G6} for $p=n+1$.

\begin{figure}
 \centering
\includegraphics[scale=0.65]{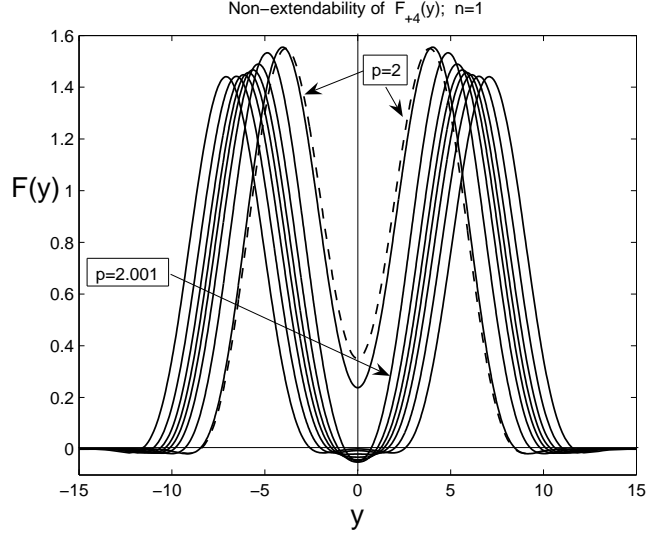}    
 \vskip -.4cm
\caption{\rm\small $F_{+4}$ for $p=n+1$ jumps to $F_{+2,2,+2}$ if
$p$ is increased by $10^{-3}$;  $n=1$.}
   \vskip -.3cm
 \label{FNN1}
\end{figure}

Concerning other, more complicated profiles for $p=n+1$, such as
$F_{+4}(y)$ and others containing such structures shown in Figures
\ref{G6} and \ref{FC1}, numerical results suggest that these
cannot be extended for $p>n+1$. For instance, Figure \ref{FNN1}
demonstrates that the profile $F_{+4}(y)$ very quickly jumps to
the type $F_{+2,2,+2}$ (cf. Figure \ref{Fp22}) for the increment
$\D p=10^{-3}$, and even for smaller $\D p$'s. This
 illustrates
the fact that
 both profiles are very close and are originated at the same branching
 point. We do not intend here to get numerically a correct bifurcation diagram
 for this delicate case. It seems that this demands more advanced
 parameter continuation techniques that are available in the standard {\tt
 MatLab} environment.
  See a similar unstable case below for $p < n+1$.


\ssk

 \noi{\bf Remark: on $\mu$-bifurcations.}
 The origin of various branches of solutions can be also seen via
  an additional  parameterization; cf.
\cite[\S~4.3]{BGW1}. Namely, we consider equation (\ref{3}), where
we introduce the parameter by replacing
  $$
  \b \mapsto \mu >0.
  $$
  Then, linearizing as above, we conclude that bifurcation points
  are (for even profiles)
   \be
   \label{8}
    \mbox{$
   \mu_l= \frac 1l, \quad l =
    2, \, 4, \, ... \, .
    $}
    \ee
Since
 $
 \b < \frac 14 < \frac 12,
  $
we have that at least two bifurcation branches originated at
$\mu_2$ and $\mu_4$, we  have a good chance to be extended to
$\mu=\b$ (this global continuation is an open mathematical
problem, but can be checked numerically), and hence generate
solutions of the original ODE (\ref{f11E}).

\section{Global blow-up similarity profiles for $p \in (1,n+1)$}
 \label{SectHS}

 For $p\in (1,n+1)$, we have $\b<0$ in
(\ref{RVarsE}), so that the similarity solutions describe
expanding as $t \to T^-$ waves with global blow-up uniformly on
any compact subset in $x$.

\subsection{Same oscillatory behaviour close to interfaces}

The ODE (\ref{f11E})  reads for $f \approx 0$ 
 $$
 -(|f|^nf)^{(4)} - \b y f'+...=0 \quad (\b < 0),
 $$
 so, replacing $y_0-y \mapsto y$, on integration for $y \approx 0$,
we have
  $$
(|f|^nf)''' =+ \b y_0 f+... \, .
 $$
 This gives  (\ref{le2}), where $\l = \b y_0 <0$ is
 reduced to $-1$ by scaling.
 Thus, for $p \in(1,n+1)$, the similarity profiles are
  oscillatory near interfaces as for $p=n+1$.
 By (\ref{as55}), the 2D asymptotic bundle is enough to
match two symmetry boundary conditions (\ref{BCs}),
and  the proof of matching remains open.

  \subsection{Profiles and branches}

\begin{figure}
 \centering
\includegraphics[scale=0.65]{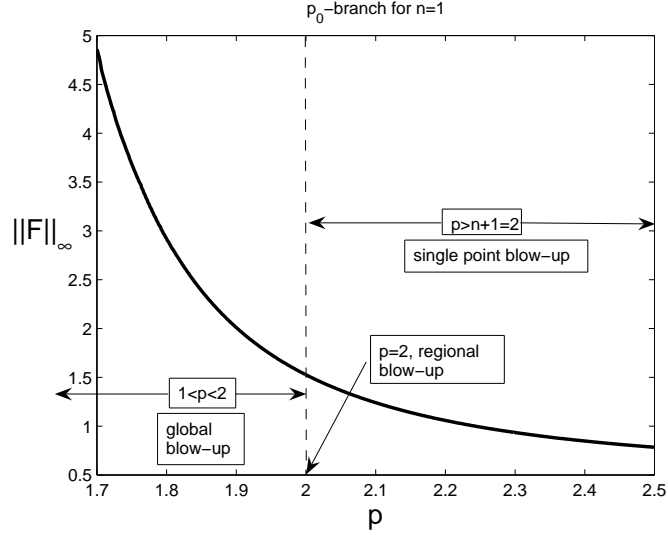}
 \vskip -.4cm
\caption{\rm\small The $p_0$-branch of $F_0(y)$ of the ODE
(\ref{2}) for $n=1$.}
   \vskip -.3cm
 \label{F03N}
\end{figure}

\begin{figure}
\includegraphics[scale=0.7]{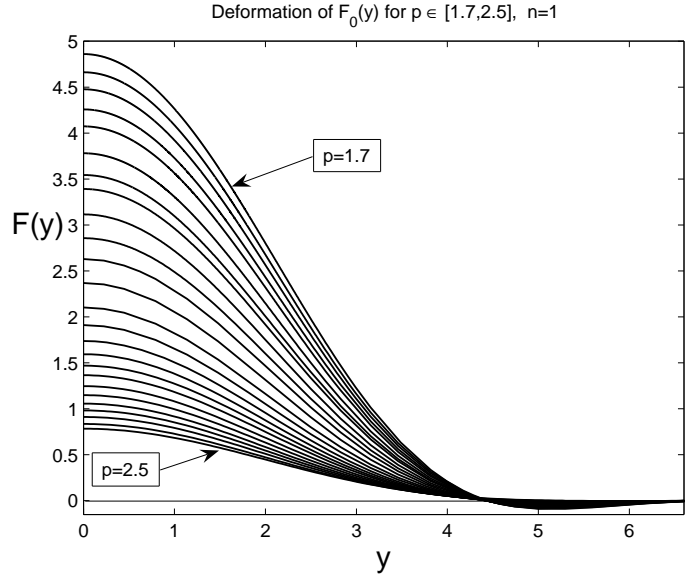}
 \vskip -.4cm
\caption{\rm\small Deformation of $F_0(y)$ from Figure
\ref{F03N}.}
   \vskip -.3cm
 \label{F04N}
\end{figure}

 It turns out that, for $1<p<n+1$, the ODE (\ref{2}) is more difficult to solve
  numerically than for $p \ge n+1$.
   In Figure \ref{F03N}, for $n=1$, we present the $p_0$-branch of
  the generic profiles $F_0$ by continuation in $p$ from $p=2.5$
  (single point blow-up)
  until $p=1.7 < n+1=2$ (global blow-up).
   The actual continuous  $p$-deformation of functions $F_0(y)$ is seen from Figure \ref{F04N}.
  For smaller $p>1$, the
  similarity profiles $F_0(y)$ become larger (and according to
  (\ref{g1N}) the corresponding constant equilibria diverge
  exponentially as $p \to 1^-$) and more oscillatory close to the
  interface.
Note   that the $p_0$-branch is expected to consist of
asymptotically (structurally) stable blow-up profiles $F_0(y)$,
but we cannot prove this even in the linearized approximation (the
linearized operator is a difficult non-self-adjoint operator with
unknown spectrum).

\begin{figure}
 \centering \subfigure[$p_0$-branch]{
\includegraphics[scale=0.52]{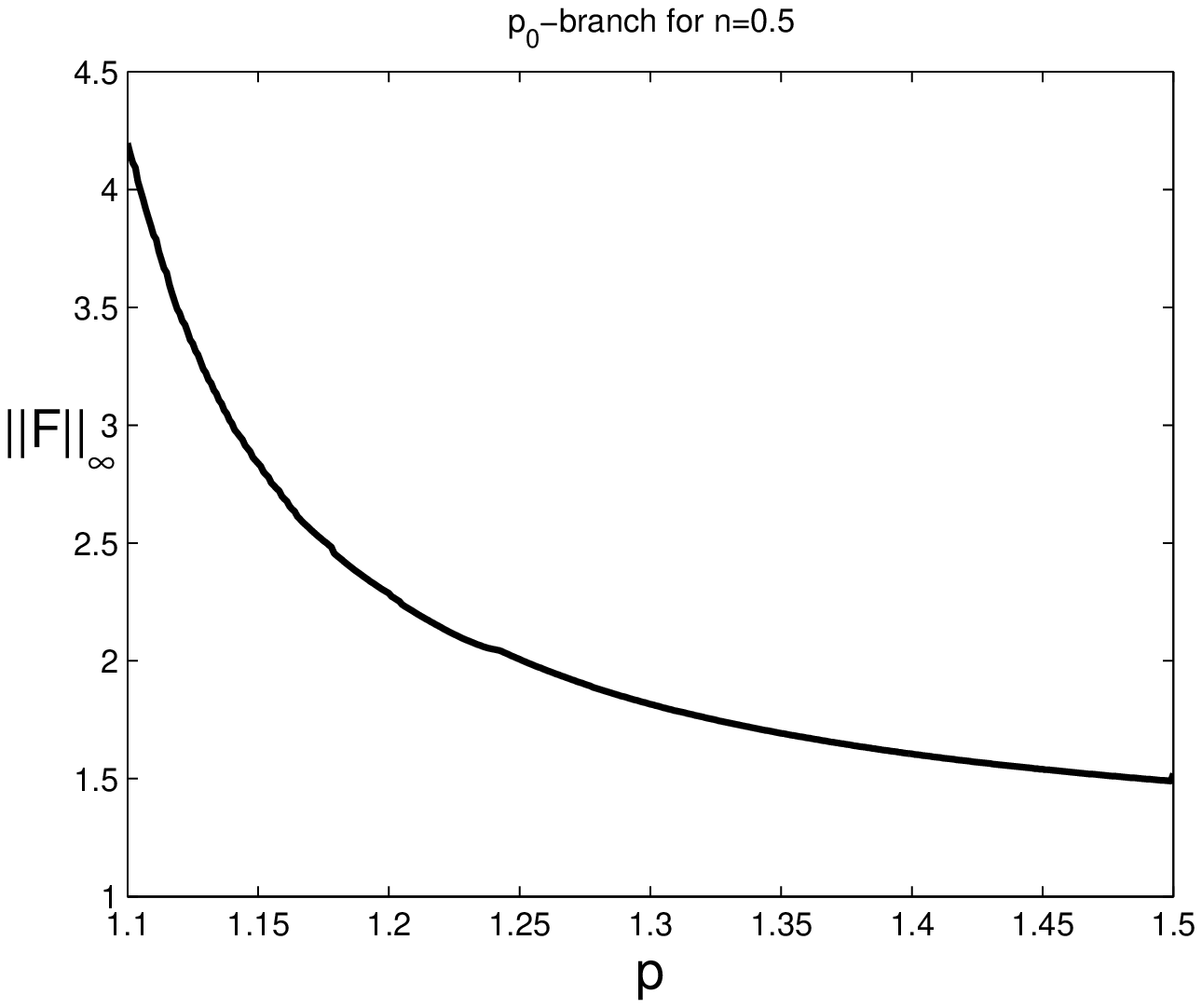}
} \centering \subfigure[$F_0$ profiles]{
\includegraphics[scale=0.52]{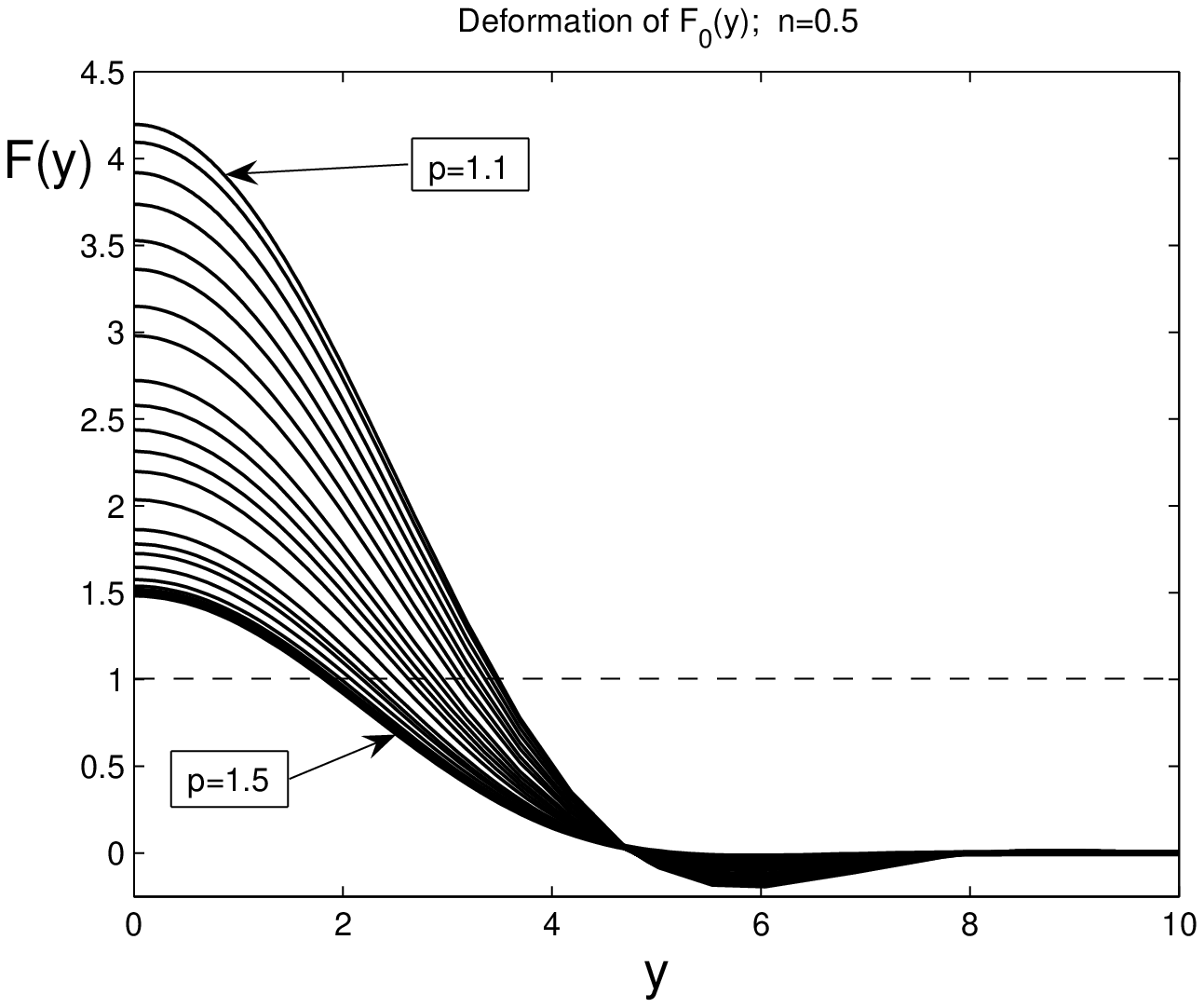}
}
 \vskip -.2cm
\caption{\rm\small  The  $p$-branch of solutions $ F_0(y)$ of
equation (\ref{2N}) for $n=0.5$ (a); corresponding deformation of
$F_{+2,2,+2}$ (b).}
 \label{F001}
\end{figure}

\begin{figure}
 \centering \subfigure[$p_2$-branch]{
\includegraphics[scale=0.52]{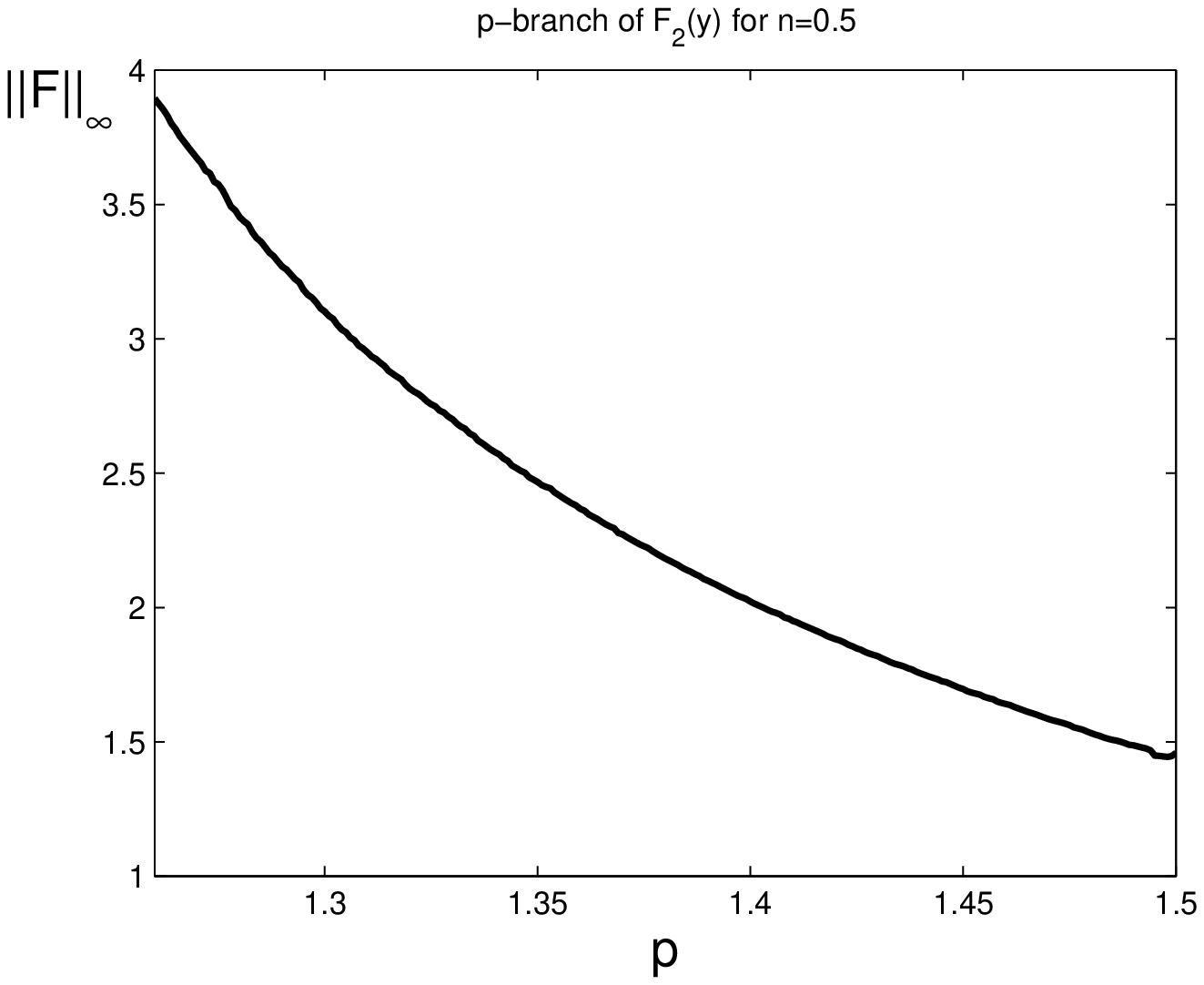}
} \centering \subfigure[$F_2$ profiles]{
\includegraphics[scale=0.52]{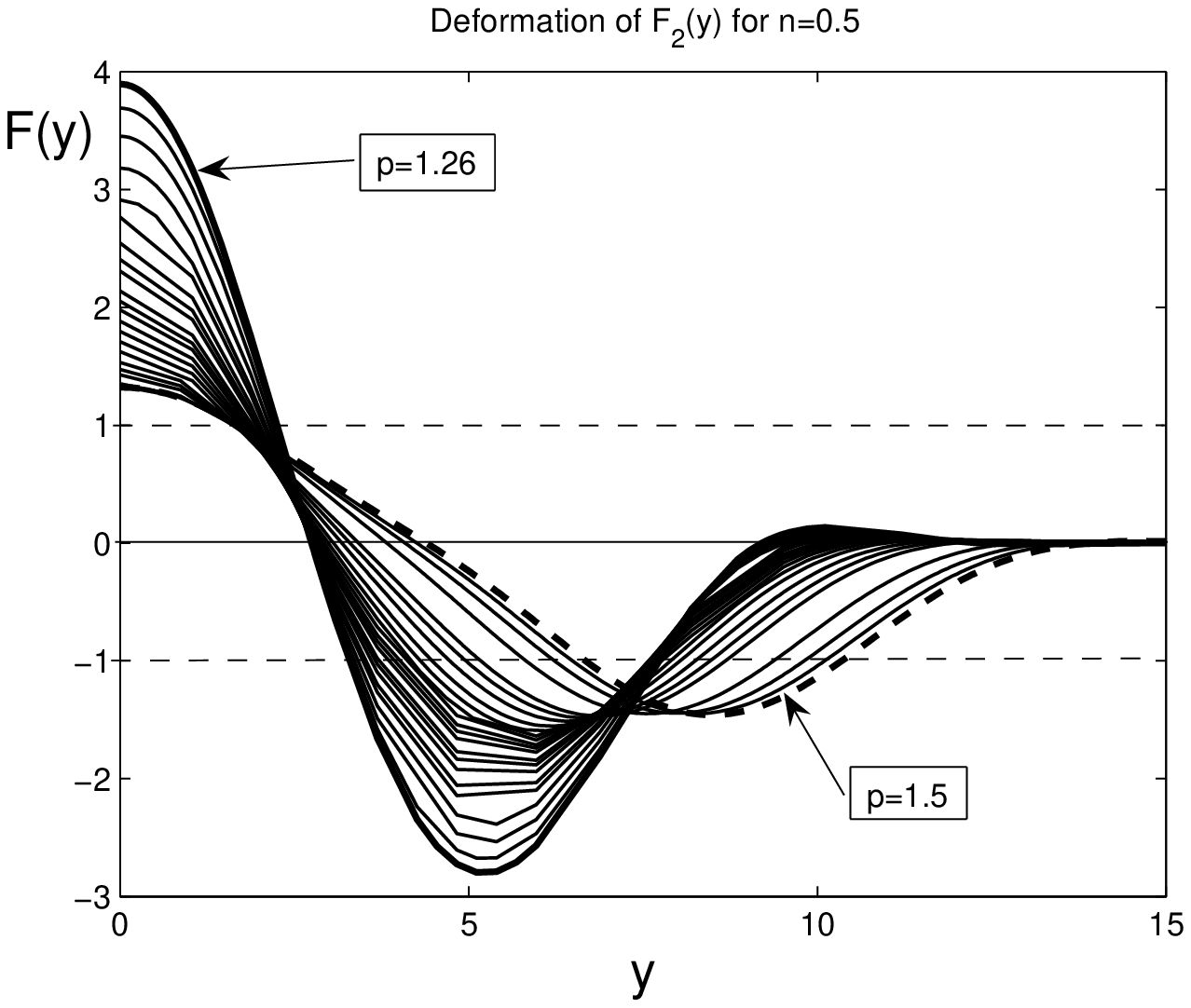}
}
 \vskip -.2cm
\caption{\rm\small  The  $p$-branch of solutions $ F_2(y)$ of
equation (\ref{2N}) for $n=0.5$ (a); corresponding deformation of
$F_{+2,2,+2}$ (b).}
 \label{F201}
\end{figure}

In order to avoid the exponential discrepancy (\ref{g1N}) of
branches as $p \to 1^-$, we now  perform in the ODE (\ref{2}) the
following additional scaling:
 \be
 \label{2NN}
  \mbox{$
  F \mapsto C F, \,\,\, y \mapsto a y, \quad C^{(1-\a)(p-1)}=
  \frac 1{p-1}, \,\,\,a^4=C^{-(1-\a)(p-1)} \quad \big(C=F_*^{\frac
  1{1-\a}}\big)
   $}
    \ee
 to get the equation with fixed equilibria $F_* = \pm 1$ and 0:
  \be
  \label{2N}
   \mbox{$
  -F^{(4)} - \b(p-1)(1-\a)|F|^{-\a}F'y- |F|^{-\a}F+
  |F|^{p(1-\a)-1}F=0 \quad \big( \a= \frac n{n+1}\big).
   $}
   \ee
In Figure \ref{F001}, we show the $p_0$-branch for the ODE
(\ref{2N}) with $n=0.5$ and $p \in [1.1, 1.5]$.
In Figure \ref{F201}, we show the $p_2$-branch of the profiles
$F_2(y) \equiv F_{-2,1,+2,1,-2}(y)$ for the ODE (\ref{2N}) with
$n=0.5$ (a small
 ``discontinuity" is seen there).

\begin{figure}
 \centering \subfigure[$p$-branch]{
\includegraphics[scale=0.52]{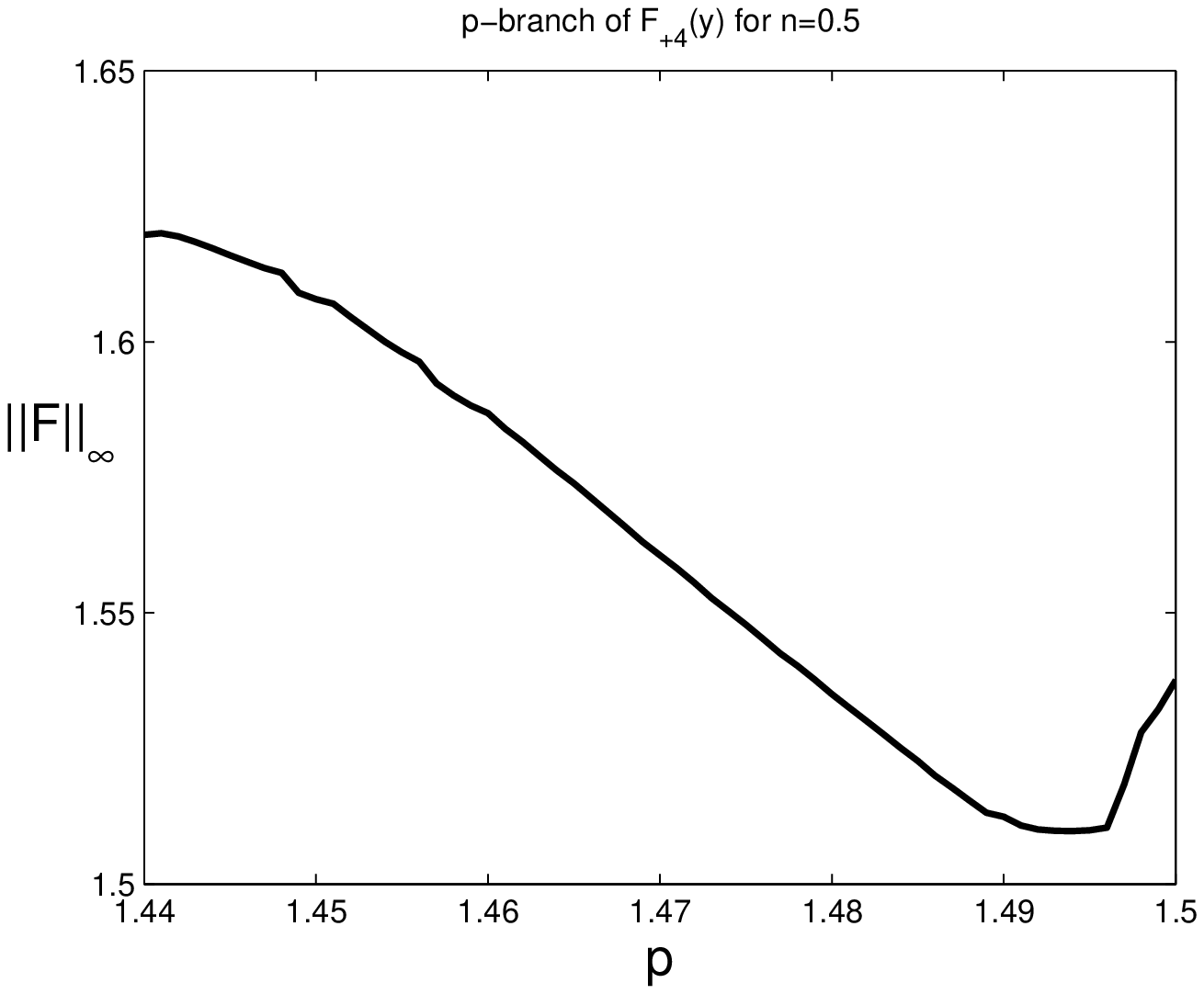}
} \centering \subfigure[$F_{+4}$ profiles]{
\includegraphics[scale=0.52]{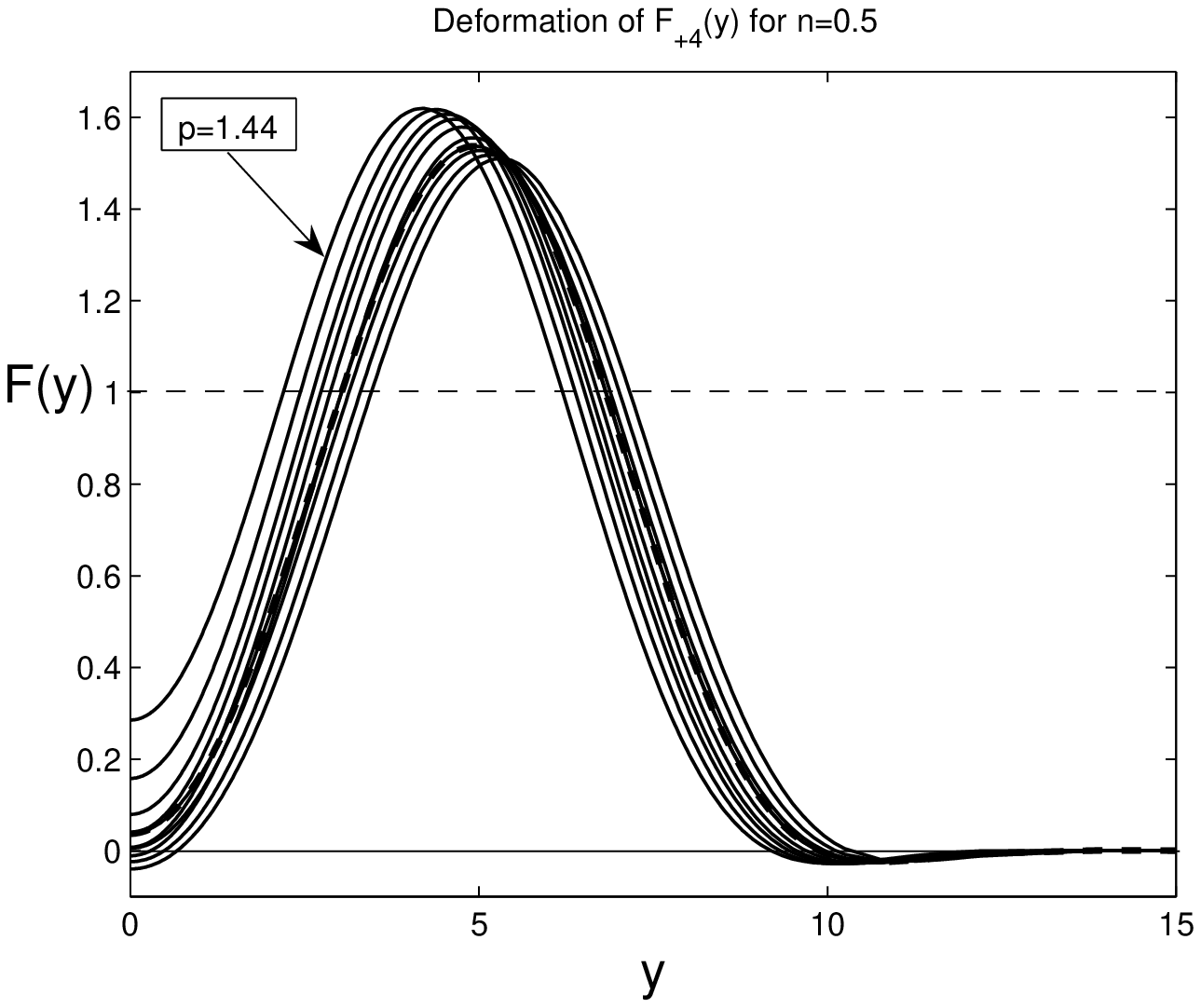}
}
 \vskip -.2cm
\caption{\rm\small  The  $p$-branch of solutions $ F_{+4}(y)$ of
equation (\ref{2N}) for $n=0.5$ (a); corresponding deformation of
$F_{+4}$ (b).}
 \label{FG1}
\end{figure}

\begin{figure}
 \centering \subfigure[$p$-branch]{
\includegraphics[scale=0.52]{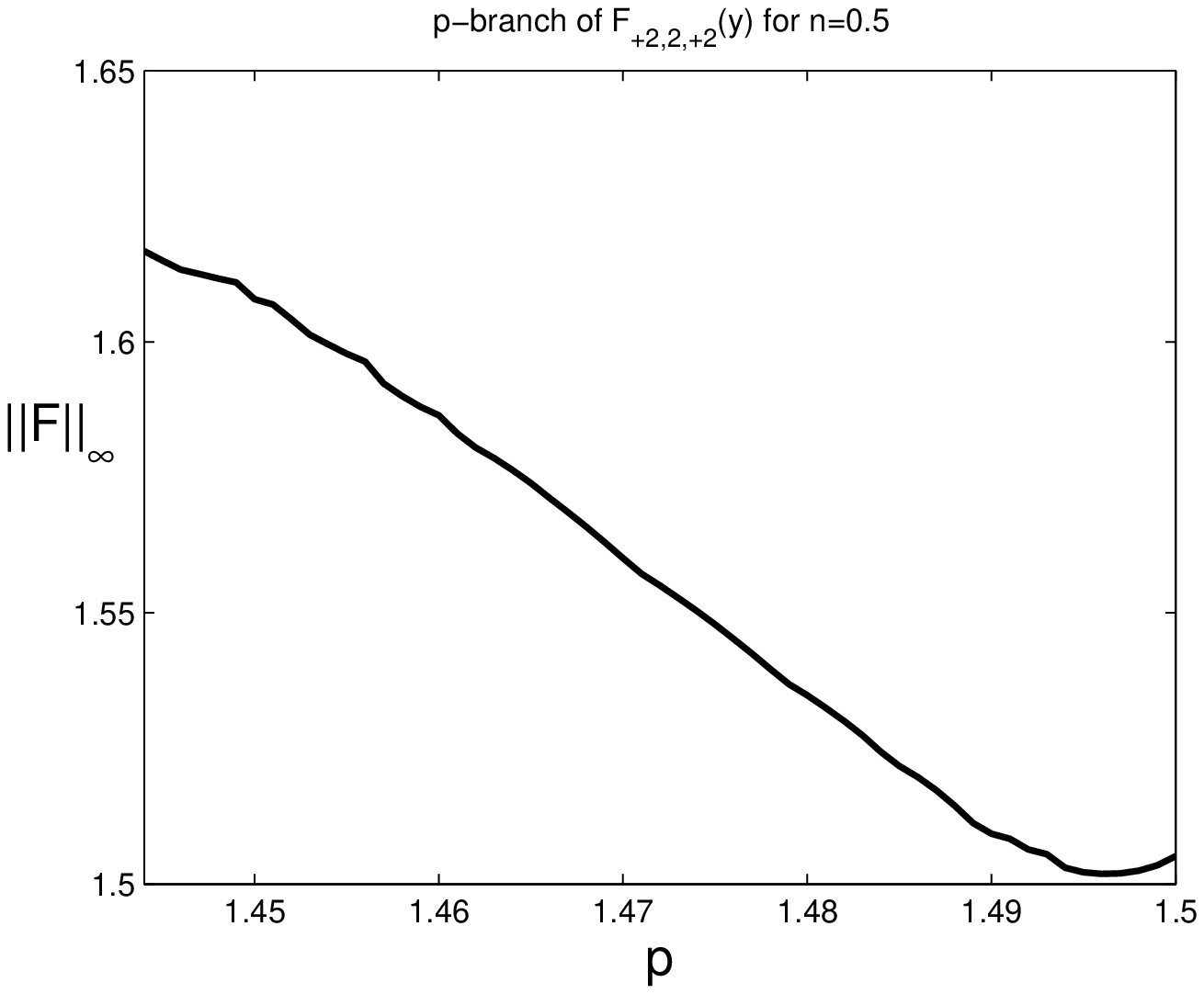}
} \centering \subfigure[$F_{+2,2,+2}$ profiles]{
\includegraphics[scale=0.52]{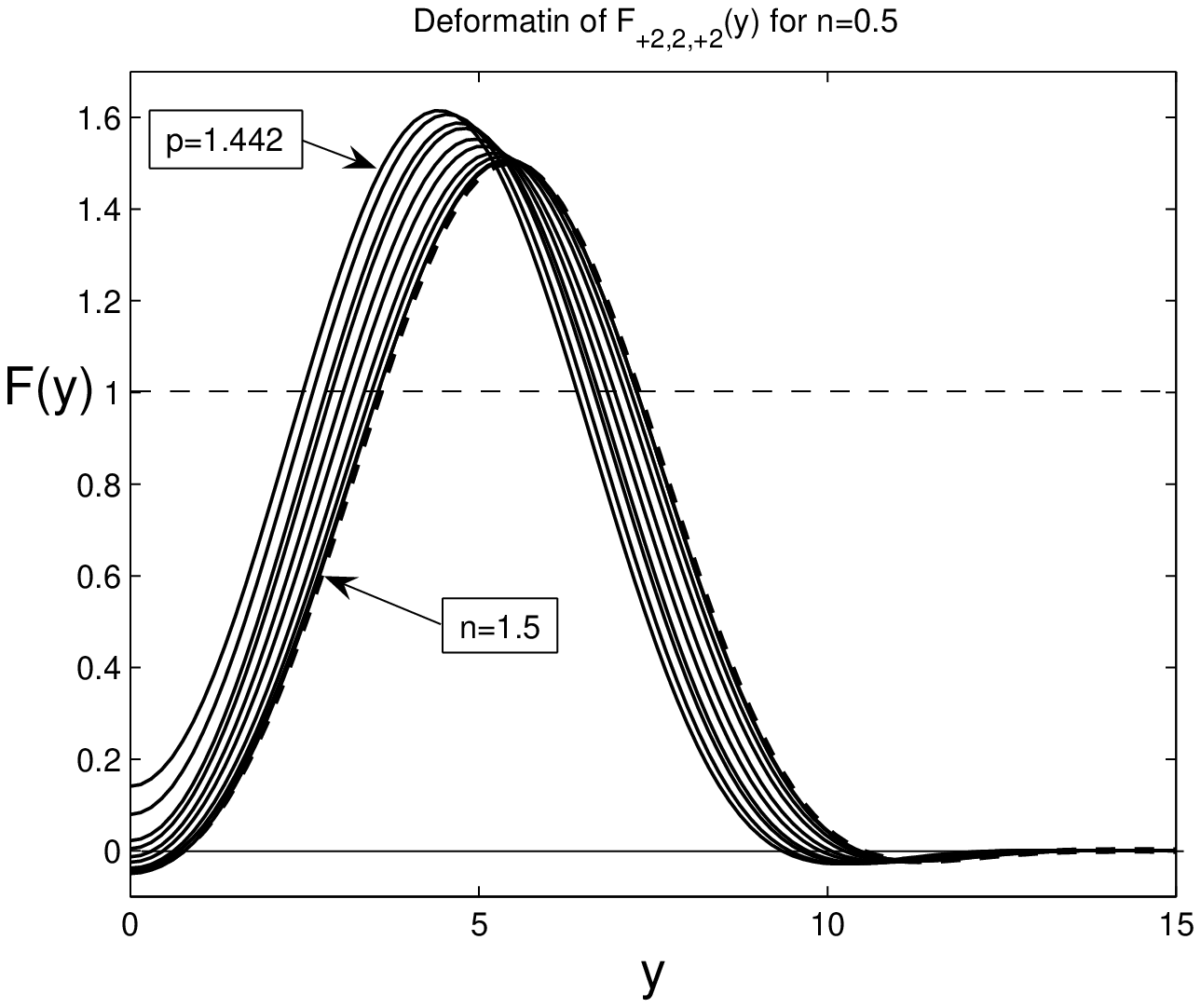}
}
 \vskip -.2cm
\caption{\rm\small  The  $p$-branch of solutions $ F_{+2,2+2}(y)$
of equation (\ref{2N}) for $n=0.5$ (a); corresponding deformation
of $F_{+2,2+2}$ (b).}
 \label{FG2}
\end{figure}

Concerning more complicated profiles not from the basic family
$\{F_l\}$, in Figures \ref{FG1} and \ref{FG2}, we show the
$p$-branches of the profiles $F_{+4}(y)$ and $F_{+2,2,+2}(y)$ for
the ODE (\ref{2N}) with $n=0.5$. We observe that both $p$-branches
(and the corresponding deformations) look very similar, and
actually, they do coincide for
 $$
 p \le p_* \approx 1.487.
 $$
We then claim that there exists a {\em branching point} $p=p_*
\in(1,n+1)$, at which the $p$-branch splits into two, thus giving
us two different profiles with similar geometric shapes; see
\cite{VainbergTr} for classic branching theory. Then this is not a
standard saddle-node bifurcation, and analytical (and numerical)
justification of such branching for the equivalent integral
equation with non-differentiable nonlinearities represents a
difficult open problem.
 Note that, according to (\ref{6}), higher-order profiles
 $F_{+2k}$ and geometrically neighbouring ones
 $F_{+2,2,...,+2}$ can also appear at bifurcations from the
  constant  equilibrium $F_*$. This needs special
 research and remains an open problem.

Finally, we must admit that our numerics did not detect the
$p$-extensions to the global blow-up  range $p <n+1$ of the odd
(anti-symmetric) basic profiles such as $F_1(y)$ ((b) in Figure
\ref{G4}), or $F_3(y)$ ((d) in Figure \ref{G4}). We always
observed a strong divergence of the method applied to the ODE
(\ref{2N}) without obvious reasons. Nevertheless, we expect that
these $p$-branches of basic blow-up profiles do exist. This formal
conclusion is associated with the question of evolution
completeness of blow-up solutions for the PDE (\ref{1.5}): indeed,
then which patterns will describe a generic evolution as $t \to
T^-$ for classes of anti-symmetric initial data? Of course, these
could be non-similarity solutions, but those hard questions are
definitely out of the present rather preliminary self-similar
consideration.


\ssk

{\bf Acknowledgement.} The  author would like to thank both the
anonymous  Reviewers for careful reading the manuscript and a
number of useful suggestions that have been used in the final
version.








\end{document}